\documentclass[a4paper,11pt]{scrartcl}
\usepackage{graphicx}
\usepackage{amsmath,amsthm,amssymb}
\usepackage{enumitem}
\usepackage{algpseudocode,algorithm}
\usepackage[square,sort&compress,numbers]{natbib}
\usepackage{type1cm}
\usepackage{url}
\usepackage{mathtools}\mathtoolsset{showonlyrefs=false}
\usepackage{booktabs}

\usepackage{subcaption}
\newcommand{\Authornote}{\renewcommand{\thefootnote}{\fnsymbol{footnote}}}
\newcommand{\authornote}{\Authornote\footnote}
\theoremstyle{plain}
\newtheorem{theorem}{Theorem}[section]

\theoremstyle{definition}

\theoremstyle{remark}
\newtheorem{remark}[theorem]{Remark}

\newcommand{\refalg}[1]{Algorithm~\ref{#1}}

\newcommand{\reffig}[1]{Figure~\ref{#1}}
\newcommand{\reftab}[1]{Table~\ref{#1}}

\newcommand{\finbox}{\nolinebreak\hfill{\small $\blacksquare$}}

\newcounter{alnum}

\newenvironment{Problem}{\begin{array}.{*{20}{l}}\}}{\end{array}}

\newcommand{\MIN}{\mathop{\mathrm{Minimize}}}

\newcommand{\ST}{\mathop{\mathrm{subject~to}}}

\newcommand{\diag}{\mathop{\mathrm{diag}}\nolimits}

\renewcommand{\Re}{\ensuremath{\mathbb{R}}}
\newcommand{\bi}[1]{\ensuremath{\boldsymbol{#1}}}

\newcommand{\rr}[1]{\ensuremath{\mathrm{#1}}}

\newcommand{\SC}{\ensuremath{\mathcal{S}}}

\begin{document}

\begin{center}
  {\Large\bfseries\sffamily% 
  Robust Truss Topology Optimization via Semidefinite }
  \par\medskip
  {\Large\bfseries\sffamily% 
  Programming with Complementarity Constraints: }
  \par\medskip
  {\Large\bfseries\sffamily% 
  A Difference-of-Convex Programming Approach 
  }%
  \par%
  \bigskip%
  {
  Yoshihiro Kanno~\authornote[2]{%
%  Corresponding author. 
    Mathematics and Informatics Center, 
    The University of Tokyo, 
    Hongo 7-3-1, Tokyo 113-8656, Japan.
    E-mail: \texttt{kanno@mi.tokyo-u.ac.jp}. 
  % Phone: +81-3-5841-6913. 
  % Fax: +81-3-5841-6886.
  }
  }
\end{center}

\begin{abstract}
  The robust truss topology optimization against the uncertain static 
  external load can be formulated as mixed-integer semidefinite 
  programming. 
  Although a global optimal solution can be computed with a 
  branch-and-bound method, it is very time-consuming. 
  This paper presents an alternative formulation, semidefinite 
  programming with complementarity constraints, and proposes an 
  efficient heuristic. 
  The proposed method is based upon the convex-concave procedure for DC 
  (difference-of-convex) programming. 
  It is shown that the method can often find a practically reasonable 
  truss design within the computational cost of solving some dozen of 
  convex optimization subproblems. 
\end{abstract}

\begin{quote}
  \textbf{Keywords}
  \par
  Robust optimization; 
  design-dependent load; 
  complementarity constraint; 
  semidefinite programming; 
  difference-of-convex programming; 
  concave-convex procedure. 
\end{quote}

\section{Introduction}

Many studies have been done on robust optimization of structures against 
uncertainty in external loads. 
A possibilistic (or bounded-but-unknown) model of uncertainty, 
assuming only the set of values that the input data can possibly take, 
might be useful when reliable statistical property of uncertainty, 
which is required for a probabilistic model of uncertainty, is 
unavailable or imprecise. 
With a possibilistic model of uncertainty, design optimization 
considering structural robustness against the uncertainty is treated 
within the framework of robust optimization \cite{BtEN09}. 

Attention of this paper is focused on robust topology optimization of 
truss structures against uncertainty in the static nodal external load.\footnote{
A truss is an assemblage of straight bars (called 
members) connected by pin-joints (called nodes) that do not transfer moment.
See section~\ref{sec:motivation} for some concrete examples. } 
Namely, we attempt to find a truss design that minimizes the worst-case 
compliance, i.e., the maximal value of the compliance among the 
specified set of external loads.\footnote{
The compliance of a truss, formally defined by 
\eqref{eq:def.compliance.0}, is equivalent to the twice strain energy of 
the truss at the equilibrium state under the prescribed boundary 
conditions. It can be regarded as a global measure of the 
displacements, and hence by minimizing the compliance the global 
stiffness of the truss is maximized.  } 
The seminal work of \citet{BtN97} shows that, based on the conventional 
ground structure method, this optimization problem can be formulated as 
{\em semidefinite programming\/} (SDP).\footnote{
The ground structure method is commonly used in truss 
topology optimization. 
It prepares an initial setting, called the ground structure, consisting 
of many members connected by nodes. 
The cross-sectional areas of the members are treated as design variables, 
while the locations of the nodes are specified. 
See section~\ref{sec:motivation} for more account. } 
SDP is a class of convex optimization, and can be solved efficiently 
with a primal-dual interior-point method \cite{AL12,WSV00}. 
Closely related formulations of continuum-based robust structural 
optimization can be found in 
\cite{CC03,CC08,TNKK11,BST12,CD08,dGAJ08,HK15,Tho16}. 
Furthermore, nonlinear SDP approaches to robust structural 
optimization have been proposed in \cite{GBZG09,GDG11,KT06,HTK15}. 
This paper discusses and deals with intrinsic difficulty in robust truss 
topology optimization. 

Suppose that uncertain external forces can possibly be 
applied to at any nodes, and that no external force is applied 
at an intermediate point of a truss member. 
If the set of exiting nodes is specified, then the robust truss 
optimization problem can be recast as SDP \cite{BtN97}. 
In this approach, all the nodes in the specified set remain in the 
obtained solution. 
Also, the obtained solution can possibly have some nodes other han the 
specified ones, but at such extra nodes no uncertain external force is 
considered. 
Thus, it is difficult to predict in advance the set of existing nodes in 
the robust optimal truss. 
In other words, the uncertainty model of external forces should be 
treated as a design-dependent model \cite{KG10}. 
This design dependency can be addressed by introducing 0-1 variables to 
represent the set of existing members in a truss design \cite{YK10}. 
The robust truss topology optimization problem is then formulated as 
{\em mixed-integer semidefinite programming problem\/} (MISDP), which 
can be solved globally with a branch-and-bound method \cite{YK10}. 
Unfortunately, due to large computational cost, this MISDP approach can 
be applied only to small-scale problem instances \cite{YK10}. 

Another issue that has not been considered in literature on robust truss 
topology optimization \cite{BtN97,YK10,KG10} is the treatment of 
parallel consecutive members in the ground structure method. 
Specifically, in robust truss topology optimization with uncertain 
external load, overlapping members in the ground structure are not redundant.\footnote{
With reference to concrete examples, we will thoroughly discuss this 
issue in section~\ref{sec:motivation}. } 
Such non-redundancy of overlapping members has also been recognized  
in truss topology optimization considering the self-weight load 
\cite{BBtZ94,KY17} and the member buckling constraints \cite{Mel14,GCO05}. 
Consider the conventional truss topology optimization. 
It is often that the optimal solution has parallel consecutive members 
that are connected by nodes  supported only in the direction of those 
members. 
A sequence of such members is sometimes called a {\em chain\/} \cite{Ach99}. 
If only the compliance is considered as the structural performance, one 
can remove the intermediate nodes and replace the chain with a single 
longer member. 
This procedure is called the {\em hinge cancellation\/} \cite{Ach99,Roz96}. 
Since the hinge cancellation does not change the compliance, overlapping 
of members in a ground structure can be removed in advance by deleting 
the longer member when two members overlap. 
In contrast, in robust truss optimization under load uncertainties, a 
solution having a chain is infeasible,\footnote{
A solution having a chain cannot be in equilibrium with uncertain loads 
applied at intermediate nodes of the chain. 
Therefore, the worst-case compliance of the solution is infinitely large. 
} 
while the one stabilized by hinge cancellation can be feasible. 
This means that, in 
general, a global optimal truss topology cannot be captured  without 
incorporating overlapping members in a ground structure. 
However, on the other hand, presence of overlapping members in a final 
truss design is not allowed from a practical point of view. 
Therefore, a special treatment is required in a robust optimization 
method to prohibit the presence of overlapping members. 

This paper addresses the two difficulties in robust truss optimization 
explained above: the design dependency of the uncertainty model of the 
external load and the necessity of incorporating overlapping members in 
a ground structure. 
Both can be dealt with by introducing, for each member, a 0-1 design 
variable indicating whether the member vanishes or exists. 
Therefore, the robust truss topology optimization can be formulated as 
MISDP; see section~\ref{sec:preliminary.node}. 
However, as mentioned before, this approach is applicable only to 
problems of small size. 
In contrast, in this paper 
we attempt to propose a heuristic that can often find 
a feasible solution with a reasonable objective value. 
Through the numerical experiments with 
problem instances having up to about $700$ members, 
it is shown that the proposed method usually converges after solving 
only a few dozen of convex optimization subproblems, and that it finds 
a practically reasonable solution. 

This paper is partially inspired by papers of 
\citet{JPW16} and \citet{LB16}. 
\citet{JPW16} present a DC (difference-of-convex) programming approach 
to finding a stationary point of linear programming with complementarity 
constraints; see also \cite{LtPd11,MDLtT12} for DC programming 
approaches to complementarity constraints. 
A function is said to be a DC (difference-of-convex) function if it can 
be represented as a difference of two convex functions. 
A DC programming problem is a minimization problem of a DC function 
under some inequality constraints, where all the constraint functions 
are DC functions. 
One of well-known local heuristics for finding a local optimal solution 
of DC programming is the concave-convex procedure\footnote{%
The concave-convex procedure is also known as 
the convex-concave procedure \cite{LB16,YCL14}. } \cite{SL09,CSWB06,FM01,NSS05}. 
\citet{LB16} show that an extension of the concave-convex procedure 
can serve as an efficient heuristic for diverse nonconvex optimization 
problems. 
For more account on the DC programming and the concave-convex procedure, 
see section~\ref{sec:preliminary.DC} and the references therein. 
In this paper, we first formulate the robust truss topology optimization 
as {\em semidefinite programming with complementarity constraints\/} (SDPCC). 
Following an idea found in \cite{JPW16}, we recast this problem as a DC 
programming problem. 
A variant of the concave-convex procedure, which is similar to the one 
in \cite{LB16}, is then applied to this DC programming formulation. 
Each iteration of the proposed method consists of solving an SDP problem. 

The paper is organized as follows. 
In section~\ref{sec:motivation} we explain intrinsic difficulties in 
robust truss topology optimization by using some illustrative examples. 
Section~\ref{sec:algorithm} provides an overview of the necessary 
background of the DC programming and the concave-convex procedure, and 
presents the general framework of the algorithm used in this paper. 
Section~\ref{sec:integer} briefly reviews the existing MISDP formulation 
for robust truss topology optimization, and extends it to the problem 
setting with a ground structure incorporating overlapping members. 
Section~\ref{sec:heuristic} presents a new formulation and solution 
method for robust truss topology optimization. 
Section~\ref{sec:ex} reports the results of numerical experiments. 
Conclusions are drawn in section~\ref{sec:conclude}. 

In our notation, 
we use $\bi{x}^{\top}$ and $X^{\top}$ to denote the transposes of vector 
$\bi{x} \in \Re^{n}$ and matrix $X \in \Re^{m \times n}$, respectively. 
For vectors $\bi{x} = (x_{i}) \in \Re^{n}$ and 
$\bi{y} = (y_{i}) \in \Re^{n}$, we write $\bi{x} \ge \bi{y}$ 
if $x_{i} \ge y_{i}$ $(i=1,\dots,n)$. 
Particularly, $\bi{x} \ge \bi{0}$ means 
$x_{i} \ge 0$ $(i=1,\dots,n)$. 
The Euclidean norm of $\bi{x}$ is denoted by 
$\| \bi{x} \| = \sqrt{\bi{x}^{\top} \bi{x}}$. 
We use $\bi{1} = (1,1,\dots,1)^{\top}$ to denote the all-ones vector. 
Let $\SC^{n}$ denote the set of $n \times n$ real symmetric matrices. 
We write $X \succeq O$ if $X \in \SC^{n}$ is positive semidefinite. 
We use $\diag(\bi{x})$ to denote a diagonal matrix, the vector of 
diagonal components of which is $\bi{x}$. 
For a finite set $T$, we use $|T|$ to denote the cardinality of 
$T$, i.e., the number of elements in $T$.

\section{Motivation}
\label{sec:motivation}

In this section, we explain intrinsic difficulties in robust truss 
topology optimization, which motivate us to develop the method proposed 
in this paper. 
Details of the examples in this section appear in 
section~\ref{sec:ex.small}. 

Suppose that uncertain static external forces are applied at all the 
nodes, and only at nodes, of a truss. 
The robust truss topology optimization is to find a truss design that 
minimizes the worst-case compliance, i.e., the maximum value of the 
compliance among possible realizations of the external load, under the 
upper bound constraint on the structural volume. 

\begin{figure}[tp]
  %%%% C:\doc\robust\topology\load_ccp\x2y1_for_fig\opt_design.m
  \centering
  \begin{subfigure}[b]{0.32\textwidth}
    \centering
    \includegraphics[scale=0.8]{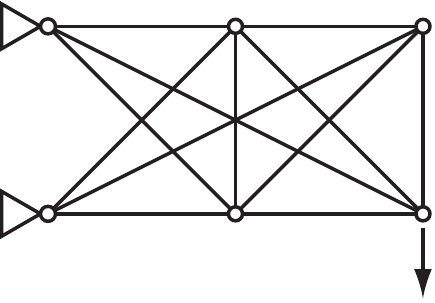}
    \caption{}
    \label{fig:gs2x1}
  \end{subfigure}
  \hfill
  \begin{subfigure}[b]{0.32\textwidth}
    \centering
    \includegraphics[scale=0.30]{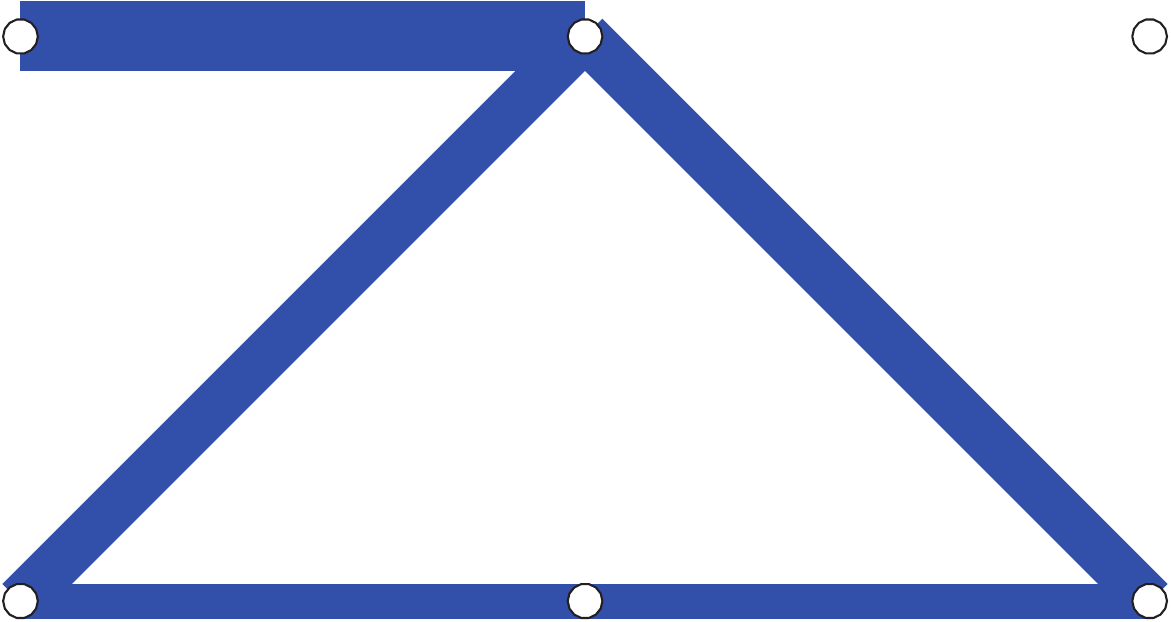}
    \caption{}
    \label{fig:x2_y1_nominal}
  \end{subfigure}
  \hfill
  \begin{subfigure}[b]{0.32\textwidth}
    \centering
    \includegraphics[scale=0.30]{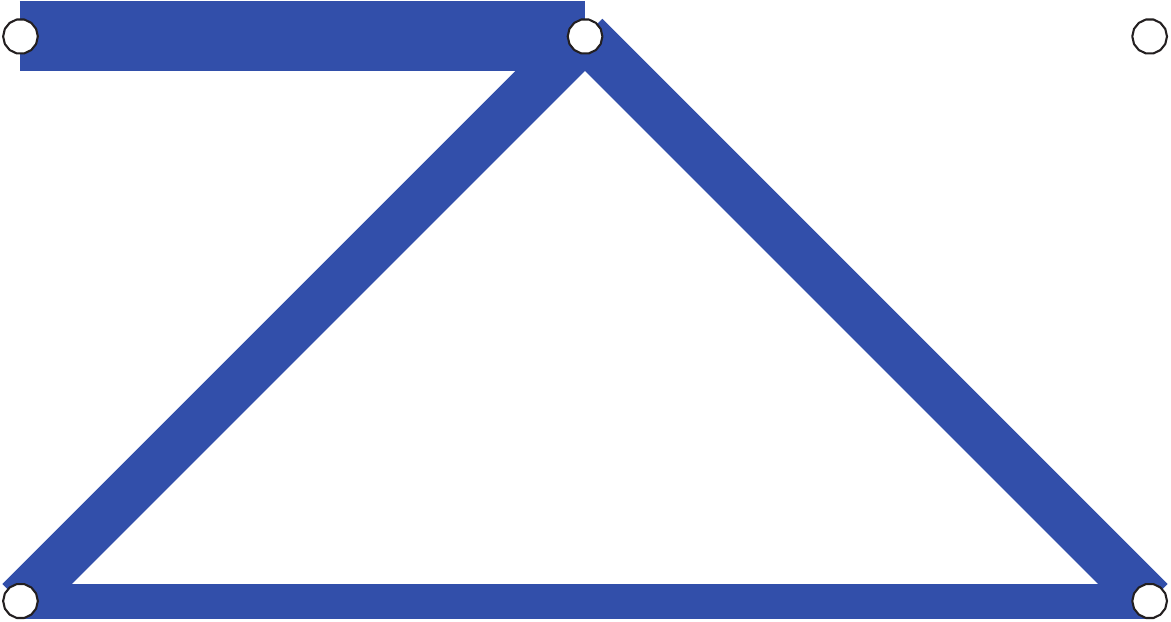}
    \caption{}
    \label{fig:x2_y1_nominal_reduced}
  \end{subfigure}
  \caption{An example of the conventional compliance minimization. 
  \subref{fig:gs2x1} The ground structure (with $12$ members); 
  \subref{fig:x2_y1_nominal} the optimal solution; and 
  \subref{fig:x2_y1_nominal_reduced} the final truss design after hinge 
  cancellation. 
  }
  \label{fig:x2_y1}
\end{figure}

\begin{figure}[tp]
  %%%% C:\doc\robust\topology\load_ccp\x2y1_for_fig\opt_design.m
  \centering
  \begin{subfigure}[b]{0.32\textwidth}
    \centering
    \includegraphics[scale=0.30]{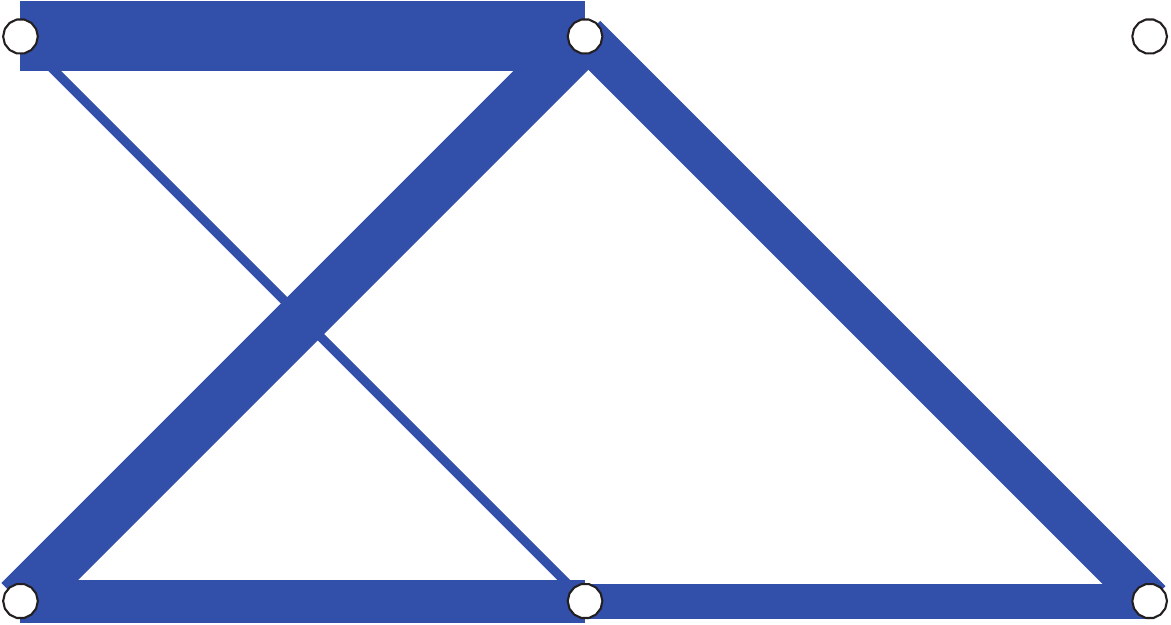}
    \caption{}
    \label{fig:x2_y1_specified}
  \end{subfigure}
  \hfill
  \begin{subfigure}[b]{0.32\textwidth}
    \centering
    \includegraphics[scale=0.8]{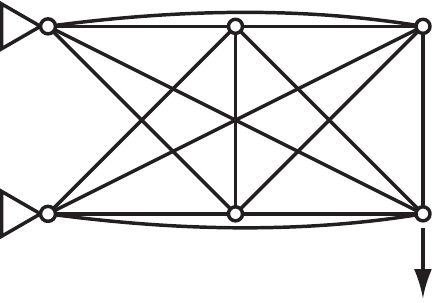}
    \caption{}
    \label{fig:gs2x1robust}
  \end{subfigure}
  \hfill
  \begin{subfigure}[b]{0.32\textwidth}
    \centering
    \includegraphics[scale=0.30]{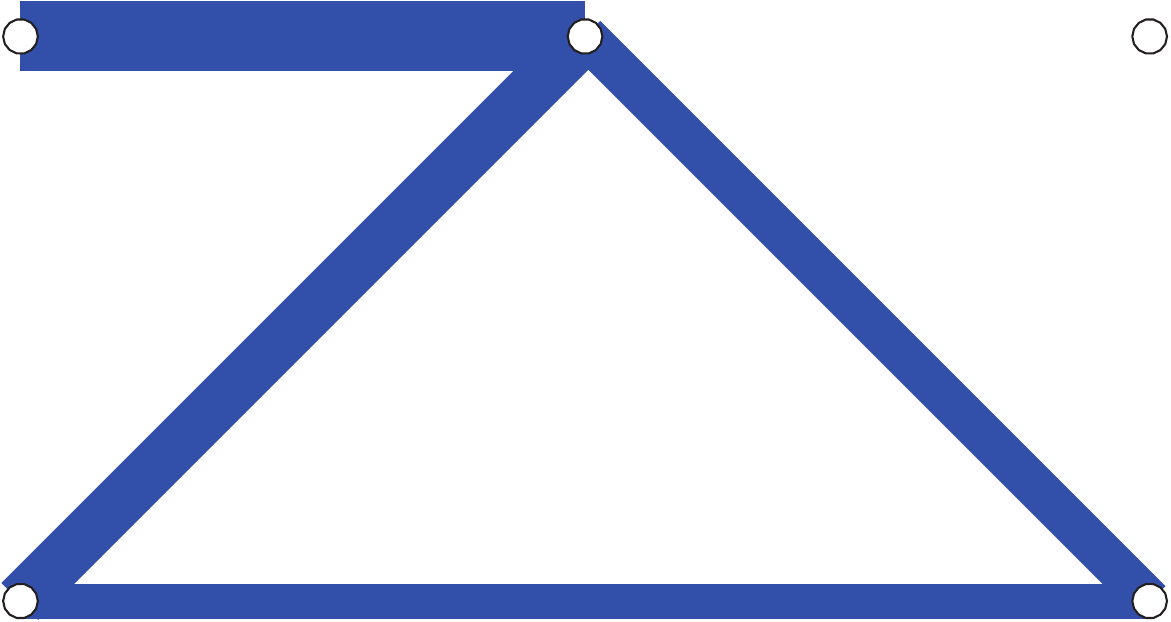}
    \caption{}
    \label{fig:x2_y1_global}
  \end{subfigure}
  \caption{The robust optimization corresponding to the example in 
  \reffig{fig:x2_y1}. 
  \subref{fig:x2_y1_specified} The optimal solution obtained from 
  the ground structure in \reffig{fig:gs2x1}; 
  \subref{fig:gs2x1robust} the ground structure including overlapping 
  members ($14$ members in total); and 
  \subref{fig:x2_y1_global} the global optimal solution. 
  }
  \label{fig:robust.x2_y1}
\end{figure}

With reference to the examples in \reffig{fig:x2_y1} and 
\reffig{fig:robust.x2_y1}, we explain that incorporating overlapping 
members into a ground structure is necessary for the robust truss 
topology optimization. 
We begin with the conventional compliance minimization, without 
considering uncertainties. 
\reffig{fig:gs2x1} shows a ground structure, which consists of $12$ 
members and has no overlapping members. 
The ground structure is an initial setting for 
truss topology optimization. 
The members consisting of a ground structure are called the candidate 
members, and their cross-sectional areas are design variables to be 
optimized. 
It is worth noting that the locations of the nodes are not treated as 
the design variables. 
If the cross-sectional area of a member becomes equal to zero as a result 
of optimization, then the member is removed from the truss. 
Thus, the connectivity of members, called the topology in this research 
area, usually changes from the ground structure. 
A vertical external load is applied at the rightmost bottom node, as 
shown in \reffig{fig:gs2x1}. 
The optimal solution is shown in \reffig{fig:x2_y1_nominal}, where 
the width of each member is proportional to its cross-sectional area. 
This solution has a chain consisting of two members. 
The hinge cancellation yields the final truss design 
shown in \reffig{fig:x2_y1_nominal_reduced}. 
It should be clear that the objective value, i.e., the compliance, of 
the solution in \reffig{fig:x2_y1_nominal_reduced} is same as the one in 
\reffig{fig:x2_y1_nominal}. 
We next consider the robust optimization. 
Since uncertain external force is supposed to be applied also at the 
intermediate node of the chain, the solution in 
\reffig{fig:x2_y1_nominal} becomes infeasible. 
As a result, the optimal solution has one additional member to stabilize 
that node, as shown in \reffig{fig:x2_y1_specified}. 
Alternatively, consider a ground structure in \reffig{fig:gs2x1robust}, 
which consists of $14$ members. 
The newly added two members, depicted as slightly 
curved lines, are in fact straight bars. 
Therefore, each of them overlaps with two shorter members, and is 
called an overlapping members. 
Moreover, the middle node is exactly located on a line forming an 
overlapping member. 
Therefore, we say that the overlapping member lies on the middle node. 
The robust optimal solution obtained from this ground structure is shown 
in \reffig{fig:x2_y1_global}.\footnote{
It should be clear that the worst-case compliance for the solution in 
\reffig{fig:x2_y1_global} is smaller than that for the solution in 
\reffig{fig:x2_y1_specified}. 
} 
Namely, at the global optimal solution, the longer member is selected 
instead of the chain and the intermediate node of the chain is removed. 
Thus, overlapping members do not mean redundancy, because connection 
of two members via an intermediate node introduces an extra uncertain 
load applied to that node. 
Similar non-redundancy of overlapping members in a ground structure can 
be observed also in, e.g., the compliance minimization of trusses under 
the self-weight loads \cite{BBtZ94,KY17}. 

A key in this robust optimization is selecting the set of nodes which 
the optimal solution has. 
This is because the uncertainty in external loads depends on the set of 
existing nodes, in a manner that 
uncertain external forces are supposed to be 
applied to all the existing nodes. 
Moreover, a node lying on an existing member should be removed. 
The next example illustrates that selecting a optimal set of nodes in 
a heuristic manner is indeed difficult. 

\begin{figure}[tp]
  %%%% C:\doc\robust\topology\load_ccp\x3y3_for_fig\opt_design.m
  \centering
  \begin{subfigure}[b]{0.32\textwidth}
    \centering
    \includegraphics[scale=0.30]{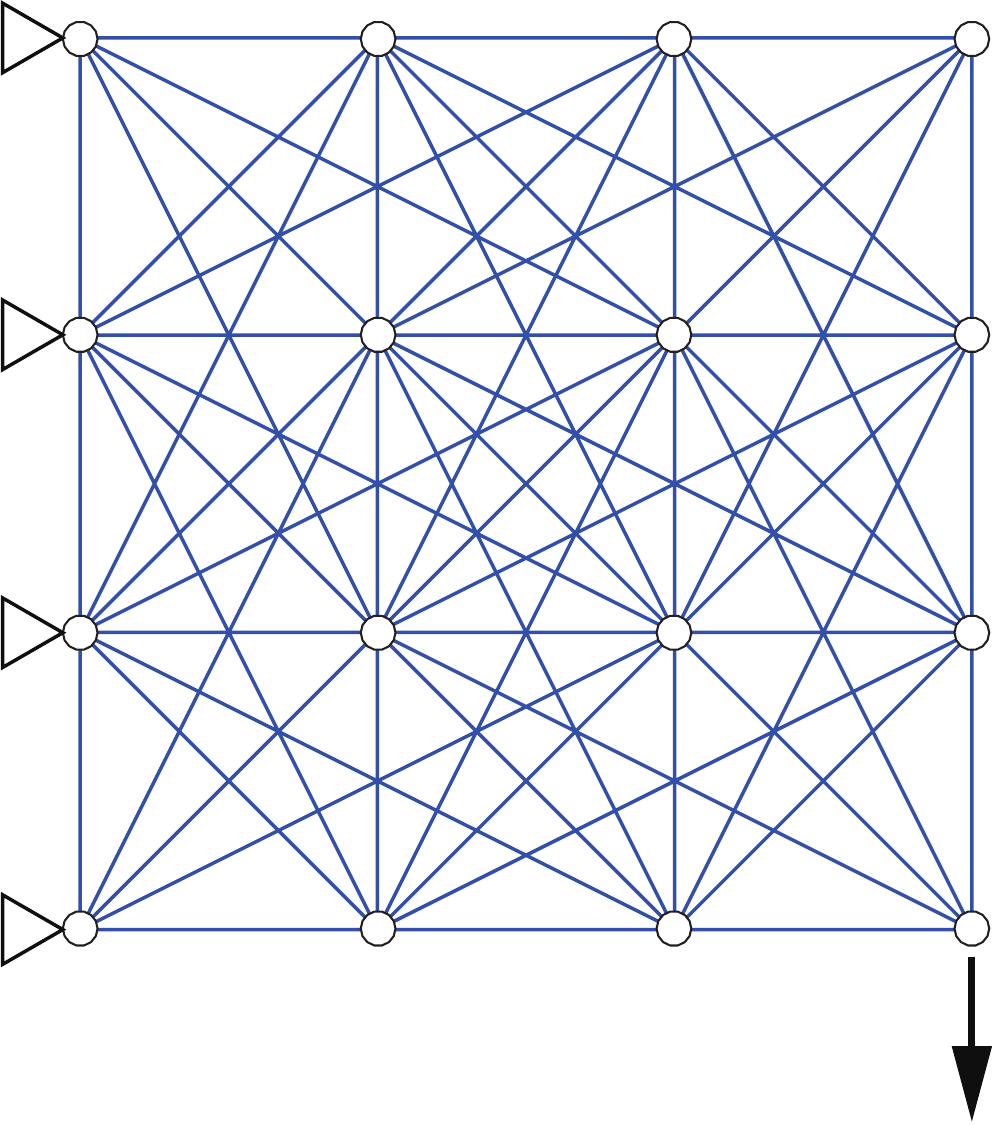}
    \caption{}
    \label{fig:gs_x3_y3}
  \end{subfigure}
  \hfill
  \begin{subfigure}[b]{0.32\textwidth}
    \centering
    \includegraphics[scale=0.35]{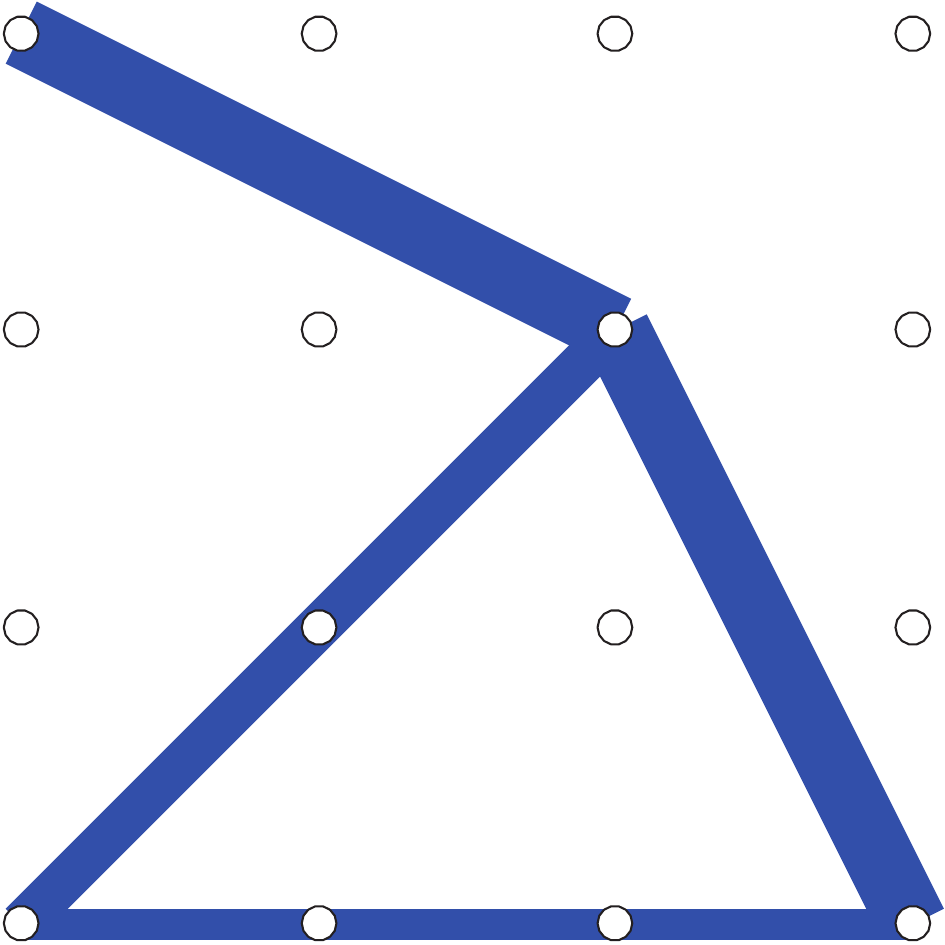}
    \caption{}
    \label{fig:x3_y3_nominal}
  \end{subfigure}
  \par\medskip
  \begin{subfigure}[b]{0.32\textwidth}
    \centering
    \includegraphics[scale=0.35]{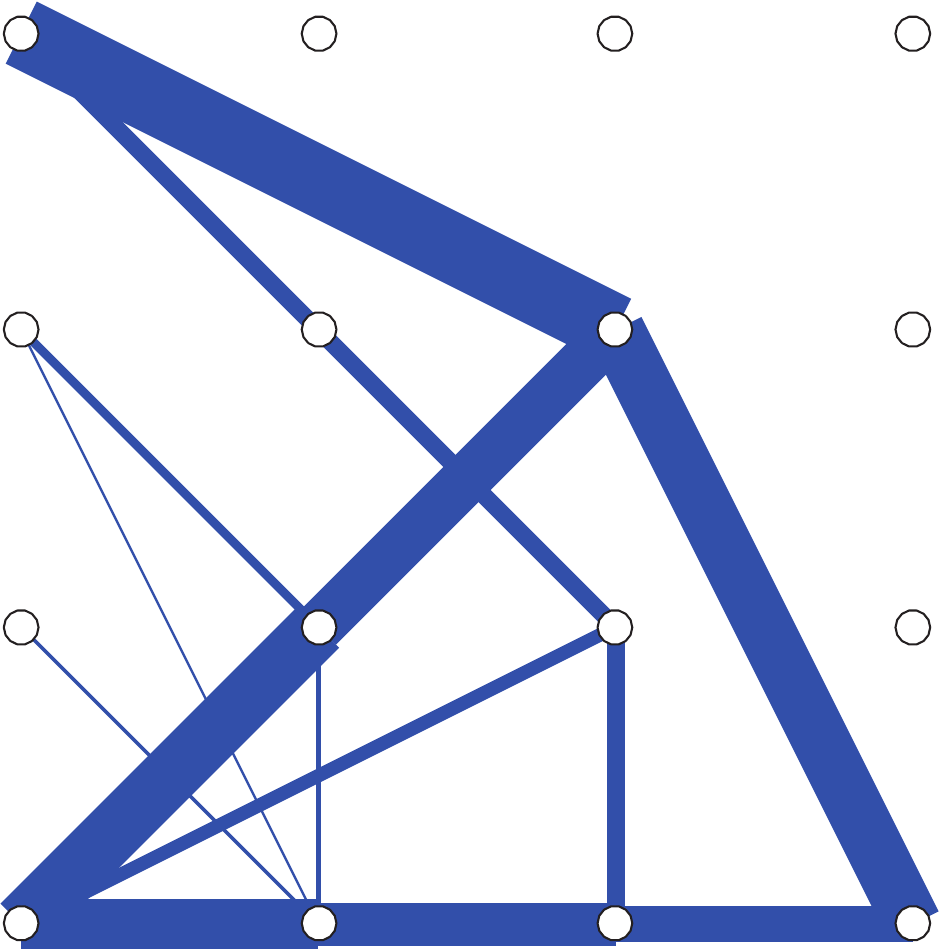}
    \caption{}
    \label{fig:x3_y3_specified1}
  \end{subfigure}
  \hfill
  \begin{subfigure}[b]{0.32\textwidth}
    \centering
    \includegraphics[scale=0.35]{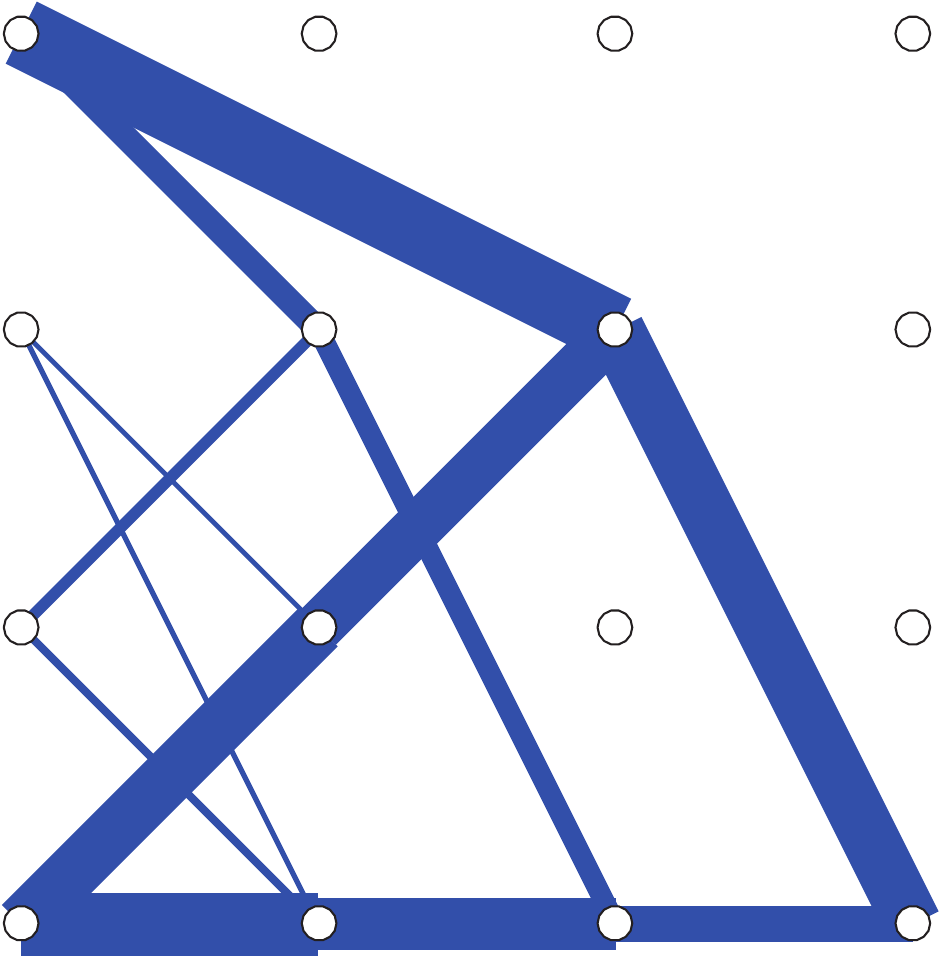}
    \caption{}
    \label{fig:x3_y3_global_nov}
  \end{subfigure}
  \hfill
  \begin{subfigure}[b]{0.32\textwidth}
    \centering
    \includegraphics[scale=0.35]{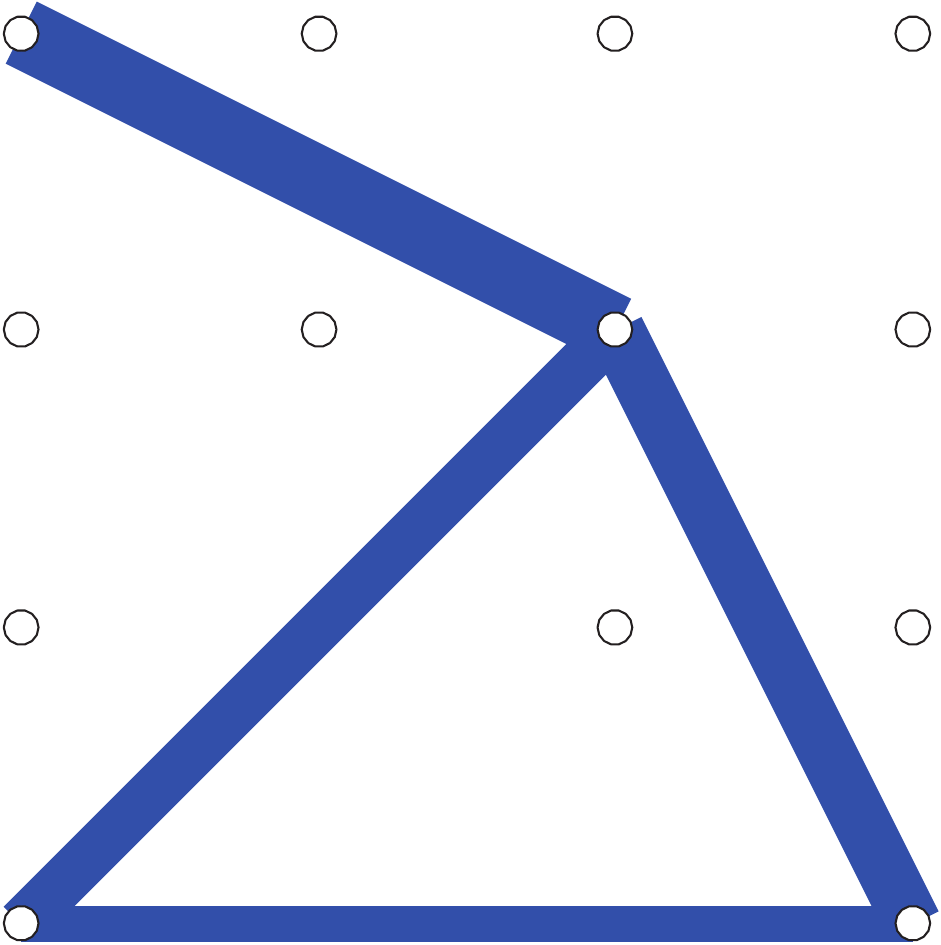}
    \caption{}
    \label{fig:x3_y3_global}
  \end{subfigure}
  \caption{Difficulties in the robust truss topology optimization. 
  \subref{fig:gs_x3_y3} The ground structure and the nominal external load; 
  \subref{fig:x3_y3_nominal} the optimal solution for the nominal load; 
  \subref{fig:x3_y3_specified1} a solution obtained for a 
  design-independent uncertainty model of the external load; 
  \subref{fig:x3_y3_global_nov} the robust optimal solution for the 
  design-dependent uncertainty model and the ground structure without 
  overlapping members (the objective value is $3259.115\,\mathrm{J}$); and 
  \subref{fig:x3_y3_global} the robust optimal solution for the 
  design-dependent uncertainty model and the ground structure with 
  overlapping members (the objective value is $2442.708\,\mathrm{J}$). 
  }
  \label{fig:robust.x3_y3}
\end{figure}

Consider the problem setting shown in \reffig{fig:gs_x3_y3}. 
\reffig{fig:x3_y3_nominal} shows the optimal solution of the nominal 
(i.e., not robust) optimization problem. 
A simple heuristic to predict a set of existing nodes in a robust 
optimal solution is to adopt the set of nodes that the nominal optimal 
solution has. 
Suppose that uncertain external forces are applied {\em only\/} to the 
five free nodes that the solution in \reffig{fig:x3_y3_nominal} has. 
The optimal solution of this robust optimization problem is shown in 
\reffig{fig:x3_y3_specified1}. 
This solution has two extra nodes that the solution in 
\reffig{fig:x3_y3_nominal} does not have. 
Therefore, the solution in \reffig{fig:x3_y3_specified1} assumes that 
external forces are {\em not\/} applied to these two nodes (in this 
sense, this solution is {\em not truly robust\/}). 
On the other hand, \reffig{fig:x3_y3_global_nov} shows the optimal 
solution of the robust topology optimization with a ground structure 
which does {\em not\/} include overlapping members. 
It is observed that one of the nodes in \reffig{fig:x3_y3_specified1} is 
missing in the solution in \reffig{fig:x3_y3_global_nov}. 
Furthermore, \reffig{fig:x3_y3_global} shows the optimal solution of the 
robust topology optimization with a ground structure {\em including\/} 
overlapping members. 
It is observed that the three intermediate nodes in 
\reffig{fig:x3_y3_nominal} are removed and the two chains are replaced by 
longer members. 
As a result, the objective value of the solution in 
\reffig{fig:x3_y3_global_nov} is more than $1.33$ times larger than that of 
the solution in \reffig{fig:x3_y3_global}. 
Thus, it is crucial to design an optimization process so as to allow 
vanishment of intermediate nodes on chains. 

\begin{figure}[tp]
  %%%% C:\doc\robust\topology\load_ccp\x3y2_for_fig\opt_design.m
  \centering
  \begin{subfigure}[b]{0.45\textwidth}
    \centering
    \includegraphics[scale=0.35]{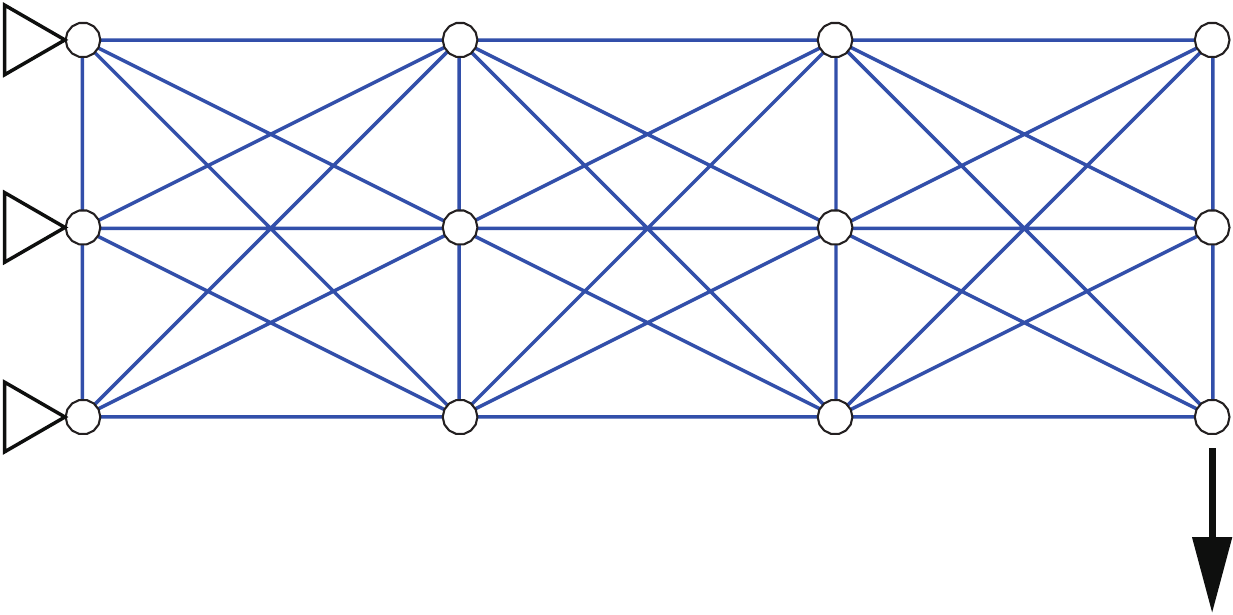}
    \caption{}
    \label{fig:gs_x3_y2}
  \end{subfigure}
  \hfill
  \begin{subfigure}[b]{0.45\textwidth}
    \centering
    \includegraphics[scale=0.45]{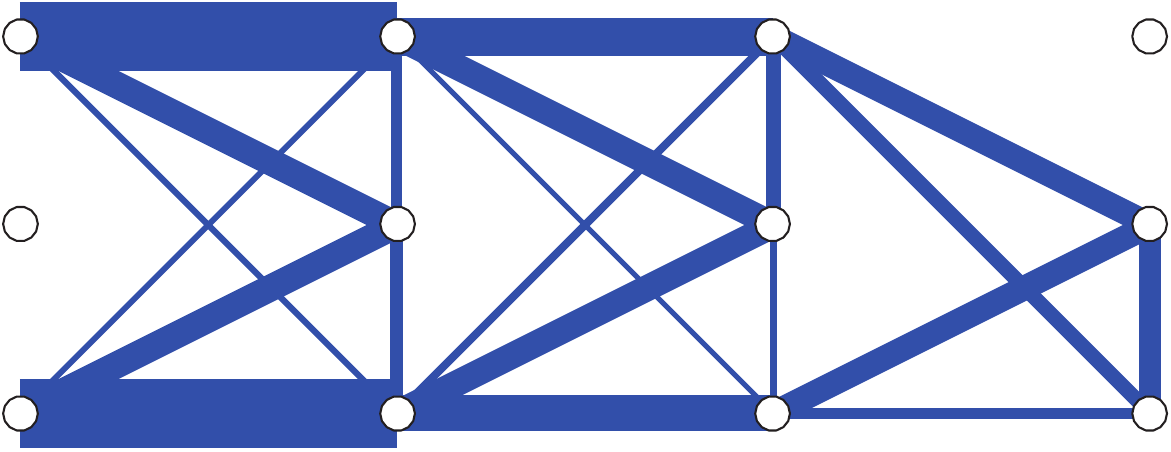}
    \caption{}
    \label{fig:x3_y2_nominal}
  \end{subfigure}
  \par\medskip
  \begin{subfigure}[b]{0.45\textwidth}
    \centering
    \includegraphics[scale=0.45]{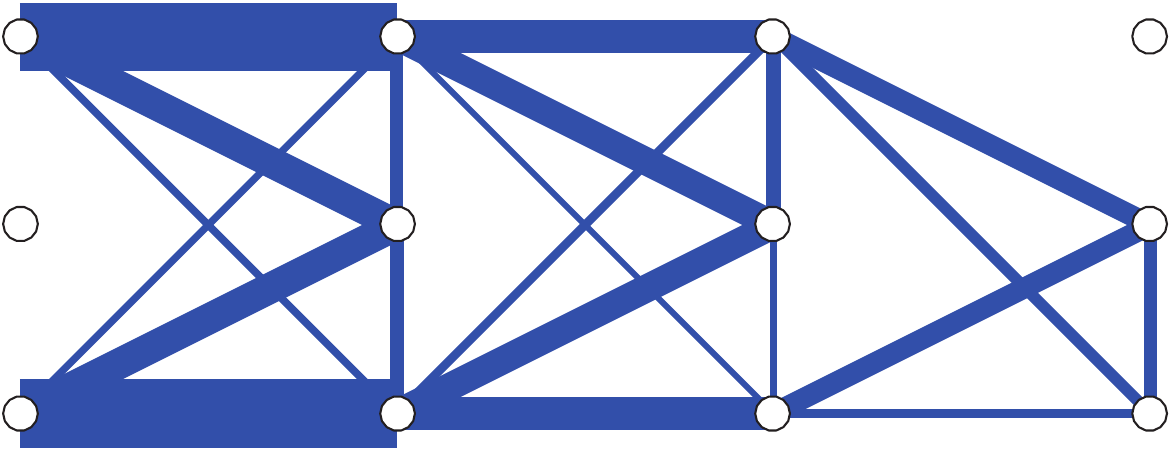}
    \caption{}
    \label{fig:x3_y2_specified}
  \end{subfigure}
  \hfill
  \begin{subfigure}[b]{0.45\textwidth}
    \centering
    \includegraphics[scale=0.45]{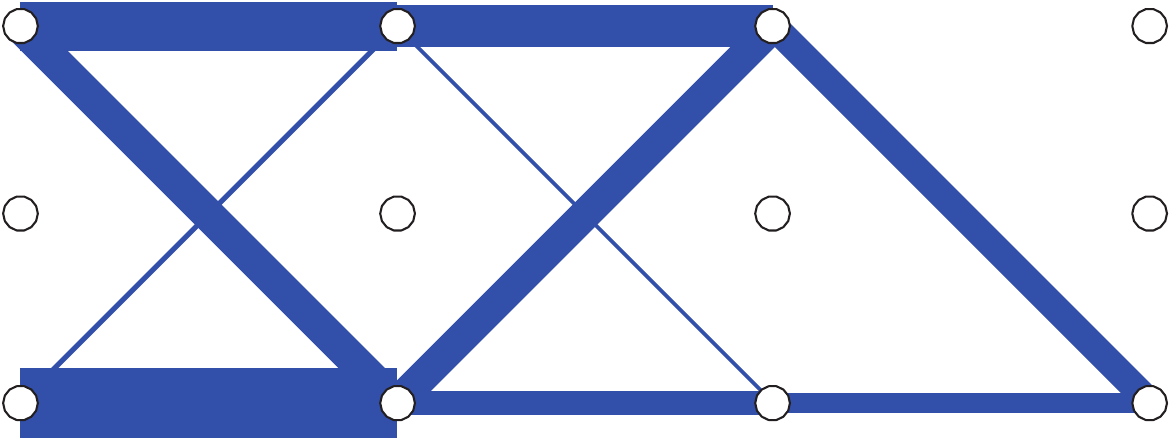}
    \caption{}
    \label{fig:x3_y2_global}
  \end{subfigure}
  \caption{A node which is not on a chain can also vanish as a result of 
  the robust topology optimization. 
  \subref{fig:gs_x3_y2} The ground structure and the nominal external load; 
  \subref{fig:x3_y2_nominal} the optimal solution for the nominal load; 
  \subref{fig:x3_y2_specified} the robust solution obtained for a 
  design-independent uncertainty model of the external load 
  (the objective value is $13934.896\,\mathrm{J}$); 
  \subref{fig:x3_y2_global} the robust optimal solution for the 
  design-dependent uncertainty model 
  (the objective value is $11093.750\,\mathrm{J}$). 
  }
  \label{fig:robust.x3_y2}
\end{figure}

It is also possible that a node which is {\em not\/} lying on a chain 
vanishes as a result of robust optimization. 
Consider the problem setting in \reffig{fig:gs_x3_y2}. 
The optimal solution of the nominal optimization problem is shown in 
\reffig{fig:x3_y2_nominal}. 
If we suppose that uncertain external forces are applied to the eight 
free nodes in \reffig{fig:x3_y2_nominal}, then the solution in 
\reffig{fig:x3_y2_specified} becomes optimal. 
Unlike the example in \reffig{fig:x3_y3_specified1}, uncertain external 
forces are supposed to be applied to all the nodes that the solution in 
\reffig{fig:x3_y2_specified} has. 
In this sense, the solution in \reffig{fig:x3_y2_specified} is a local 
optimal solution of the robust topology optimization. 
In contrast, the global optimal solution of the robust topology 
optimization is shown in \reffig{fig:x3_y2_global}.\footnote{
In this example, overlapping longer members are {\em not\/} incorporated 
into the ground structure, because with overlapping members 
the global optimization method (YALMIP~\cite{Lof04}) did not converge 
within realistic computational time. }
It is observed that three nodes in \reffig{fig:x3_y2_specified} are 
missing in \reffig{fig:x3_y2_global}. 
The objective value of the solution in \reffig{fig:x3_y2_specified} is 
more than $1.25$ times larger than that of the solution in 
\reffig{fig:x3_y2_global}. 
Thus, it is crucial to design an optimization algorithm that can deal 
with the design-dependent uncertainty model of the external load. 

In section~\ref{sec:integer} and section~\ref{sec:heuristic}, we 
propose a formulation and an algorithm for overcoming the difficulties 
discussed in this section.

\section{Algorithmic framework}
\label{sec:algorithm}

As preliminaries, section~\ref{sec:preliminary.DC} briefly introduces 
the notion of DC programming and the concave-convex procedure for 
solving it. 
In section~\ref{sec:algorithm.heuristic}, we introduce the optimization 
problem that we consider in this paper, and present an extension of the 
concave-convex procedure. 

\subsection{Fundamentals: DC programming and concave-convex procedure}
\label{sec:preliminary.DC}

Let $f_{i}$, $g_{i} : \Re^{n} \to \Re$ $(i=0,1,\dots,m)$ be convex. 
The optimization problem having the following form is called a DC 
programming problem: 
\begin{subequations}\label{P.DC}%
  \begin{alignat}{3}
    & \MIN
    &{\quad}& 
    f_{0}(\bi{x}) - g_{0}(\bi{x}) \\
    & \ST &&
    f_{i}(\bi{x}) - g_{i}(\bi{x}) \le 0 , 
    \quad i=1,\dots,m . 
  \end{alignat}
\end{subequations}
For simplicity, we assume that $g_{0}$, $g_{1},\dots,g_{m}$ are 
differentiable. 

The concave-convex procedure is known as a heuristic for finding a local 
optimal solution of problem \eqref{P.DC}. 
Let $\bi{x}^{(k)} \in \Re^{n}$ denote the (feasible) incumbent value of 
$\bi{x}$ at the $k$th iteration. 
Define $\hat{g}_{i}(\,\cdot\, ; \bi{x}^{(k)}) : \Re^{n} \to \Re$ 
$(i=0,1,\dots,m)$ by 
\begin{align}
  \hat{g}_{i}(\bi{x}; \bi{x}^{(k)})
  = g_{i}(\bi{x}^{(k)}) + \nabla g_{i}(\bi{x}^{(k)})^{\top} (\bi{x} - \bi{x}^{(k)})  . 
  \label{eq.def.hat.g.linear}
\end{align}
The concave-convex procedure updates the solution by letting 
$\bi{x}^{(k+1)}$ be the optimal solution of the following optimization problem: 
\begin{subequations}\label{P.CCP.1}%
  \begin{alignat}{3}
    & \MIN  &{\quad}& 
    f_{0}(\bi{x}) - \hat{g}_{0}(\bi{x}; \bi{x}^{(k)}) \\
    & \ST && 
    f_{i}(\bi{x}) - \hat{g}_{i}(\bi{x}; \bi{x}^{(k)}) \le 0 , 
    \quad i=1,\dots,m . 
  \end{alignat}
\end{subequations}
It is worth noting that this subproblem is convex. 

For the sequence, $\{ \bi{x}^{(k)} \}$, generated by the concave-convex 
procedure, it is known that the objective value of \eqref{P.DC}, i.e., 
$\{ f_{0}(\bi{x}^{(k)}) - g_{0}(\bi{x}^{(k)}) \}$, converges. 
However, $\{ \bi{x}^{(k)} \}$ does not necessarily converges to a local 
optimal solution; see, e.g., \cite[\S~1.3]{LB16}. 
Applications of the concave-convex procedure include 
transductive support vector machines (SVMs) \cite{CSWB06,FM01}, 
feature selection in SVMs \cite{NSS05}, etc. 

The concave-convex procedure can be considered as a version of 
DCA (difference of convex algorithm) \cite{PdLt97,PdLt98}; see 
\cite{LB16,SL09} for accounts of this fact. 
For DCA and its applications we direct the reader to 
\cite{LtPd05,PdLt14}; an application of DCA in structural engineering 
can be found in \cite{SP96}. 
As shown in \cite{SL09},  the concave-convex procedure can also be 
viewed as a variant of  MM algorithms (majorization-minimization algorithms) 
\cite{HL04,LCZ14}.\footnote{In fact, 
$\hat{g}_{i}(\cdot \,; \bi{x}^{(k)})$ defined by 
\eqref{eq.def.hat.g.linear} is a majorization function of $g_{i}$. }
The MM algorithm is a generalization of the well-known 
EM algorithm (expectation-maximization algorithm) \cite{DLR77}. 
The MM algorithms have been frequently employed in machine learning and 
image processing as seen in, e.g, \cite{HL05,HL00,FBdN07,STL11,SBP17}. 

MMA (the method of moving asymptotes) \cite{Sva87,Sva02,Zil01}, 
which is frequently used for continuum-based topology 
optimization \cite{BS03}, also solves a sequence of convex optimization 
approximations of the original problem. 
MMA approximates the objective function and the constraint functions by 
convex linear fractional functions by using the function values and 
gradients as well as some parameters controlling the vertical asymptotes 
of the generated functions. 
Also, sequential parametric convex approximation methods with 
application to truss optimization can be found in \cite{BBtT10,BtJKNZ00}.

\subsection{Heuristic for convex optimization with complementarity constraints}
\label{sec:algorithm.heuristic}

In this paper we attempt to solve problems having the following form: 
\begin{subequations}\label{P.convex.CC}%
  \begin{alignat}{3}
    & \MIN  &{\quad}& 
    f(\bi{\xi},\bi{y},\bi{z}) \\
    & \ST && 
    (\bi{\xi}, \bi{y}, \bi{z}) \in \varOmega , \\
    & && 
    \bi{y} \ge \bi{0} , \\
    & && 
    \bi{z} \ge \bi{0} , \\
    & && 
    \bi{y}^{\top} \bi{z} = 0 . 
    \label{P.convex.CC.5}
  \end{alignat}
\end{subequations}
Here, $f : \Re^{l} \times \Re^{n} \times \Re^{n} \to \Re$ is convex, 
$\varOmega \subseteq \Re^{l} \times \Re^{n} \times \Re^{n}$ is closed 
and convex, and the optimization variables are 
$\bi{\xi} \in \Re^{l}$, $\bi{y} \in \Re^{n}$, and $\bi{z} \in \Re^{n}$. 
Problem \eqref{P.convex.CC} is convex optimization with complementarity 
constraints. 

Following the idea in \cite{JPW16}, we can reduce problem 
\eqref{P.convex.CC} to a DC programming problem as follows. 
Consider a differentiable function 
$\phi : \Re^{n} \times \Re^{n} \to \Re$ satisfying 
\begin{align*}
  \phi(\bi{y}, \bi{z}) = 0 
  &\quad\Leftrightarrow\quad
  \bi{y}^{\top} \bi{z} = 0 , \\
  \phi(\bi{y}, \bi{z}) \ge 0 
  &\quad\Leftarrow\quad
  \bi{y} \ge \bi{0} , \ \bi{z} \ge \bi{0} . 
\end{align*}
The complementarity constraints, \eqref{P.convex.CC.5}, can be replaced 
by a penalization term as follows: 
\begin{subequations}\label{P.convex.CC.penalty.1}%
  \begin{alignat}{3}
    & \MIN  &{\quad}& 
    f(\bi{\xi},\bi{y},\bi{z}) + \rho \phi(\bi{y},\bi{z}) \\
    & \ST && 
    (\bi{\xi}, \bi{y}, \bi{z}) \in \varOmega , \\
    & && 
    \bi{y} \ge \bi{0} , \\
    & && 
    \bi{z} \ge \bi{0} . 
  \end{alignat}
\end{subequations}
Here, $\rho > 0$ is a penalty parameter. 
For sufficiently large $\rho$, problem \eqref{P.convex.CC.penalty.1} is 
equivalent to problem \eqref{P.convex.CC}. 
We next assume that $\phi$ can be decomposed as 
\begin{align*}
  \phi(\bi{y},\bi{z}) 
  = \phi_{+}(\bi{y},\bi{z}) - \phi_{-}(\bi{y},\bi{z}) , 
\end{align*}
where $\phi_{+}$ and $\phi_{-}$ are convex. 
Then problem \eqref{P.convex.CC.penalty.1} is reduced to the following 
form: 
\begin{subequations}\label{P.convex.CC.penalty.2}%
  \begin{alignat}{3}
    & \MIN  &{\quad}& 
    (f(\bi{\xi},\bi{y},\bi{z}) + \rho \phi_{+}(\bi{y},\bi{z}) )
    - \rho \phi_{-}(\bi{y},\bi{z}) \\
    & \ST && 
    (\bi{\xi}, \bi{y}, \bi{z}) \in \varOmega , \\
    & && 
    \bi{y} \ge \bi{0} , \\
    & && 
    \bi{z} \ge \bi{0} . 
  \end{alignat}
\end{subequations}
This is a DC programming problem, because 
$f(\bi{\xi},\bi{y},\bi{z}) + \rho \phi_{+}(\bi{y},\bi{z})$ and 
$\rho \phi_{-}(\bi{y},\bi{z})$ are convex. 
There exist several different choices for $\phi$, $\phi_{+}$, and 
$\phi_{-}$ \cite{JPW16,LtPd11}. 
In this paper we adopt 
\begin{align}
  \phi(\bi{y},\bi{z}) 
  &= \| \bi{y} + \bi{z} \|^{2} - \| \bi{y} - \bi{z} \|^{2} , 
  \label{def.function.phi} \\
  \phi_{+}(\bi{y},\bi{z}) 
  &= \| \bi{y} + \bi{z} \|^{2} , 
  \label{def.function.phi.+} \\
  \phi_{-}(\bi{y},\bi{z}) 
  &= \| \bi{y} - \bi{z} \|^{2} . 
  \label{def.function.phi.-}
\end{align}

To solve problem \eqref{P.convex.CC.penalty.1}, we apply the 
concave-convex procedure to problem \eqref{P.convex.CC.penalty.2} with 
gradually increasing the penalty parameter, $\rho$. 
For point 
$(\bi{y}^{(k)},\bi{z}^{(k)}) \in \Re^{n} \times \Re^{n}$, 
define $\hat{\phi}_{-}(\,\cdot\,; \bi{y}^{(k)},\bi{z}^{(k)}) : 
  \Re^{n} \times \Re^{n} \to \Re$ by
\begin{align}
  & \hat{\phi}_{-}(\bi{y},\bi{z}; \bi{y}^{(k)},\bi{z}^{(k)}) 
  = \phi_{-}(\bi{y}^{(k)},\bi{z}^{(k)}) \notag\\
  &\qquad
  + \nabla_{\bi{y}}\phi_{-}(\bi{y}^{(k)},\bi{z}^{(k)})^{\top}(\bi{y} - \bi{y}^{k})
  + \nabla_{\bi{z}}\phi_{-}(\bi{y}^{(k)},\bi{z}^{(k)})^{\top}(\bi{z} - \bi{z}^{k}).
  \label{def.function.hat.phi.-}
\end{align}
The proposed algorithm updates the solution by letting 
$(\bi{\xi}^{(k+1)},\bi{y}^{(k+1)},\bi{z}^{(k+1)})$ 
be an optimal solution of the following convex optimization problem: 
\begin{subequations}\label{P.subproblem.1}%
  \begin{alignat}{3}
    & \MIN  &{\quad}& 
    f(\bi{\xi},\bi{y},\bi{z}) 
    + \rho_{k} \phi_{+}(\bi{y},\bi{z}) 
    - \rho_{k} \hat{\phi}_{-}(\bi{y},\bi{z}; \bi{y}^{(k)},\bi{z}^{(k)}) \\
    & \ST && 
    (\bi{\xi}, \bi{y}, \bi{z}) \in \varOmega , \\
    & && 
    \bi{y} \ge \bi{0} , \\
    & && 
    \bi{z} \ge \bi{0} . 
  \end{alignat}
\end{subequations}
A reasonable stopping criterion is that the residual of the 
complementarity constraints is small enough, i.e., 
\begin{align}
  \phi(\bi{y}^{(k+1)}, \bi{z}^{(k+1)}) \le \epsilon_{1} , 
  \label{eq.stopping.inequality.1}
\end{align}
and the update of the incumbent solution is small enough, i.e., 
\begin{align}
  \| (\bi{\xi}^{(k+1)},\bi{y}^{(k+1)},\bi{z}^{(k+1)}) 
  - (\bi{\xi}^{(k)},\bi{y}^{(k)},\bi{z}^{(k)}) \| \le \epsilon_{2} ,
  \label{eq.stopping.inequality.2}
\end{align}
where $\epsilon_{1}$, $\epsilon_{2} > 0$ are thresholds. 
The algorithm is formally stated in \refalg{alg:concave.convex.penalty}.\footnote{
Choice of an initial point in the numerical experiments 
is explained in section~\ref{sec:ex}. } 

\begin{algorithm}
  \caption{penalty concave-convex procedure for convex optimization with 
  complementarity constraints}
  \label{alg:concave.convex.penalty}
  \begin{algorithmic}[1]
    \Require
    $\bi{\xi}^{(0)} \in \Re^{l}$, $\bi{y}^{(0)} \in \Re^{n}$, 
    $\bi{z}^{(0)} \in \Re^{n}$, $\rho_{0} > 0$, 
    $\rho_{\rr{max}} > \rho_{0}$, and $\mu > 1$. 
    \\
    $k \gets 0$. 
    \Repeat
    \State
    Let $(\bi{\xi}^{(k+1)},\bi{y}^{(k+1)},\bi{z}^{(k+1)})$
    be an optimal solution of problem \eqref{P.subproblem.1}. 
    \State
    $\rho_{k+1} := \min\{ \mu \rho_{k} , \rho_{\rr{max}} \}$. 
    \State
    Set $k \gets k+1$. 
    \Until
    stopping criterion is satisfied. 
  \end{algorithmic}
\end{algorithm}

\begin{remark}
  \refalg{alg:concave.convex.penalty} is designed essentially based on 
  the algorithm proposed by \citet{LB16} for solving problem \eqref{P.DC}. 
  In their algorithm, the following subproblem is solved to update 
  $\bi{x}^{(k)}$: 
%\begin{subequations}
  \begin{alignat*}{3}
    & \MIN  &{\quad}& 
    f_{0}(\bi{x}) - \hat{g}_{0}(\bi{x}; \bi{x}^{(k)}) 
    + \rho_{k} \sum_{i=1}^{m} s_{i} \\
    & \ST && 
    f_{i}(\bi{x}) - \hat{g}_{i}(\bi{x}; \bi{x}^{(k)}) \le s_{i} , 
    \quad i=1,\dots,m , \\
    & && 
    s_{i} \ge 0 , 
    \quad i=1,\dots,m . 
  \end{alignat*}
%\end{subequations}
  Thus, penalization terms for all the constraints are added to the 
  objective function by using the $\ell_{1}$-exact penalty function. 
  In contrast, in \refalg{alg:concave.convex.penalty}  only the 
  complementarity constraints are penalized, and the other constraints 
  of the original optimization problem 
  are satisfied at the solution of the subproblem. 
  \finbox
\end{remark}

\section{Mixed-integer semidefinite programming formulation for robust 
  truss topology optimization}
\label{sec:integer}

In section~\ref{sec:preliminary.robust}, we recall the existing MISDP 
formulation for robust truss topology optimization, without considering 
overlapping members in the ground structure. 
Section~\ref{sec:preliminary.node} presents treatment of overlapping 
members within the framework of MISDP. 

\subsection{Review of existing formulation}
\label{sec:preliminary.robust}

In this section we briefly review an MISDP formulation of the robust 
truss topology optimization under the load uncertainty \cite{YK10}; see 
also \cite{BtN97}. 

Following the ground structure method, consider a truss consisting of 
candidate members connected by nodes. 
Let $m$ and $d$ denote the number of the members and the number of 
degrees of freedom of the nodal displacements,\footnote{
The degrees of freedom of a truss is the possible 
components of the nodal displacements that define the configuration of 
the truss. }   respectively.
We use $x_{i}$ $(i=1,\dots,m)$ to denote the member cross-sectional 
areas, which are design variables to be optimized. 
Throughout the paper, we assume small deformation and linear elasticity. 

Let $t_{i} \in \{ 0,1 \}$ be a variable that serves as an indicator of 
existence of member $i$ such that $t_{i}=1$ means that member $i$ exists 
and $t_{i}=0$ means that it vanishes. 
We use $\overline{x} > 0$ and 
$\underline{x} \in [0, \overline{x}]$ 
to denote the specified upper and lower bounds for the cross-sectional 
area of an existing member, i.e., $x_{i}$ should satisfy 
$x_{i} \in \{ 0 \} \cup [\underline{x}, \overline{x}]$. 
This constraint can be written by using $t_{i}$ as 
\begin{align}
  \underline{x} t_{i} &\le x_{i} \le \overline{x} t_{i} .  
  \label{eq.member.integer.1}
\end{align}

We next introduce $s_{j} \in \{ 0,1 \}$ $(j=1,\dots,d)$ to represent 
the existence of the $j$th degree of freedom. 
A node in a ground structure is removed if and only if all the members 
connected to the nodes vanish. 
Let $s_{j}=1$ mean that the node having the $j$th degree of freedom 
exists, and $s_{j}=0$ mean that it vanishes. 
We use $I(j) \subseteq \{ 1,\dots,m \}$ to denote the set of indices of 
the members connected to the node having the $j$th degree of freedom. 
Then $s_{j}$ is related to $t_{1},\dots,t_{m}$ as follows: 
\begin{align}
  t_{i}  &\le s_{j} , 
  \quad \forall i \in I(j) . 
  \label{eq.node.integer.1}
\end{align}

Let $K(\bi{x}) \in \Re^{d \times d}$ denote the stiffness matrix of a 
truss, which is a (matrix-valued) linear function of $\bi{x}$. 
For a given external load, denoted $\bi{p} \in \Re^{d}$, the compliance 
of the truss is defined by 
\begin{align}
  \pi(\bi{x};\bi{p}) 
  = \sup \{ 2\bi{p}^{\top} \bi{u} - \bi{u}^{\top} K(\bi{x}) \bi{u} 
  \mid \bi{u} \in \Re^{d} \} . 
  \label{eq:def.compliance.0}
\end{align}
The conventional compliance minimization problem is formulated in 
variables $\bi{x}$ as follows: 
\begin{subequations}\label{P.nominal.compliance}%
  \begin{alignat}{3}
    & \MIN  &{\quad}& 
    \pi(\bi{x};\bi{p}) \\
    & \ST && 
    \bi{x} \ge \bi{0} , 
    \label{P.nominal.compliance.2} \\
    & &&
    \bi{c}^{\top} \bi{x} \le \overline{c} . 
    \label{P.nominal.compliance.3}
  \end{alignat}
\end{subequations}
Here, $c_{i}$ is the undeformed length of member $i$, 
$\bi{c} = (c_{1},\dots,c_{m})^{\top}$, and 
$\overline{c} > 0$ is the specified upper bound for the structural volume. 

The uncertainty model of the external load is defined as follows: 
Let $\tilde{\bi{p}} \in \Re^{d}$ denote the nominal value (or the best 
estimate) of the external load. 
Define a constant matrix $Q \in \Re^{d \times d}$ by 
\begin{align*}
  Q = 
  \left[
  \begin{array}{@{}c|c|c|c@{}}
    &  &  &  \\
    \tilde{\bi{p}} & \alpha \bi{q}_{1} & \cdots & \alpha \bi{q}_{d-1} \\
    &  &  &  \\
  \end{array}
  \right] , 
\end{align*}
where $\bi{q}_{1},\dots,\bi{q}_{d-1} \in \Re^{d}$ are the orthonormal 
basis vectors of the orthogonal complement of $\tilde{\bi{p}}$, and 
$\alpha > 0$ is a constant representing the level of uncertainty. 
Then the uncertainty set of the external load, i.e., the set of all 
possible realizations of the external load, is defined by 
\begin{align}
  P(\bi{s}) 
  = \{ \diag(\bi{s}) Q \bi{e} \mid \| \bi{e} \| \le 1  \} . 
  \label{eq.def.uncertainty.set.s}
\end{align}
For example, suppose that $\tilde{\bi{p}}$ has only one nonzero 
component. 
Then, without loss of generality we can assume 
$\tilde{p}_{1} \not= 0$, and we have that 
\begin{align}
  \tilde{\bi{p}} = 
  \begin{bmatrix}
    \tilde{p}_{1} \\
    0 \\
    0 \\
    \vdots \\
    0 \\
  \end{bmatrix}  , 
  \quad
  Q = 
  \begin{bmatrix}
    \tilde{p}_{1} & 0 & 0 & \cdots & 0 \\
    0 & \alpha & 0 & \cdots  & 0 \\
    0 & 0 & \alpha & \cdots  & 0 \\
    \vdots & \vdots & \vdots  & \ddots  & \vdots  \\
    0 & 0 & 0 & \cdots & \alpha \\
  \end{bmatrix}
  . 
  \label{eq.ex.uncertainty.1}
\end{align}

Let $J_{\rr{f}} \subseteq \{ 1,\dots,d \}$ denote the set of indices of 
nonzero components of $\tilde{\bi{p}}$, i.e., 
\begin{align*}
  J_{\rr{f}} 
  = \{ j \in \{ 1,\dots,d \} \mid \tilde{p}_{j} \not= 0 \} . 
\end{align*}
The nodes to which the nominal external load, $\tilde{\bi{p}}$, is 
applied should not be removed in the course of optimization. 
Therefore, we impose the following constraint: 
\begin{align}
  s_{j} = 1 , 
  \quad \forall j \in J_{\rr{f}} .
  \label{eq.fix.loading.node}
\end{align}

In robust optimization, we attempt to find a truss design that minimizes 
the maximal compliance (i.e., the worst-case compliance) when the 
external load can take any value in $P(\bi{s})$. 
With reference to \eqref{eq.member.integer.1}, 
\eqref{eq.node.integer.1}, \eqref{P.nominal.compliance}, 
\eqref{eq.def.uncertainty.set.s}, and \eqref{eq.fix.loading.node}, 
we can see that this optimization problem can be formulated as follows: 
\begin{subequations}\label{P.robust.0}%
  \begin{alignat}{3}
    & \MIN  &{\quad}& 
    \sup\{ \pi(\bi{x};\bi{p}) \mid \bi{p} \in P(\bi{s}) \} \\
    & \ST && 
    s_{j} = 1 , 
    \quad \forall j \in J_{\rr{f}} , \\
    & &&
    t_{i} \le s_{j} , 
    \quad  \forall i \in I(j); \ j=1,\dots,d , \\
    & &&
    \underline{x} \bi{t} \le \bi{x} \le \overline{x} \bi{t} , \\
    & &&
    \bi{c}^{\top} \bi{x} \le \overline{c} , \\
    & &&
    \bi{s} \in \{ 0,1 \}^{d} , \\
    & &&
    \bi{t} \in \{ 0,1 \}^{m} . 
  \end{alignat}
\end{subequations}

For $\bi{x} \in \Re^{m}$ $(\bi{x} \ge \bi{0})$, $\bi{s} \in \{ 0,1 \}^{d}$, 
and $w \in \Re$, define $W(\bi{x},\bi{s},w) \in \SC^{d+1}$ by 
\begin{align*}
  W(\bi{x},\bi{s},w) = 
  \begin{bmatrix}
    w I & (\diag(\bi{s}) Q)^{\top} \\
    \diag(\bi{s}) Q & K(\bi{x}) \\
  \end{bmatrix}
  . 
\end{align*}
It is shown in \cite[Lemma~2.2]{BtN97} that $w \in \Re$ satisfies 
\begin{align*}
  w \ge  \sup\{ \pi(\bi{x};\bi{p}) \mid \bi{p} \in P(\bi{s}) \} 
\end{align*}
if and only if 
\begin{align}
  \begin{bmatrix}
    w I & (\diag(\bi{s}) Q)^{\top} \\
    \diag(\bi{s}) Q & K(\bi{x}) \\
  \end{bmatrix}
  \succeq O
  \label{eq.matrix.inequality.lemma}
\end{align}
holds. 
Consequently, problem~\eqref{P.robust.0} is equivalently rewritten as 
follows: 
\begin{subequations}\label{P.robust.1}%
  \begin{alignat}{3}
    & \MIN  &{\quad}& w \\
    & \ST && 
    W(\bi{x},\bi{s},w) \succeq O , 
    \label{P.robust.1.SDP} \\
    & &&
    s_{j} = 1 , 
    \quad \forall j \in J_{\rr{f}} , \\
    & &&
    t_{i} \le s_{j} , 
    \quad  \forall i \in I(j); \ j=1,\dots,d , 
    \label{P.robust.1.2} \\
    & &&
    \underline{x} \bi{t} \le \bi{x} \le \overline{x} \bi{t} , 
    \label{P.robust.1.4} \\
    & &&
    \bi{c}^{\top} \bi{x} \le \overline{c} , \\
    & &&
    \bi{s} \in \{ 0,1 \}^{d} , 
    \label{P.robust.1.3} \\
    & &&
    \bi{t} \in \{ 0,1 \}^{m} . 
    \label{P.robust.1.6}
  \end{alignat}
\end{subequations}

Problem \eqref{P.robust.1} is an MISDP problem. 
By relaxing the 0-1 constraints into linear inequality constraints, we 
obtain an SDP relaxation. 
Since SDP can be solved efficiently with a primal-dual interior-point 
method, we can find a global optimal solution of problem 
\eqref{P.robust.1} with a branch-and-bound method \cite{YK10}.

\subsection{Treatment of members lying on a line}
\label{sec:preliminary.node}

As explained in section~\ref{sec:motivation}, for the robust truss 
topology optimization it is necessary to incorporate overlapping members 
to a ground structure. 
Since the existence of overlapping members in a final 
truss design is not accepted, it is required to incorporate the 
constraints prohibiting the presence of overlapping members in a truss design. 
To the best of the author's knowledge, such a consideration cannot be 
found in literature on robust truss topology optimization. 

\begin{figure}[tp]
  \centering
  \includegraphics[scale=0.55]{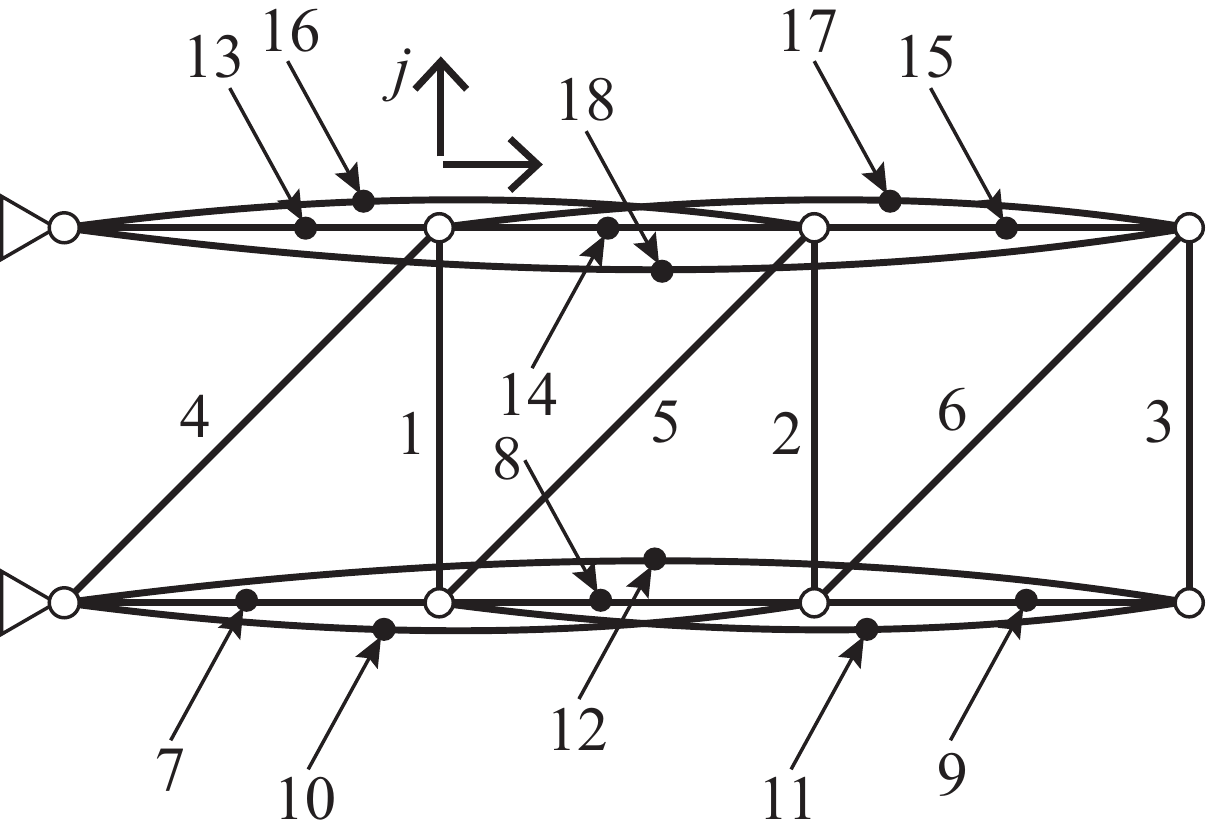}
  \caption{A ground structure consisting of $18$ members. 
  }
  \label{fig:gs_overlap}
\end{figure}

Recall that, for each $j=1,\dots,d$, $s_{j}=1$ means that the node having 
the $j$th degree of freedom exists, and $s_{j}=0$ means that it vanishes. 
Also, $t_{i}=1$ means that member $i$ exists and $t_{i}=0$ means that 
it vanishes. 
Let $L(j) \subseteq \{ 1,\dots,m \}$ denote the set of indices of the 
members lying on the node having the $j$th degree of freedom. 
\reffig{fig:gs_overlap} shows an example of ground structure with 
overlapping members. 
It has $m=18$ members and $d=12$ degrees of freedom of nodal 
displacements. 
In this example, we have $L(j)=\{ 16,18 \}$ and 
$I(j)=\{ 1,4,13,14,17 \}$. 
If $s_{j}=1$, then all the members in $L(j)$ cannot exist, i.e., 
$t_{i}=0$ $(\forall i \in L(J))$. 
Also, if there exists an $i \in L(j)$ such that $t_{i}=1$, then the 
corresponding node cannot exist, i.e., $s_{j}=0$. 
These two conditions can be formulated as 
\begin{align}
  s_{j} \le 1-t_{i} , 
  \quad \forall i \in L(j) . 
  \label{eq.additional.1}
\end{align}
In the following, we add \eqref{eq.additional.1} to problem 
\eqref{P.robust.1}. 

It is worth noting that a truss design involving a chain is infeasible 
for the presented robust optimization problem, because it is unstable 
and cannot be in equilibrium with uncertain forces applied to an 
intermediate node of the chain. 
Therefore, a solution obtained by the proposed method does not involve a 
member that is longer than the maximum member length of the ground 
structure. 
This might be considered an advantage of the robust topology 
optimization, because in a conventional truss topology optimization an 
optimal solution may possibly have a long chain and special treatment, 
such as the local buckling constraints \cite{GCO05,Mel14}, 
is necessary for avoiding presence of a too long member.

\section{Simple heuristic for robust truss topology optimization}
\label{sec:heuristic}

The MISDP approach presented in 
section~\ref{sec:integer} is only applicable to small-size instances. 
Alternatively, in this section we present a heuristic having small 
computational cost. 
In section~\ref{sec:heuristic.complementarity}, we reformulate our 
robust truss topology optimization as SDPCC. 
A concave-convex procedure is then applied to this formulation in 
section~\ref{sec:heuristic.algorithm}.

\subsection{Formulation as semidefinite programming with complementarity constraints}
\label{sec:heuristic.complementarity}

In this section, we reformulate problem \eqref{P.robust.1} with 
constraint \eqref{eq.additional.1} as SDPCC, which is well suited for 
applying the concave-convex procedure. 

We begin with constraints \eqref{P.robust.1.2} and \eqref{P.robust.1.4}, 
which describe the relation between $s_{j}$ and  $\bi{x}$. 
Let $r_{j}$ $(j=1,\dots,d)$ denote the sum of the cross-sectional 
areas of the members that are connected to the node having the $j$th 
degree of freedom, i.e., 
\begin{align}
  r_{j} = \sum_{i \in I(j)} x_{i} . 
  \label{eq.vanishing.2}
\end{align}
Observe that $r_{j}>0$ implies $s_{j}=1$, because at least one member 
connected to the corresponding node exists and, hence, 
the node should exist. 
Also, $s_{j}=0$ implies that the corresponding node vanishes, and hence 
$r_{j}=0$, i.e., all the members connected to the node should vanish. 
These two assertions can be written as 
\begin{align}
  (1 - s_{j}) r_{j} = 0 . 
  \label{eq.transform.to.compatibility.1}
\end{align}
Namely, constraint \eqref{P.robust.1.2} can be replaced with  
\eqref{eq.transform.to.compatibility.1}. 
For notational simplicity, in the following we write 
\eqref{eq.vanishing.2} as 
\begin{align*}
  \bi{r} = R \bi{x} 
\end{align*}
with a constant matrix $R \in \Re^{d \times m}$. 

We next consider constraints \eqref{P.robust.1.4} and 
\eqref{eq.additional.1}, which describe the relation between 
$s_{j}$ and $\bi{x}$. 
Let $v_{j}$ $(j=1,\dots,d)$ denote the sum of the cross-sectional areas 
of the members that lying across the node having the $j$th degree of 
freedom, i.e., 
\begin{align}
  v_{j} = \sum_{i \in L(j)} x_{i} . 
  \label{eq.vanishing.4}
\end{align}
Observe that $v_{j}>0$ implies $s_{j}=0$, because at least one member 
lying across the corresponding node exists and, thence, the node should 
vanish. 
Also, $s_{j}=1$ implies that the corresponding node exists, and hence 
all the members lying across the node should vanish, i.e., $v_{j}=0$. 
These two assertions can be written as 
\begin{align}
  s_{j} v_{j}  = 0 . 
  \label{eq.transform.to.compatibility.2}
\end{align}
Namely, constraint \eqref{eq.additional.1} can be replaced with 
\eqref{eq.transform.to.compatibility.2}, 
For notational simplicity, we write \eqref{eq.vanishing.4} as 
\begin{align*}
  \bi{v} = V \bi{x} 
\end{align*}
by using a constant matrix $V \in \Re^{d \times m}$. 

Finally, consider constraint \eqref{P.robust.1.4}. 
For each $i=1,\dots,m$, we introduce a new variable $z_{i} \in \Re$ so 
that $\underline{x}-z_{i}$ corresponds to the lower bound for the 
cross-sectional area of member $i$. 
If $x_{i}>0$, then member $i$ exists and the lower bound should be 
$\underline{x}$, which means $z_{i}=0$. 
Also, when the lower bound becomes smaller than $\underline{x}$ (i.e., 
$z_{i}>0$), then member $i$ should vanish (i.e., $x_{i}=0$), and hence we 
set $z_{i}=\underline{x}$. 
These relations can be written as follows: 
\begin{align}
  \underline{x} - z_{i} 
  &\le x_{i} \le \overline{x} , 
  \label{eq.transform.to.compatibility.3} \\
  0 &\le z_{i} \le \underline{x} , 
  \label{eq.transform.to.compatibility.4} \\
  x_{i} z_{i} &= 0 . 
  \label{eq.transform.to.compatibility.5}
\end{align}
Consequently, constraints \eqref{P.robust.1.4}  and \eqref{P.robust.1.6} 
 can be replaced with 
\eqref{eq.transform.to.compatibility.3}, 
\eqref{eq.transform.to.compatibility.4}, and 
\eqref{eq.transform.to.compatibility.5}.

The upshot is that problem \eqref{P.robust.1} incorporating constraint 
\eqref{eq.additional.1} is equivalently rewritten as follows: 
\begin{subequations}\label{P.complementarity.1}%
  \begin{alignat}{3}
    & \MIN  &{\quad}& w \\
    & \ST && 
    W(\bi{x},\bi{s},w) \succeq O , 
    \label{P.complementarity.1.psd} \\
    & &&
    s_{j} = 1 , 
    \quad \forall j \in J_{\rr{f}} , \\
    & &&
    \bi{r} = R \bi{x} , \\
    & &&
    \bi{v} = V \bi{x} , \\
    & &&
    \bi{0} \le \bi{s} \le \bi{1} , 
    \label{P.complementarity.1.s} \\
    & &&
    \underline{x} \bi{1} - \bi{z} \le \bi{x} \le \overline{x} \bi{1} , \\
    & &&
    \bi{0} \le \bi{z} \le \underline{x} \bi{1} , \\
    & &&
    \bi{c}^{\top} \bi{x} \le \overline{c} , \\
    & &&
    (1 - s_{j}) r_{j} = 0 , 
    &{\quad}&  j=1,\dots,d , 
    \label{P.complementarity.1.4} \\
    & &&
    s_{j} v_{j} = 0 , 
    &{\quad}&  j=1,\dots,d , 
    \label{P.complementarity.1.6} \\
    & &&
    x_{i} z_{i} = 0 , 
    &{\quad}&  i=1,\dots,m . 
    \label{P.complementarity.1.8}
  \end{alignat}
\end{subequations}
Here, $\bi{x}$, $\bi{z}$, $\bi{s}$, $\bi{r}$, $\bi{v}$, and 
$w$ are variables to be optimized. 
observe that any feasible solution of problem \eqref{P.complementarity.1} 
satisfies 
$1-s_{j} \ge 0$, $r_{j} \ge 0$, $s_{j} \ge 0$, $v_{j} \ge 0$, 
$x_{i} \ge 0$, and $z_{i} \ge 0$. 
Therefore, constraints 
\eqref{P.complementarity.1.4}, \eqref{P.complementarity.1.6}, and 
\eqref{P.complementarity.1.8} are complementarity constraints. 
Constraint \eqref{P.complementarity.1.psd} is a linear matrix inequality 
constraint in terms of $\bi{x}$, $\bi{s}$, and $w$. 
Thus, problem \eqref{P.complementarity.1} has the form of SDPCC.

\begin{remark}
  Since the algorithm presented in this paper consists 
  of sequential approximation, adding some linear inequalities 
  may possibly limit the search spae and enhance the convergence. 
  In the following, we consider linear valid inequalities, which 
  naturally stem from the complementarity constraints and can be handled 
  effectively in the numerical solution. 
  Suppose that two variables, $\alpha$, $\beta \in \Re$, are subjected 
  to the complementarity constraint and their upper bounds are given, 
  i.e., 
  \begin{align}
    0 &\le \alpha \le \bar{\alpha} , 
    \label{eq.valid.inequality.1} \\
    0 &\le \beta \le \bar{\beta} , 
    \label{eq.valid.inequality.2} \\
    \alpha \beta  &= 0 ,
    \label{eq.valid.inequality.3}
  \end{align}
  where $\bar{\alpha}$ and $\bar{\beta}$ are positive constants. 
  It is known that the inequality 
  \begin{align*}
    \bar{\beta} \alpha + \bar{\alpha} \beta 
    \le \bar{\alpha} \bar{\beta} 
  \end{align*}
  serves as a valid constraint for \eqref{eq.valid.inequality.1}, 
  \eqref{eq.valid.inequality.2}, and \eqref{eq.valid.inequality.3} 
  \cite{MPY12,YMP16}. 
  In the same manner, we can construct valid constraints for problem 
  \eqref{P.complementarity.1}. 
  Concerning \eqref{P.complementarity.1.4}, observe that we obtain 
  \begin{align*}
    1 - s_{j}  &\le 1 , \\
    r_{j}  &\le \overline{x} |I(j)|
  \end{align*}
  from \eqref{P.complementarity.1.s} and \eqref{eq.vanishing.2}, 
  respectively. 
  Therefore, the constraints 
  \begin{alignat}{2}
    {-}\overline{x} |I(j)| s_{j} + r_{j} 
    &\le 0 , 
    &{\quad}& j=1,\dots,d 
    \label{eq.valid.add.1}
  \end{alignat}
  are valid for \eqref{P.complementarity.1.4}. 
  Similarly, inequalities 
  \begin{alignat}{2}
    \overline{x} |L(j)| s_{j} + v_{j} 
    & \le \overline{x} |L(j)| , 
    &{\quad}& j=1,\dots,d, 
    \label{eq.valid.add.2} \\
    \underline{x} x_{i} + \overline{x} z_{i} 
    & \le \underline{x} \overline{x} , 
    &{\quad}& i=1,\dots,m
    \label{eq.valid.add.3}
  \end{alignat}
  are valid constraints for \eqref{P.complementarity.1.6} and 
  \eqref{P.complementarity.1.8}, respectively. 
  In the following, we add constraints \eqref{eq.valid.add.1}, 
  \eqref{eq.valid.add.2}, and \eqref{eq.valid.add.3} to problem 
  \eqref{P.complementarity.1}. 
  \finbox
\end{remark}

\subsection{Penalty concave-convex procedure for robust truss topology optimization}
\label{sec:heuristic.algorithm}

Problem \eqref{P.complementarity.1} has the form of problem 
\eqref{P.convex.CC} studied in section~\ref{sec:algorithm.heuristic}. 
To see this, it is convenient to rewrite problem 
\eqref{P.complementarity.1} as follows: 
\begin{subequations}\label{P.complementarity.2}%
  \begin{alignat}{3}
    & \MIN  &{\quad}& w \\
    & \ST && 
    (\bi{x},\bi{z},\bi{s},\bi{r},\bi{v},w)  \in F , \\
    & &&
    \bi{1}-\bi{s} \ge \bi{0} , \quad
    \bi{r} \ge \bi{0} , \quad
    (\bi{1}-\bi{s})^{\top} \bi{r} = 0 , 
    \label{P.complementarity.2.3} \\
    & &&
    \bi{s} \ge \bi{0} , \quad
    \bi{v} \ge \bi{0} , \quad
    \bi{s}^{\top} \bi{v} = 0 , 
    \label{P.complementarity.2.4} \\
    & &&
    \bi{x} \ge \bi{0} , \quad
    \bi{z} \ge \bi{0} , \quad
    \bi{x}^{\top} \bi{z} = 0 . 
    \label{P.complementarity.2.5}
  \end{alignat}
\end{subequations}
Here, $F$ is defined by 
\begin{align*}
  F = \{ (\bi{x},\bi{z},\bi{s},\bi{r},\bi{v},w) 
  \mid\, & 
  W(\bi{x},\bi{s},w) \succeq O , \
  \bi{r} = R \bi{x} , \
  \bi{v} = V \bi{x} ,  \notag\\
  & 
  \bi{0} \le \bi{s} \le \bi{1} , \
  s_{j} = 1 \ (\forall j \in J_{\rr{f}}) , \notag\\
  & \underline{x} \bi{1} - \bi{z} \le \bi{x} \le \overline{x} \bi{1} , \
  \bi{0} \le \bi{z} \le \underline{x} \bi{1} , \
  \bi{c}^{\top} \bi{x} \le \overline{c} , \notag\\
  & {-}\overline{x} |I(j)| s_{j} + r_{j} \le 0 \ (j=1,\dots,d), \notag\\
  & \overline{x} |L(j)| s_{j} + v_{j}  \le \overline{x} |L(j)| 
  \ (j=1,\dots,d), \notag\\
  & \underline{x} \bi{x} + \overline{x} \bi{z} 
  \le \underline{x} \overline{x} \bi{1} 
  \} , 
\end{align*}
which is a convex set. 
The complementarity constraints in \eqref{P.complementarity.2.3}, 
\eqref{P.complementarity.2.4}, and \eqref{P.complementarity.2.5} can be 
replaced with penalization terms added to the objective function as 
follows: 
\begin{subequations}\label{P.complementarity.3}%
  \begin{alignat}{3}
    & \MIN  &{\quad}& 
    w + \rho \phi(\bi{1}-\bi{s}, \bi{r}) + \rho \phi(\bi{s}, \bi{v}) 
    + \rho \phi(\bi{x}, \bi{z}) \\
    & \ST && 
    (\bi{x},\bi{z},\bi{s},\bi{r},\bi{v},w)  \in F . 
  \end{alignat}
\end{subequations}
Here, $\rho > 0$ is a sufficiently large penality parameter, 
and $\phi$ has been defined by \eqref{def.function.phi}. 
Problem \eqref{P.complementarity.3} has the same form as problem 
\eqref{P.convex.CC.penalty.1}. 

We are now in position to apply \refalg{alg:concave.convex.penalty} to 
problem \eqref{P.complementarity.3}. 
\refalg{alg:concave.convex.penalty} solves the subproblem in 
\eqref{P.subproblem.1}, which is explicitly written as follows:  
\begin{subequations}\label{P.complementarity.4}%
  \begin{alignat}{3}
    & \MIN  &{\quad}& 
    w + \rho_{k} \phi_{+}(\bi{1}-\bi{s}, \bi{r}) 
    + \rho_{k} \phi_{+}(\bi{s}, \bi{v}) 
    + \rho_{k} \phi_{+}(\bi{x}, \bi{z})  \notag\\
    & && \quad
    - \rho_{k} \hat{\phi}_{-}(\bi{1}-\bi{s},\bi{r}; \bi{1}-\bi{s}^{(k)},\bi{r}^{(k)})
    - \rho_{k} \hat{\phi}_{-}(\bi{s},\bi{v}; \bi{s}^{(k)},\bi{v}^{(k)}) 
    \notag\\
    & && \quad
    - \rho_{k} \hat{\phi}_{-}(\bi{x},\bi{z}; \bi{x}^{(k)},\bi{z}^{(k)})
    \\
    & \ST && 
    (\bi{x},\bi{z},\bi{s},\bi{r},\bi{v},w)  \in F . 
  \end{alignat}
\end{subequations}
Here, $\phi_{+}$ and $\hat{\phi}_{-}$ are defined by 
\eqref{def.function.phi.+} and \eqref{def.function.hat.phi.-}, 
respectively. 
Since the constant terms in the objective function can be neglected, 
problem \eqref{P.complementarity.4} can be reduced to the following 
problem: 
\begin{subequations}\label{P.DC.1}%
  \begin{alignat}{3}
    & \MIN  &{\quad}& 
    w + \rho_{k} ( \| \bi{x} + \bi{z} \|^{2}
    + \| \bi{1} - \bi{s} + \bi{r} \|^{2} + \| \bi{s} + \bi{v} \|^{2})
    \notag\\
    & && \quad
    -2\rho_{k} (\bi{x}^{(k)} - \bi{z}^{(k)})^{\top} \bi{x} 
    -2\rho_{k} (\bi{z}^{(k)} - \bi{x}^{(k)})^{\top} \bi{z} 
    \notag\\
    & && \quad
    - 2\rho_{k} (2 \bi{s}^{(k)} + \bi{r}^{(k)} - \bi{v}^{(k)} - \bi{1})^{\top} \bi{s} 
    \notag\\
    & && \quad
    - 2\rho_{k} (\bi{s}^{(k)} + \bi{r}^{(k)} -\bi{1})^{\top} \bi{r}  
    - 2\rho_{k} (\bi{v}^{(k)} - \bi{s}^{(k)})^{\top} \bi{v} \\
    & \ST && 
    (\bi{x},\bi{z},\bi{s},\bi{r},\bi{v},w)  \in F . 
  \end{alignat}
\end{subequations}
Problem \eqref{P.DC.1} is a minimization problem of a convex quadratic 
function under a linear matrix inequality constraint. 
Hence, this problem can be recast as SDP. 
Thus, at each iteration of \refalg{alg:concave.convex.penalty} we solve 
an SDP problem. 

As mentioned in section~\ref{sec:algorithm.heuristic}, a reasonable 
stopping criterion is that \eqref{eq.stopping.inequality.1} and 
\eqref{eq.stopping.inequality.2} are satisfied. 
In practice, however, we might use a relaxed criterion, which may save 
some iterations before convergence. 
Specifically, we terminate the algorithm when either 
\eqref{eq.stopping.inequality.1} or 
\begin{align}
  \| \bi{x}^{(k+1)} - \bi{x}^{(k)} \| \le \epsilon_{2}
  \label{eq.stopping.inequality.3}
\end{align}
is satisfied. 
Then, by using the obtained solution we can fix the set of existing 
members and the set of existing nodes. 
Fixing these sets means that one variable in all the complementarity 
constraints in problem \eqref{P.complementarity.1} is fixed. 
Therefore, the problem is now becomes SDP, which is to be 
solved as the post process. 
The solutions presented in section~\ref{sec:ex} are obtained in this 
manner.

\section{Numerical experiments}
\label{sec:ex}

This section reports three numerical experiments. 

The proposed algorithm was implemented in MATLAB ver.~9.0. 
At each iteration we solved an SDP problem in \eqref{P.DC.1} by using CVX, 
a MATLAB package for specifying and solving convex optimization 
problems \cite{GB08,CVX}. 
SDPT3 ver.~4.0 \cite{TTT03} was used as the solver. 
Computation was carried out on a $2.2\,\mathrm{GHz}$ Intel Core i5 
processor with $8\,\mathrm{GB}$ RAM. 
The Young modulus of the trusses in the following numerical examples is 
$20\,\mathrm{GPa}$. 

The initial point for \refalg{alg:concave.convex.penalty} is chosen as 
follows. 
We first solve problem \eqref{P.nominal.compliance} for the 
nominal external load $\tilde{\bi{p}}$, i.e., the compliance 
minimization without considering uncertainties,\footnote{
Problem \eqref{P.nominal.compliance} is convex. 
Various reformulations are known in literature; see, e.g., 
\cite{ABBTZ92,JKZ98}. 
For example, replacing $\diag(\bi{s})Q$ in 
\eqref{eq.matrix.inequality.lemma} with $\tilde{\bi{p}}$, one can 
readily obtain SDP that minimizes $w$ under constraint 
$
  \begin{bmatrix}
    w & \tilde{\bi{p}}^{\top} \\
    \tilde{\bi{p}} & K(\bi{x}) \\
  \end{bmatrix}
  \succeq O
$, \eqref{P.nominal.compliance.2}, and \eqref{P.nominal.compliance.3}. 
This formulation was used in the numerical experiments. 
It should be clear that a ground structure with overlapping members is 
used for generating the initial point, $\bi{x}^{(0)}$. 
} 
and let $\bi{x}^{(0)}$ be the obtained optimal solution. 
The initial values for the other variables are given by 
$\bi{z}^{(0)}=\bi{0}$, $\bi{s}^{(0)}=\bi{1}/2$, 
$\bi{r}^{(0)}=R \bi{x}^{(0)}$, and $\bi{v}^{(0)}=V \bi{x}^{(0)}$. 
The parameters of \refalg{alg:concave.convex.penalty} are 
$\rho_{0}=10^{-2}$, $\rho_{\rr{max}}=10^{6}$, and $\mu=1.5$. 
We terminate \refalg{alg:concave.convex.penalty} if either 
\eqref{eq.stopping.inequality.1} or \eqref{eq.stopping.inequality.3} is 
satisfied, where $\epsilon_{1}=2m\times 10^{-2} \,\mathrm{mm}^{2}$ and 
$\epsilon_{2}=10^{-2}\,\mathrm{mm}^{2}$. 
Then, as explained in section~\ref{sec:heuristic.algorithm}, 
we fix one variable in all the complementarity constraints in 
problem \eqref{P.complementarity.1}, and solve the resulting SDP problem 
to obtain the final solution. 
The settings of the initial point and the parameters 
explained above were determined by preliminary numerical experiments. 
In section~\ref{sec:ex.small}, we consider three problem instances that 
could be solved with a global optimization method. 
Cantilever truss examples with two different loading conditions, which 
are frequently solved in structural optimization, are considered in 
sections~\ref{sec:ex.eva4} and \ref{sec:ex.eva3}.

\subsection{Example (I): Comparison with global optimization}
\label{sec:ex.small}

\begin{table}[bp]
  \centering
  \caption{Characteristics of the problem instances in example (I).}
  \label{tab:ex.0.data}
  \begin{tabular}{lrrrrr}
    \toprule
    Problem & $m$ & $d$ & $\overline{c}$ $(\mathrm{mm}^{3})$ & Rob.\ opt.\ (J) & Nom.\ opt.\ (J) \\
    \midrule
    \reffig{fig:robust.x2_y1} & 14 &  8 & $0.4\times 10^{6}$ & $8984.375$ & $8000.000$ \\
    \reffig{fig:robust.x3_y3} & 98 & 24 & $1.8\times 10^{6}$ & $2442.708$ & $2006.944$ \\
    \reffig{fig:robust.x3_y2} & 35 & 18 & $1.2\times 10^{6}$ & $11093.750$ & $9375.000$ \\
    \bottomrule
  \end{tabular}
\end{table}

\begin{table}[bp]
  \centering
  \caption{Computational costs for example (I).}
  \label{tab:ex.0.result}
  \begin{tabular}{lrrrr}
    \toprule
    Problem & \multicolumn{2}{c}{Proposed method} & \multicolumn{2}{c}{YALMIP} \\
    \cmidrule(lr){2-3} \cmidrule(l){4-5}
    & {\#}iter.\ & Time (s) & {\#}iter.\ & Time (s) \\
    \midrule
    \reffig{fig:robust.x2_y1} & 3 & 4.1 & 13 & 2.7 \\
    \reffig{fig:robust.x3_y3} & 15 & 39.8 & 92 & 54.7 \\
    \reffig{fig:robust.x3_y2} & 47 & 59.5 & 1141 & 300.3 \\
    \bottomrule
  \end{tabular}
\end{table}

In this section, we consider the small-size instances presented in 
section~\ref{sec:motivation}. 
For comparison, the MISDP formulation, i.e., problem \eqref{P.robust.1} 
with constraint \eqref{eq.additional.1}, is solved with 
YALMIP~\cite{Lof04}. 
YALMIP finds a global optimal solution of an MISDP problem with a 
branch-and-bound method, at each iteration of which an SDP problem is 
solved. 
We used YALMIP with the default setting, where SDP subproblems are 
solved with SeDuMi ver.~1.3 \citep{Pol05,Stu99}. 

Consider the problem settings in 
\reffig{fig:robust.x2_y1}, \reffig{fig:robust.x3_y3}, and 
\reffig{fig:robust.x3_y2}. 
\reftab{tab:ex.0.data} lists the number of members ($m$), the number of 
degrees of freedom of displacements ($d$), and the upper bound for the 
structural volume ($\overline{c}$). 
In \reffig{fig:gs2x1robust} and \reffig{fig:gs_x3_y3}, the nodes are 
aligned on a $1\,\mathrm{m} \times 1\,\mathrm{m}$ grid. 
In \reffig{fig:gs_x3_y2}, we use a 
$1\,\mathrm{m} \times 0.5\,\mathrm{m}$ grid. 
In \reffig{fig:gs_x3_y3}, the ground structure has all possible members 
connecting two nodes but are no longer than $3\,\mathrm{m}$. 
The nominal external load, $\tilde{\bi{p}}$, is applied as shown in 
\reffig{fig:gs2x1robust}, \reffig{fig:gs_x3_y3}, and 
\reffig{fig:gs_x3_y2}. 
The uncertainty model of the external load is defined by using 
\eqref{eq.ex.uncertainty.1} with $\tilde{p}_{1}=100\,\mathrm{kN}$, 
$\alpha = 0.75 \tilde{p}_{1}$ for 
\reffig{fig:robust.x2_y1} and \reffig{fig:robust.x3_y3}, and 
$\alpha = 0.5\tilde{p}_{1}$ for \reffig{fig:robust.x3_y2}. 
The lower and upper bounds for the member cross-sectional areas are 
$\underline{x}=1\,\mathrm{mm^{2}}$ and 
$\overline{x}=700\,\mathrm{mm^{2}}$, respectively. 
In \reftab{tab:ex.0.data}, ``rob.\ opt.''\ reports the optimal value, 
obtained by YALMIP, of the robust optimization problem. 
The obtained solutions are shown in 
\reffig{fig:x2_y1_global}, \reffig{fig:x3_y3_global}, and 
\reffig{fig:x3_y2_global}. 
For reference, the optimal solutions of the (not robust) compliance 
minimization with the nominal external load, $\tilde{\bi{p}}$, are shown 
in \reffig{fig:x2_y1_nominal}, \reffig{fig:x3_y3_nominal}, 
\reffig{fig:x3_y2_nominal}.\footnote{
Ground structures without overlapping members are used to obtain the 
solutions in \reffig{fig:x2_y1_nominal}, \reffig{fig:x3_y3_nominal}, and 
\reffig{fig:x3_y2_nominal}. 
} 
The optimal values are listed in ``nom.\ opt.''\ of 
\reftab{tab:ex.0.data}. 

It is remarkable that, for every instance, the solution obtained by the 
proposed algorithm coincides wit the global optimal solution (obtained 
by YALMIP). 
\reftab{tab:ex.0.result} reports the computational costs of the two 
methods, where ``{\#}iter.''\ is the number of iterations, and 
``time'' is the required computational time. 
Note that the computational cost of the proposed method does not include 
the ones for generation of an initial point and for the post-processing. 
It is also worth noting that the problem size of SDP solved at each 
iteration of the proposed method is larger than tat of YALMIP, and hence 
the computational time per an iteration required by the proposed method 
is larger than that of YALMIP. 
For every instance, the number of SDP problems solved by the proposed 
method is smaller than that of YALMIP. 

In the experiments in this section, it has been observed that the 
proposed method converges to a global optimal solution for a small-size 
problem instance. 
In section~\ref{sec:ex.eva4} and section~\ref{sec:ex.eva3}, we examine 
large-scale instances that cannot be solved with a global optimization 
method within realistic computational time.

\subsection{Example (II)}
\label{sec:ex.eva4}

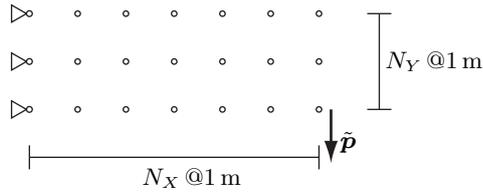
\begin{figure}[tp]
  \centering
  \includegraphics[scale=0.50]{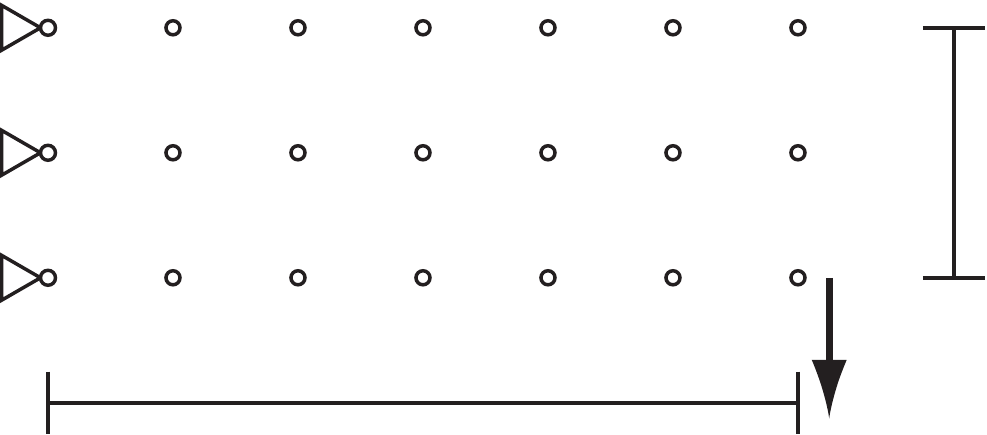}
  \begin{picture}(0,0)
    \put(-150,-50){
    \put(53,44){{\footnotesize $N_{X}${\,}@$1\,\mathrm{m}$}}
    \put(144,88){{\footnotesize $N_{Y}${\,}@$1\,\mathrm{m}$}}
    \put(127,58){{\footnotesize $\tilde{\bi{p}}$}}
    }
  \end{picture}
  \medskip
  \caption{Example (II). 
  The problem setting for $(N_{X},N_{Y})=(6,2)$. }
  \label{fig:gs6x2_bottom}
\end{figure}

\begin{figure}[tp]
  %%%% C:\doc\robust\topology\load_ccp\eva4\opt_design.m
  \centering
  \begin{subfigure}[b]{0.15\textwidth}
    \centering
    \includegraphics[scale=0.42]{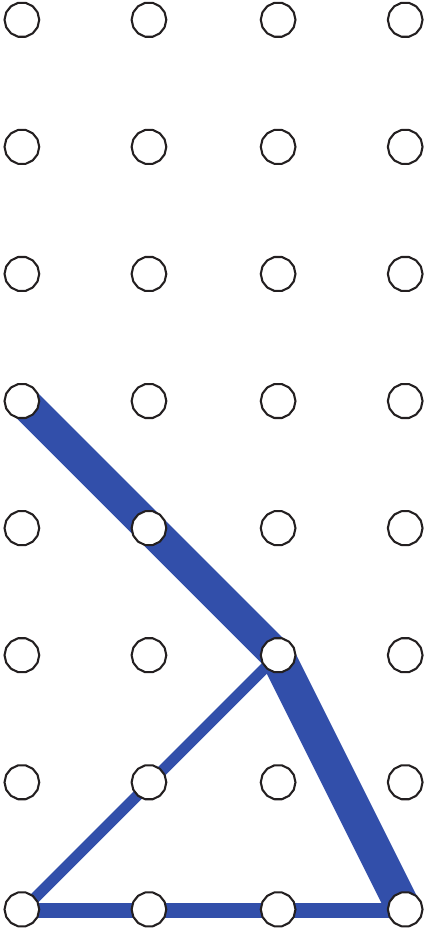}
    \caption{}
    \label{fig:x3_y7_nominal_bottom}
  \end{subfigure}
  \hfill
  \begin{subfigure}[b]{0.20\textwidth}
    \centering
    \includegraphics[scale=0.36]{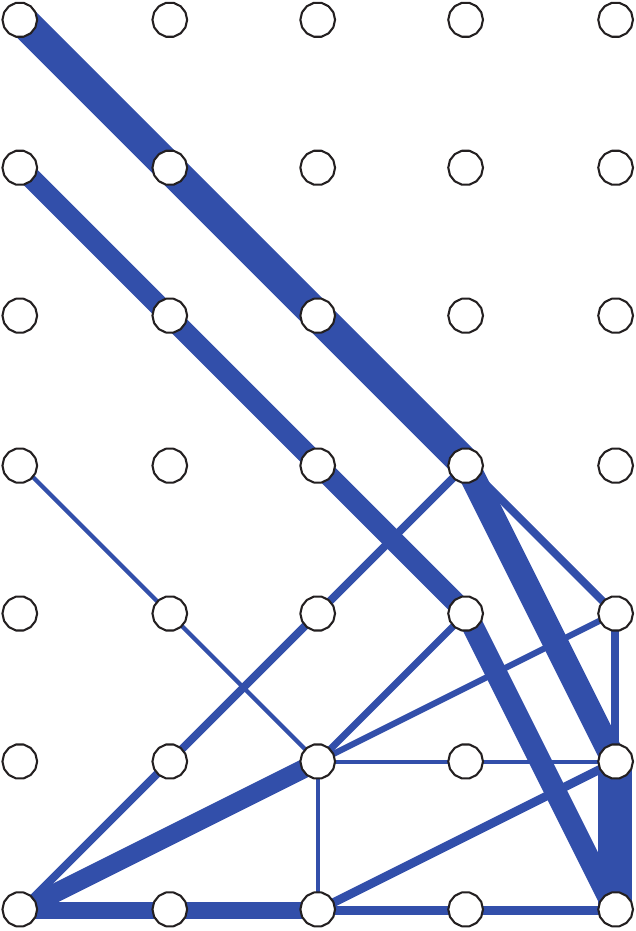}
    \caption{}
    \label{fig:x4_y6_nominal_bottom}
  \end{subfigure}
  \hfill
  \begin{subfigure}[b]{0.25\textwidth}
    \centering
    \includegraphics[scale=0.30]{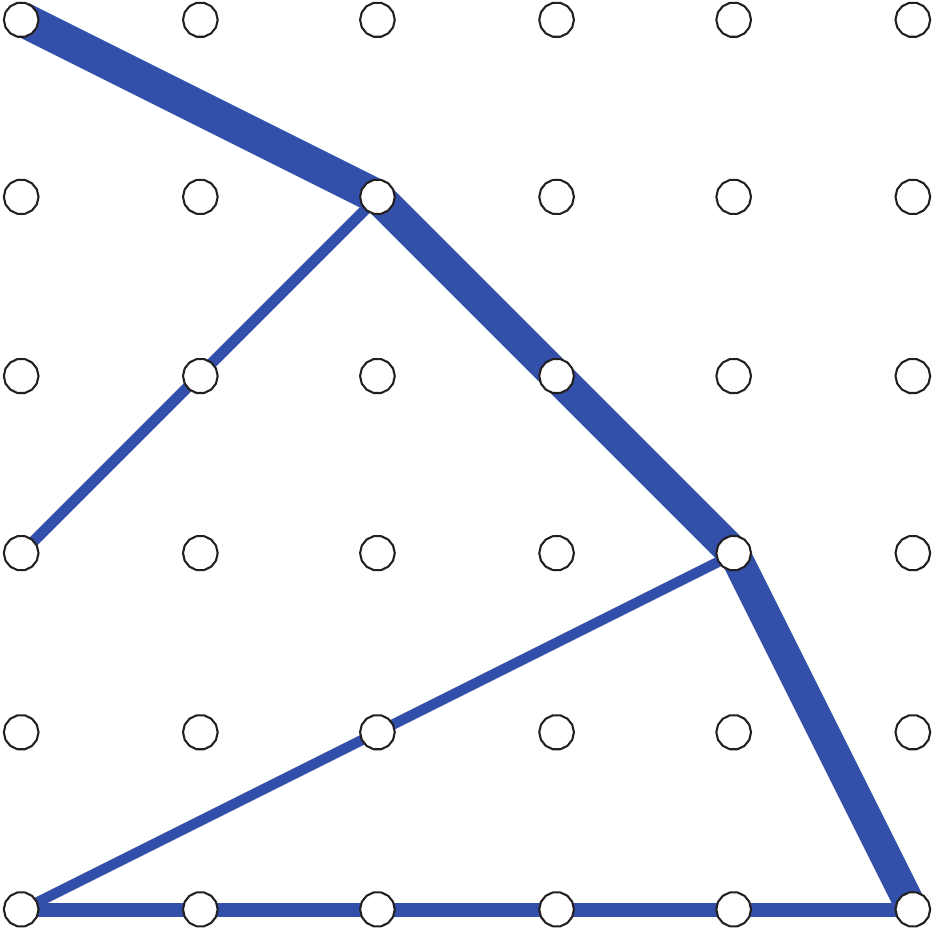}
    \caption{}
    \label{fig:x5_y5_nominal_bottom}
  \end{subfigure}
  \hfill
  \begin{subfigure}[b]{0.30\textwidth}
    \centering
    \includegraphics[scale=0.30]{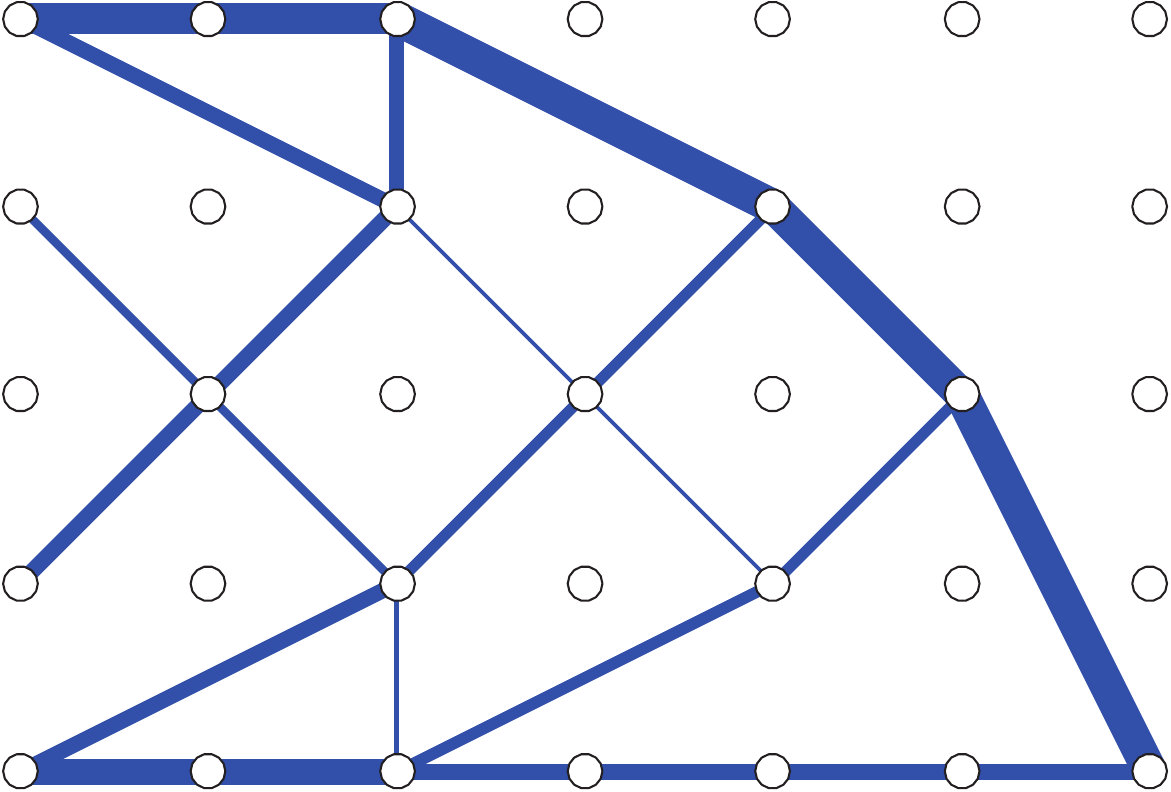}
    \caption{}
    \label{fig:x6_y4_nominal_bottom}
  \end{subfigure}
  \par\medskip
  \begin{subfigure}[b]{0.35\textwidth}
    \centering
    \includegraphics[scale=0.36]{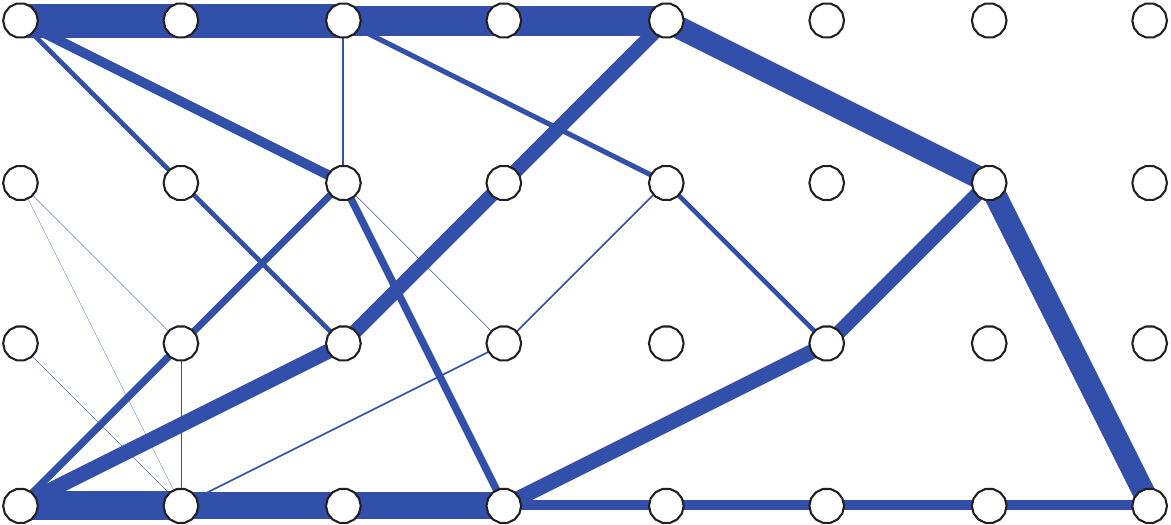}
    \caption{}
    \label{fig:x7_y3_nominal_bottom}
  \end{subfigure}
  \hfill
  \begin{subfigure}[b]{0.40\textwidth}
    \centering
    \includegraphics[scale=0.42]{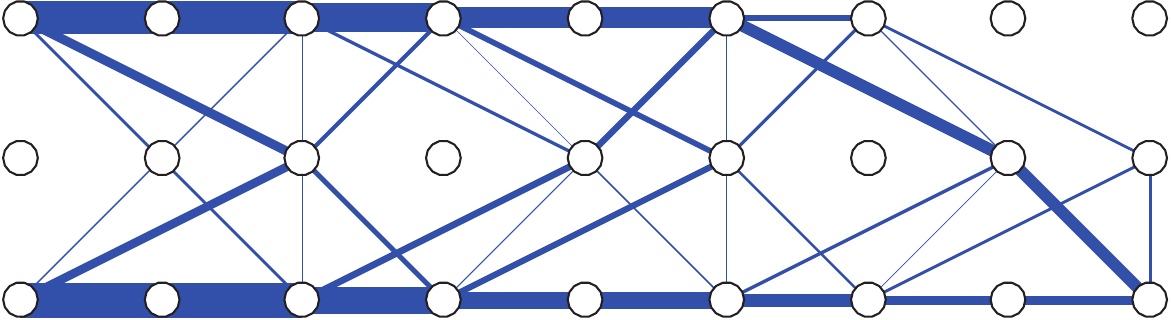}
    \caption{}
    \label{fig:x8_y2_nominal_bottom}
  \end{subfigure}
  \caption{Example (II). 
  The optimal solutions of the compliance minimization for the nominal 
  external load. 
  \subref{fig:x3_y7_nominal_bottom} $(N_{X},N_{Y})=(3,7)$; 
  \subref{fig:x4_y6_nominal_bottom} $(4,6)$; 
  \subref{fig:x5_y5_nominal_bottom} $(5,5)$; 
  \subref{fig:x6_y4_nominal_bottom} $(6,4)$; 
  \subref{fig:x7_y3_nominal_bottom} $(7,3)$; and 
  \subref{fig:x8_y2_nominal_bottom} $(8,2)$. 
  }
  \label{fig:bottom_load_nominal}
%\end{figure}
%
  \bigskip
%\begin{figure}[tp]
  %%%% C:\doc\robust\topology\load_ccp\eva4\opt_design.m
  \centering
  \begin{subfigure}[b]{0.15\textwidth}
    \centering
    \includegraphics[scale=0.42]{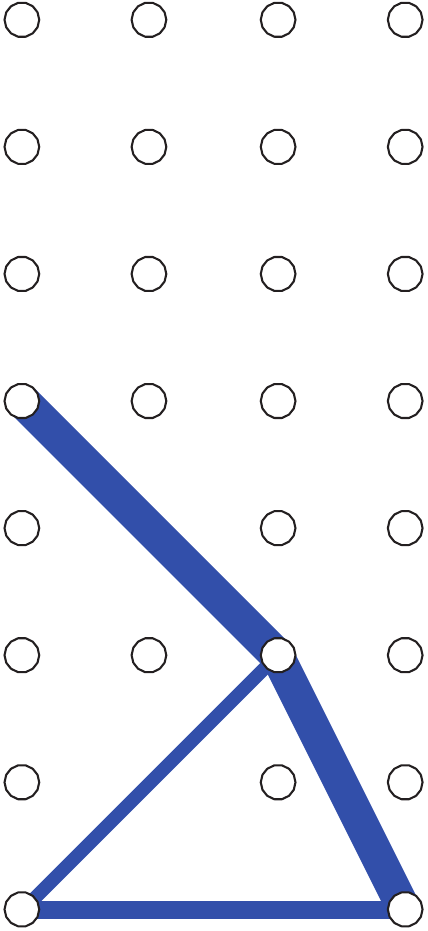}
    \caption{}
    \label{fig:x3_y7_robust_bottom}
  \end{subfigure}
  \hfill
  \begin{subfigure}[b]{0.20\textwidth}
    \centering
    \includegraphics[scale=0.36]{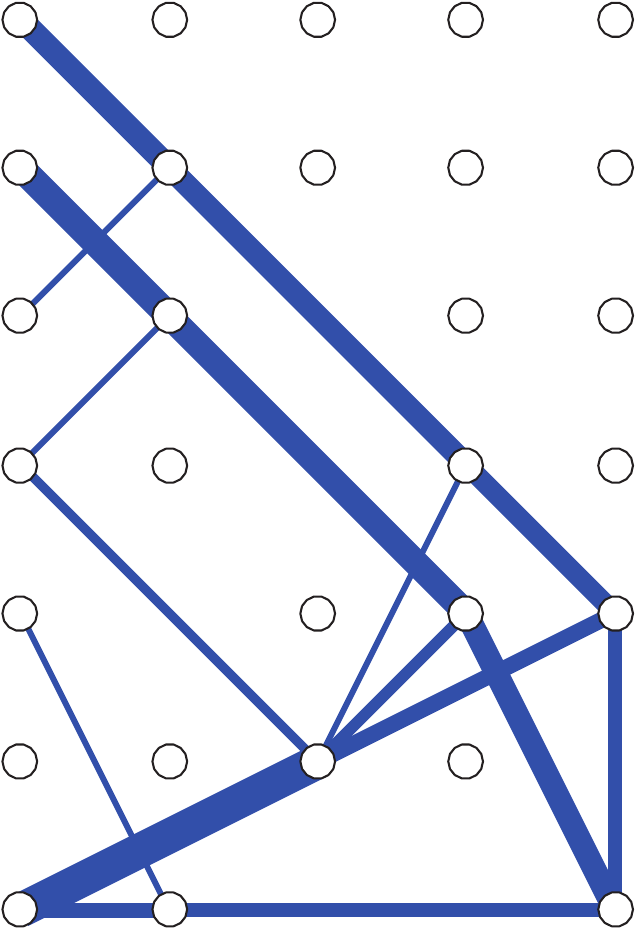}
    \caption{}
    \label{fig:x4_y6_robust_bottom}
  \end{subfigure}
  \hfill
  \begin{subfigure}[b]{0.25\textwidth}
    \centering
    \includegraphics[scale=0.30]{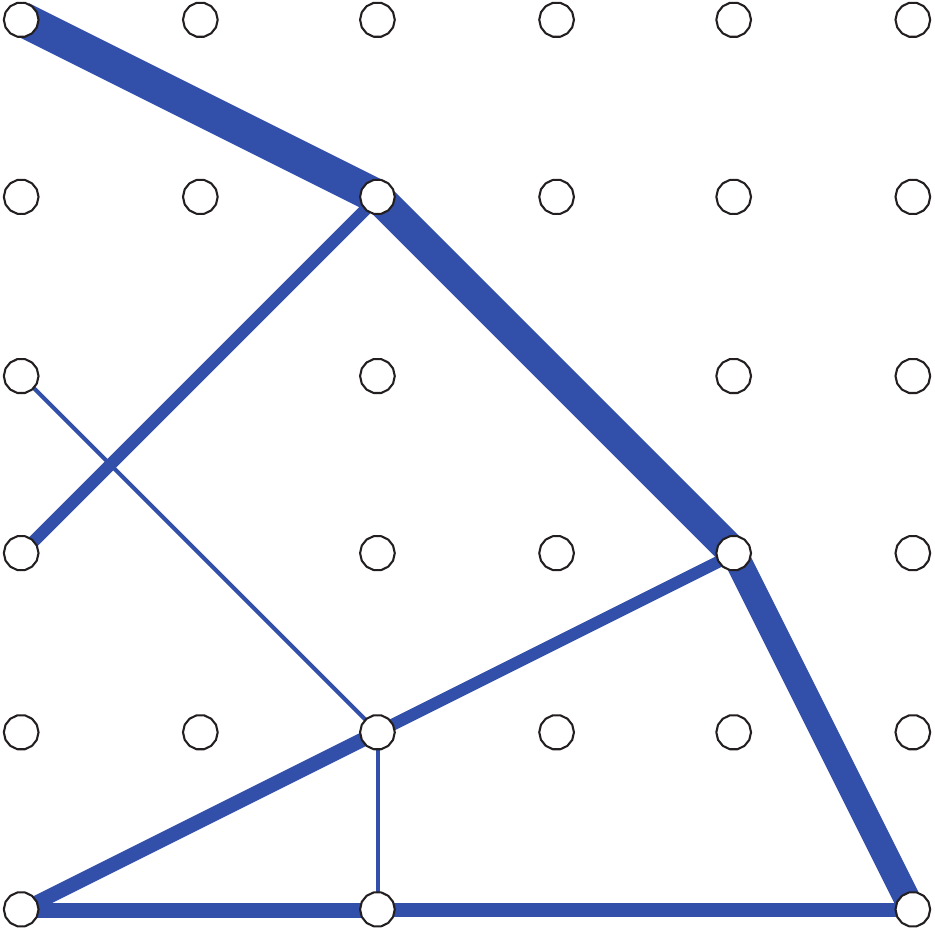}
    \caption{}
    \label{fig:x5_y5_robust_bottom}
  \end{subfigure}
  \hfill
  \begin{subfigure}[b]{0.30\textwidth}
    \centering
    \includegraphics[scale=0.30]{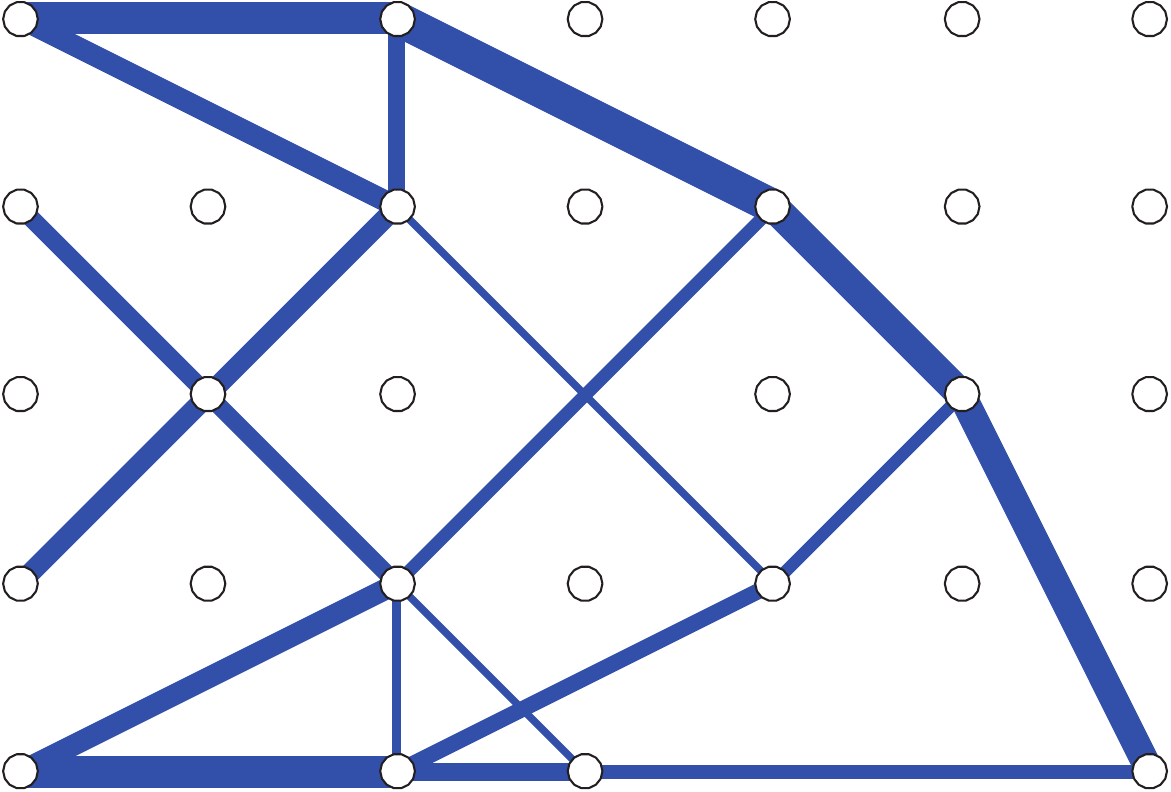}
    \caption{}
    \label{fig:x6_y4_robust_bottom}
  \end{subfigure}
  \par\medskip
  \begin{subfigure}[b]{0.35\textwidth}
    \centering
    \includegraphics[scale=0.36]{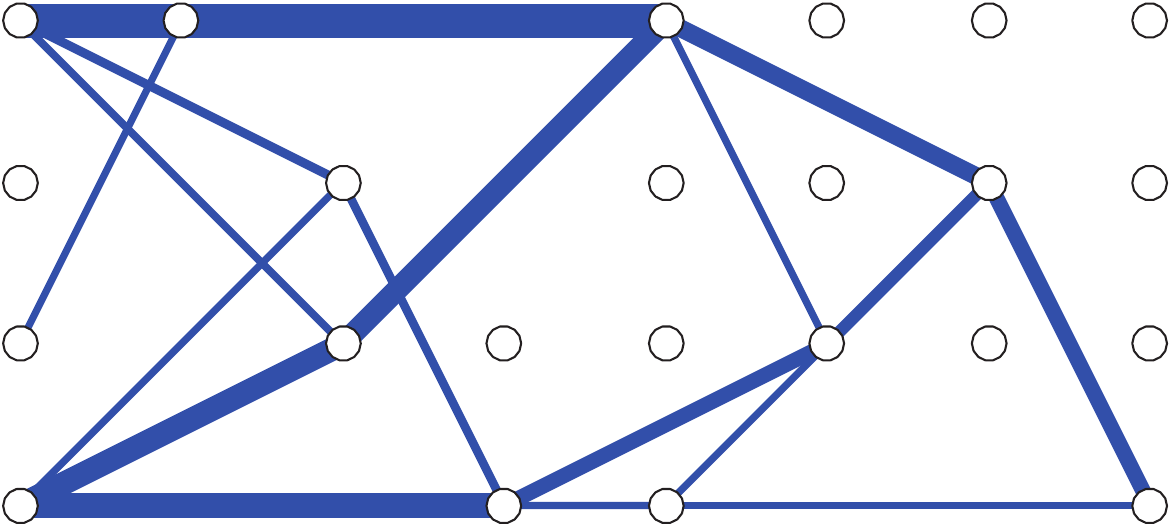}
    \caption{}
    \label{fig:x7_y3_robust_bottom}
  \end{subfigure}
  \hfill
  \begin{subfigure}[b]{0.40\textwidth}
    \centering
    \includegraphics[scale=0.42]{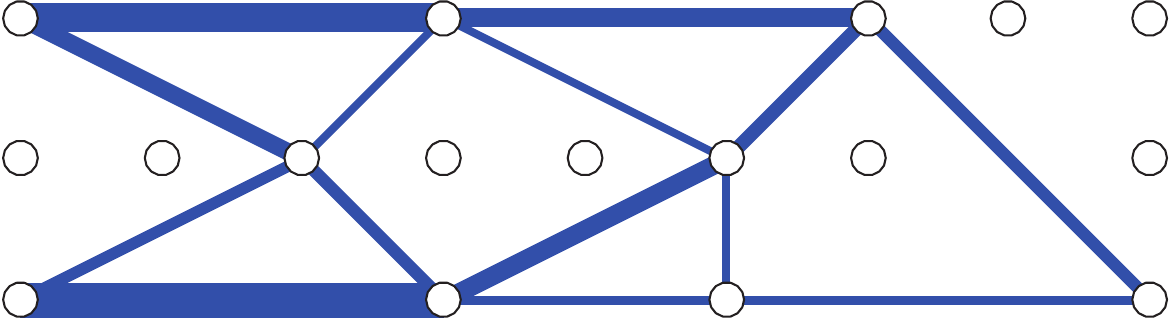}
    \caption{}
    \label{fig:x8_y2_robust_bottom}
  \end{subfigure}
  \caption{Example (II). 
  The solutions obtained by the proposed method for the robust 
  optimization under the load uncertainty. 
  \subref{fig:x3_y7_robust_bottom} $(N_{X},N_{Y})=(3,7)$; 
  \subref{fig:x4_y6_robust_bottom} $(4,6)$; 
  \subref{fig:x5_y5_robust_bottom} $(5,5)$; 
  \subref{fig:x6_y4_robust_bottom} $(6,4)$; 
  \subref{fig:x7_y3_robust_bottom} $(7,3)$; and 
  \subref{fig:x8_y2_robust_bottom} $(8,2)$. 
  }
  \label{fig:bottom_load_robust}
\end{figure}

\begin{figure}[tp]
  %%%% C:\doc\robust\topology\load_ccp\eva4\opt_design.m
  \centering
  \begin{subfigure}[b]{0.15\textwidth}
    \centering
    \includegraphics[scale=0.42]{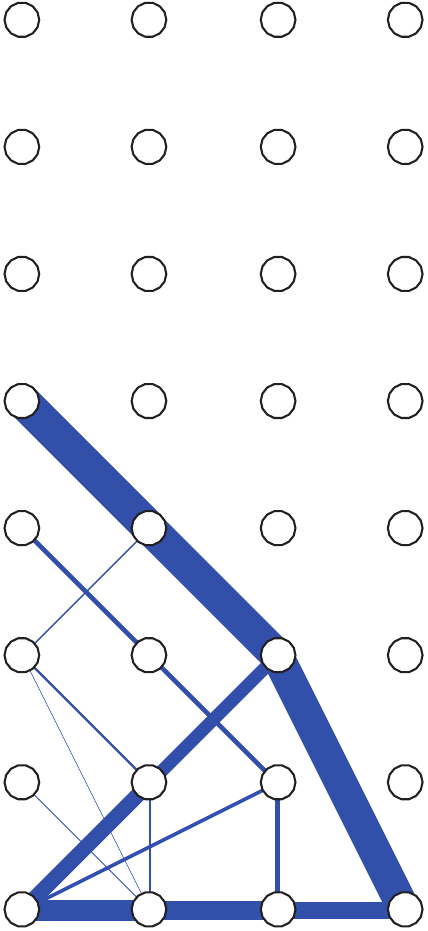}
    \caption{}
    \label{fig:x3_y7_specified_bottom}
  \end{subfigure}
  \hfill
  \begin{subfigure}[b]{0.20\textwidth}
    \centering
    \includegraphics[scale=0.36]{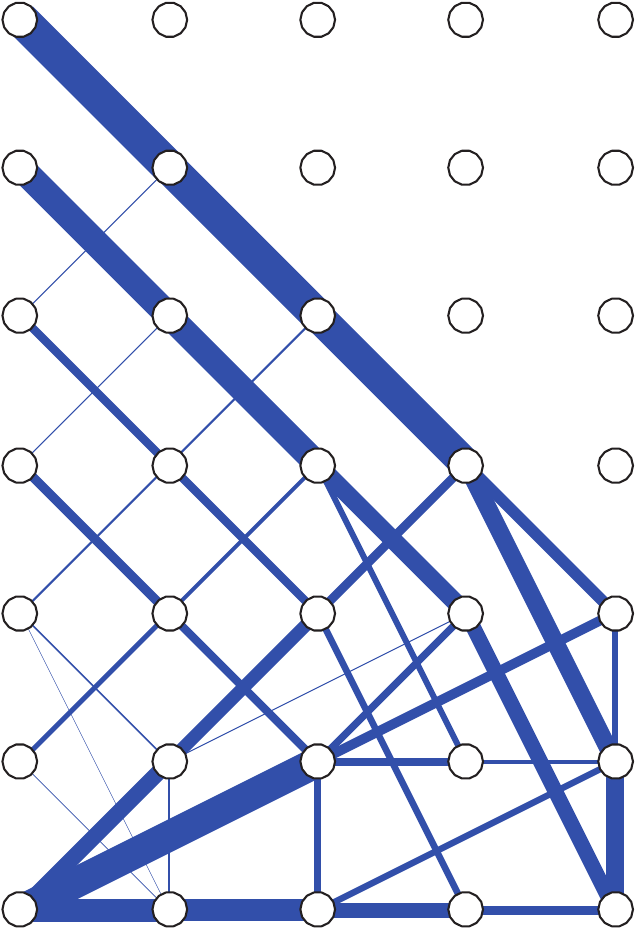}
    \caption{}
    \label{fig:x4_y6_specified_bottom}
  \end{subfigure}
  \hfill
  \begin{subfigure}[b]{0.25\textwidth}
    \centering
    \includegraphics[scale=0.30]{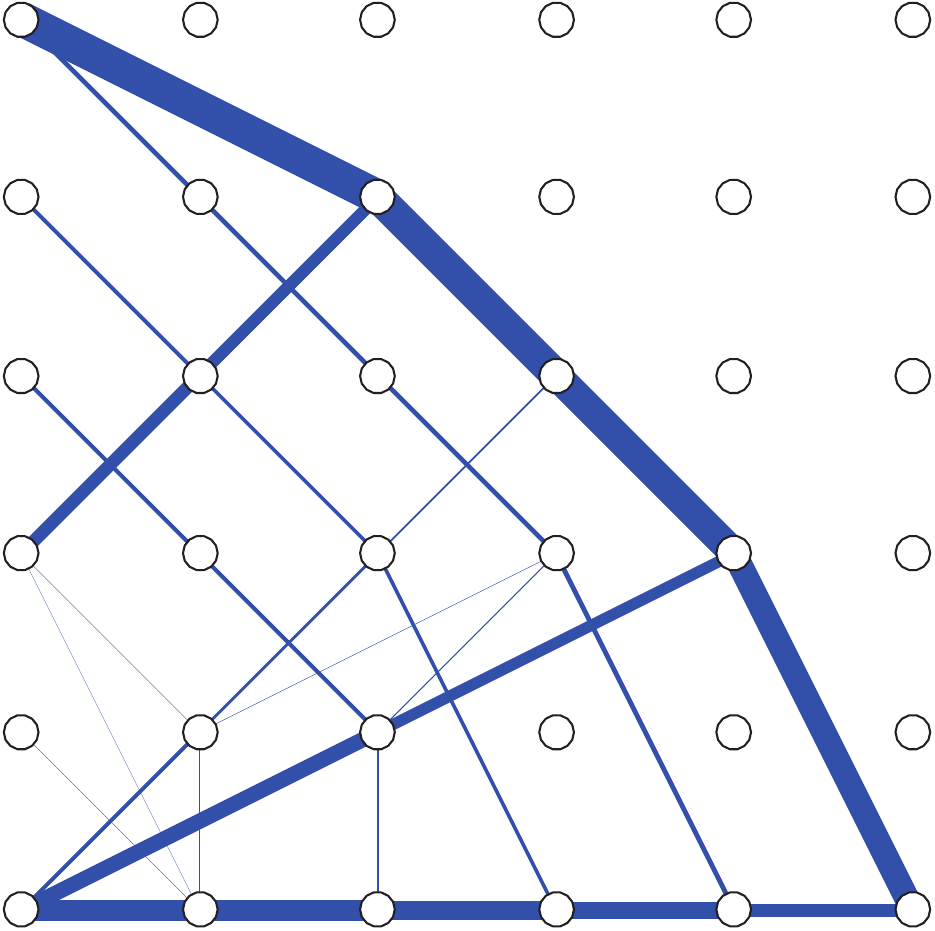}
    \caption{}
    \label{fig:x5_y5_specified_bottom}
  \end{subfigure}
  \hfill
  \begin{subfigure}[b]{0.30\textwidth}
    \centering
    \includegraphics[scale=0.30]{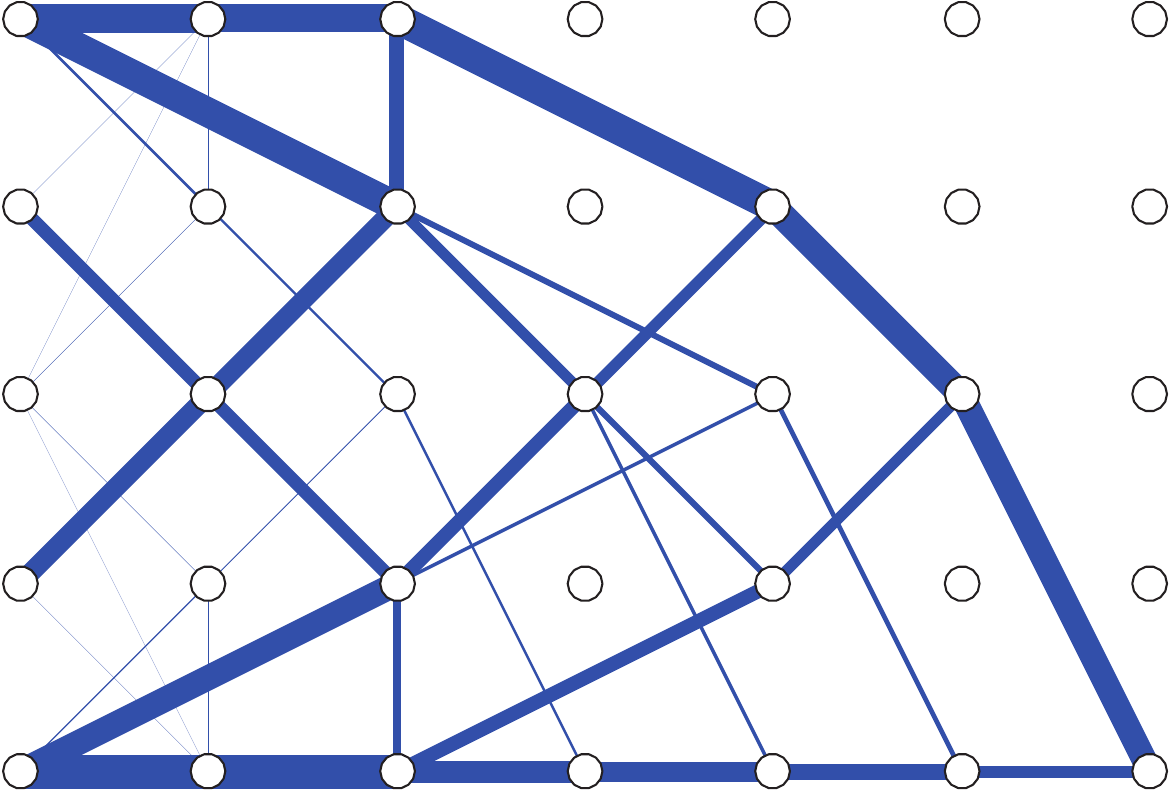}
    \caption{}
    \label{fig:x6_y4_specified_bottom}
  \end{subfigure}
  \par\medskip
  \begin{subfigure}[b]{0.35\textwidth}
    \centering
    \includegraphics[scale=0.36]{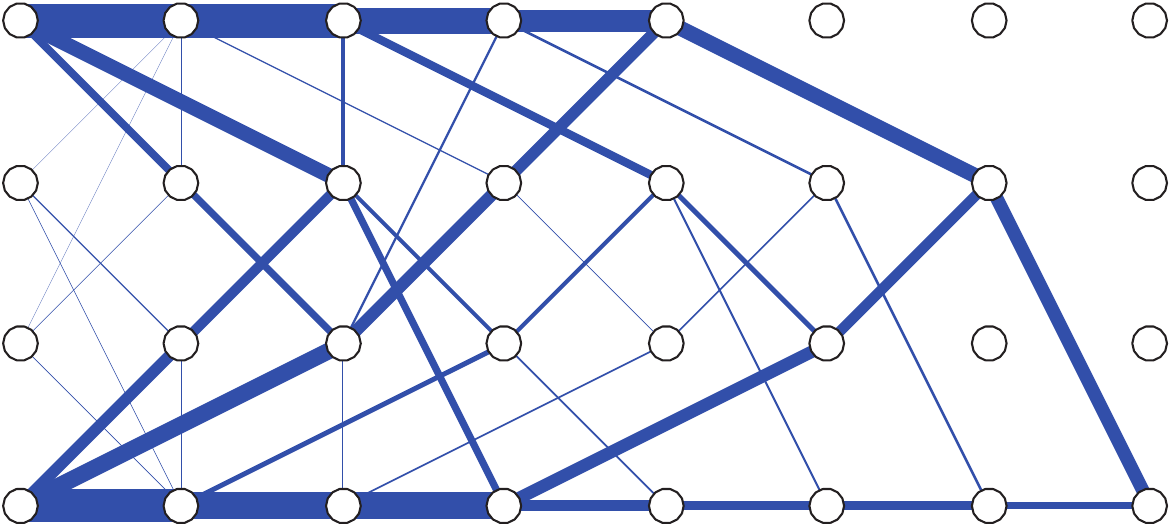}
    \caption{}
    \label{fig:x7_y3_specified_bottom}
  \end{subfigure}
  \hfill
  \begin{subfigure}[b]{0.40\textwidth}
    \centering
    \includegraphics[scale=0.42]{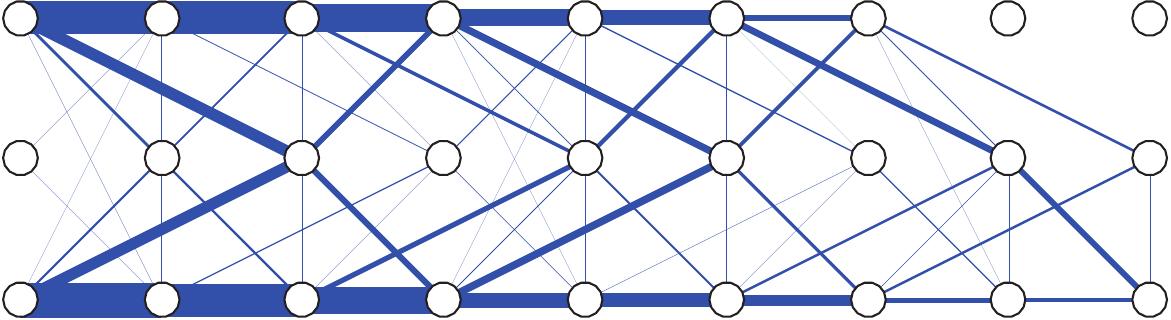}
    \caption{}
    \label{fig:x8_y2_specified_bottom}
  \end{subfigure}
  \caption{The optimal solution of problem \eqref{P.fixed.s}, where 
  $\bar{\bi{s}}$ is defined from the solutions in 
  \reffig{fig:bottom_load_nominal}. }
  \label{fig:bottom_load_specified}
\end{figure}

\newcommand{\ThreeSevenBottomM}{\ensuremath{250}}
\newcommand{\ThreeSevenBottomD}{\ensuremath{48}}
\newcommand{\ThreeSevenBottomIter}{\ensuremath{9}}
\newcommand{\ThreeSevenBottomTime}{\ensuremath{46.4}}
\newcommand{\ThreeSevenBottomObj}{\ensuremath{836.310}}
\newcommand{\ThreeSevenBottomTilde}{\ensuremath{766.518}}
\newcommand{\ThreeSevenBottomNominal}{\ensuremath{761.905}}
\newcommand{\ThreeSevenBottomRatio}{\ensuremath{1.006}}
\newcommand{\ThreeSevenBottomSpecified}{\ensuremath{986.442}}
\newcommand{\ThreeSevenBottomVolume}{\ensuremath{4.2 \times 10^{6}}}

\newcommand{\FourSixBottomM}{\ensuremath{292}}
\newcommand{\FourSixBottomD}{\ensuremath{56}}
\newcommand{\FourSixBottomIter}{\ensuremath{39}}
\newcommand{\FourSixBottomTime}{\ensuremath{242.9}}
\newcommand{\FourSixBottomObj}{\ensuremath{1807.714}}
\newcommand{\FourSixBottomTilde}{\ensuremath{1360.425}}
\newcommand{\FourSixBottomNominal}{\ensuremath{1185.185}}
\newcommand{\FourSixBottomRatio}{\ensuremath{1.148}}
\newcommand{\FourSixBottomSpecified}{\ensuremath{2534.505}}
\newcommand{\FourSixBottomVolume}{\ensuremath{4.8 \times 10^{6}}}

\newcommand{\FiveFiveBottomM}{\ensuremath{306}}
\newcommand{\FiveFiveBottomD}{\ensuremath{60}}
\newcommand{\FiveFiveBottomIter}{\ensuremath{35}}
\newcommand{\FiveFiveBottomTime}{\ensuremath{249.6}}
\newcommand{\FiveFiveBottomObj}{\ensuremath{2382.377}}
\newcommand{\FiveFiveBottomTilde}{\ensuremath{2034.270}}
\newcommand{\FiveFiveBottomNominal}{\ensuremath{1929.012}}
\newcommand{\FiveFiveBottomRatio}{\ensuremath{1.055}}
\newcommand{\FiveFiveBottomSpecified}{\ensuremath{3017.593}}
\newcommand{\FiveFiveBottomVolume}{\ensuremath{5.0 \times 10^{6}}}

\newcommand{\SixFourBottomM}{\ensuremath{292}}
\newcommand{\SixFourBottomD}{\ensuremath{60}}
\newcommand{\SixFourBottomIter}{\ensuremath{21}}
\newcommand{\SixFourBottomTime}{\ensuremath{142.1}}
\newcommand{\SixFourBottomObj}{\ensuremath{5913.978}}
\newcommand{\SixFourBottomTilde}{\ensuremath{4427.633}}
\newcommand{\SixFourBottomNominal}{\ensuremath{4143.551}}
\newcommand{\SixFourBottomRatio}{\ensuremath{1.069}}
\newcommand{\SixFourBottomSpecified}{\ensuremath{7032.673}}
\newcommand{\SixFourBottomVolume}{\ensuremath{4.8 \times 10^{6}}}

\newcommand{\SevenThreeBottomM}{\ensuremath{250}}
\newcommand{\SevenThreeBottomD}{\ensuremath{56}}
\newcommand{\SevenThreeBottomIter}{\ensuremath{40}}
\newcommand{\SevenThreeBottomTime}{\ensuremath{193.9}}
\newcommand{\SevenThreeBottomObj}{\ensuremath{14912.232}}
\newcommand{\SevenThreeBottomTilde}{\ensuremath{10960.621}}
\newcommand{\SevenThreeBottomNominal}{\ensuremath{9918.356}}
\newcommand{\SevenThreeBottomRatio}{\ensuremath{1.105}}
\newcommand{\SevenThreeBottomSpecified}{\ensuremath{17717.408}}
\newcommand{\SevenThreeBottomVolume}{\ensuremath{4.2 \times 10^{6}}}

\newcommand{\EightTwoBottomM}{\ensuremath{180}}
\newcommand{\EightTwoBottomD}{\ensuremath{48}}
\newcommand{\EightTwoBottomIter}{\ensuremath{32}}
\newcommand{\EightTwoBottomTime}{\ensuremath{103.9}}
\newcommand{\EightTwoBottomObj}{\ensuremath{43467.983}}
\newcommand{\EightTwoBottomTilde}{\ensuremath{34515.627}}
\newcommand{\EightTwoBottomNominal}{\ensuremath{34515.626}}
\newcommand{\EightTwoBottomRatio}{\ensuremath{1.000}}
\newcommand{\EightTwoBottomSpecified}{\ensuremath{71121.097}}
\newcommand{\EightTwoBottomVolume}{\ensuremath{3.2 \times 10^{6}}}

\begin{table}[bp]
  %%%% C:\doc\robust\topology\load_ccp\eva4\opt_design.m
  \centering
  \caption{Characteristics of the problem instances in example (II).}
  \label{tab:ex.I.data}
  \begin{tabular}{lrrr}
    \toprule
    $(N_{X},N_{Y})$ & $m$ & $d$ & $\overline{c}$ $(\mathrm{mm}^{3})$ \\
    \midrule
    (3,7) 
    & \ThreeSevenBottomM & \ThreeSevenBottomD & \ThreeSevenBottomVolume \\
    (4,6) 
    & \FourSixBottomM & \FourSixBottomD & \FourSixBottomVolume \\
    (5,5) 
    & \FiveFiveBottomM & \FiveFiveBottomD & \FiveFiveBottomVolume \\
    (6,4) 
    & \SixFourBottomM & \SixFourBottomD & \SixFourBottomVolume \\
    (7,3) 
    & \SevenThreeBottomM & \SevenThreeBottomD & \SevenThreeBottomVolume \\
    (8,2) 
    & \EightTwoBottomM & \EightTwoBottomD & \EightTwoBottomVolume \\
    \bottomrule
  \end{tabular}
%\end{table}
%
%\begin{table}[bp]
  %%%% C:\doc\robust\topology\load_ccp\eva4\opt_design.m
  \centering
  \caption{Computational results of example (II).}
  \label{tab:ex.I.result}
  \begin{tabular}{lrrrrrr}
    \toprule
    $(N_{X},N_{Y})$ & Obj.\ (J) & {\#}iter.\ 
    & Time (s) & $\tilde{w}$ (J) & Nom.\ opt.\ (J) & Fixed $\bi{s}$ (J) \\
    \midrule
    (3,7) 
    & \ThreeSevenBottomObj & \ThreeSevenBottomIter & \ThreeSevenBottomTime & \ThreeSevenBottomTilde & \ThreeSevenBottomNominal & \ThreeSevenBottomSpecified \\
    (4,6)
    & \FourSixBottomObj & \FourSixBottomIter & \FourSixBottomTime & \FourSixBottomTilde & \FourSixBottomNominal & \FourSixBottomSpecified \\
    (5,5)
    & \FiveFiveBottomObj & \FiveFiveBottomIter & \FiveFiveBottomTime & \FiveFiveBottomTilde & \FiveFiveBottomNominal & \FiveFiveBottomSpecified \\
    (6,4)
    & \SixFourBottomObj & \SixFourBottomIter & \SixFourBottomTime & \SixFourBottomTilde & \SixFourBottomNominal & \SixFourBottomSpecified \\
    (7,3)
    & \SevenThreeBottomObj & \SevenThreeBottomIter & \SevenThreeBottomTime & \SevenThreeBottomTilde & \SevenThreeBottomNominal & \SevenThreeBottomSpecified \\
    (8,2)
    & \EightTwoBottomObj & \EightTwoBottomIter & \EightTwoBottomTime & \EightTwoBottomTilde & \EightTwoBottomNominal & \EightTwoBottomSpecified \\
    \bottomrule
  \end{tabular}
\end{table}

Consider the ground structures shown in \reffig{fig:gs6x2_bottom}. 
The nodes are aligned on a $1\,\mathrm{m} \times 1\,\mathrm{m}$ grid, 
and the number of the nodes is $(N_{X}+1)(N_{Y}+1)$. 
The leftmost nodes are pin-supported. 
The candidate members are defined as follows. 
We first generate all possible members such that any two nodes are 
connected by a member. 
Then we remove members that are longer than $3\,\mathrm{m}$. 
It should be clear that the ground structure retains overlapping members. 

As for the nominal external load, $\tilde{\bi{p}}$, a vertical force is 
applied to the bottom rightmost node as shown in 
\reffig{fig:gs6x2_bottom}. 
The uncertainty model of the external load is given as explained in 
\eqref{eq.ex.uncertainty.1}, where 
$\tilde{p}_{1}=100\,\mathrm{kN}$ and $\alpha= 0.5 \tilde{p}_{1}$. 
The lower and upper bounds for the member cross-sectional areas are 
$\underline{x}=50\,\mathrm{mm^{2}}$ and 
$\overline{x}=700\,\mathrm{mm^{2}}$. 
The upper bound for the structural volume is 
$\overline{c}=2N_{X}N_{Y} \times 10^{5}\,\mathrm{mm}^{3}$. 

As for problem sizes, we consider six cases, 
$(N_{X},N_{Y})=(3,7)$, $(4,6)$, $(5,5)$, $(6,4)$, $(7,3)$, and $(8,2)$. 
\reftab{tab:ex.I.data} lists the number of members, the number of 
degrees of freedom of the nodal displacements, and the upper bound for 
the structural volume. 
\reffig{fig:bottom_load_nominal} collects the optimal solutions of the 
conventional (i.e., not robust) compliance minimization in 
\eqref{P.nominal.compliance} for the nominal 
external load, $\tilde{\bi{p}}$. 
For the robust optimization, the solutions obtained by the proposed 
method (\refalg{alg:concave.convex.penalty}) are shown in 
\reffig{fig:bottom_load_robust}. 
It is observed in \reffig{fig:x3_y7_robust_bottom} that the three chains 
in \reffig{fig:x3_y7_nominal_bottom} are replaced with longer single 
members and four intermediate nodes are removed. 
The length of the bottom horizontal chain in 
\reffig{fig:x5_y5_nominal_bottom} is $5\,\mathrm{m}$. 
This chain is replaced with two members in 
\reffig{fig:x5_y5_robust_bottom}, because the maximum member length in 
the ground structure is $3\,\mathrm{m}$. 
Similar observation can be made also in 
\reffig{fig:x6_y4_robust_bottom} and 
\reffig{fig:x7_y3_robust_bottom}. 
The nominal optimal solution in \reffig{fig:x8_y2_nominal_bottom} has 
many thin members as well as many nodes. 
In contrast, the robust solution in \reffig{fig:x8_y2_robust_bottom} has 
simple topology, which may be considered practically preferable. 
In all the solutions in \reffig{fig:bottom_load_robust}, the longest 
member is no longer than the maximum member length in the ground 
structure, as explained in section~\ref{sec:preliminary.robust}. 

The computational results are listed in \reftab{tab:ex.I.result}. 
Here, ``obj.''\ reports the objective value of the solution obtained by 
the proposed method, and $\tilde{w}$ is the compliance of this solution 
for the nominal external load, $\tilde{\bi{p}}$. 
The optimal value of problem \eqref{P.nominal.compliance} for 
$\tilde{\bi{p}}$ is listed as ``nom.\ opt.'' 
Therefore, $\tilde{w}$ is no smaller than the value of 
``nom.\ opt.'' 
It is observed in \reftab{tab:ex.I.result} that these two values are 
very close. 
Namely, in these examples, robustness can be achieved with 
compensation of only very small increase of the nominal compliance. 

As explained in section~\ref{sec:motivation}, one of difficulties of the 
robust truss topology optimization is that the uncertainty model of the 
external load depends on the existing nodes, i.e., on $\bi{s}$. 
For comparison, we fix $\bi{s}$ to obtain a robust solution. 
As a simple heuristic, we construct an estimate of $\bi{s}$, denoted 
$\bar{\bi{s}}$, from the existing nodes of a solution in 
\reffig{fig:bottom_load_nominal}. 
Then we solve the following robust optimization problem: 
\begin{subequations}\label{P.fixed.s}%
  \begin{alignat*}{3}
    & \MIN  &{\quad}& w \\
    & \ST && 
    \sup\{ \pi(\bi{x};\bi{p}) \mid \bi{p} \in P(\bar{\bi{s}}) \} \\
    & &&
    \bi{x} \ge \bi{0} , \\
    & &&
    \bi{c}^{\top} \bi{x} \le \overline{c} , 
  \end{alignat*}
\end{subequations}
which can be recast as SDP. 
For simplicity, the lower bound constraints on the cross-sectional areas 
of the existing members are omitted. 
The obtained solutions are shown in \reffig{fig:bottom_load_specified}. 
The optimal value is reported in ``fixed $\bi{s}$'' of 
\reftab{tab:ex.I.result}. 
It should be clear that, at each solution in 
\reffig{fig:bottom_load_specified}, uncertain external forces are 
applied only to the nodes that the corresponding solution in 
\reffig{fig:bottom_load_nominal} has. 
Nevertheless, the objective value of a solution in 
\reffig{fig:bottom_load_specified} is much larger than that of the 
corresponding solution in \reffig{fig:bottom_load_robust}. 
In other words, 
the solution obtained by the proposed method has quite high quality.

\subsection{Example (III)}
\label{sec:ex.eva3}

\begin{figure}[tp]
  \centering
  \includegraphics[scale=0.50]{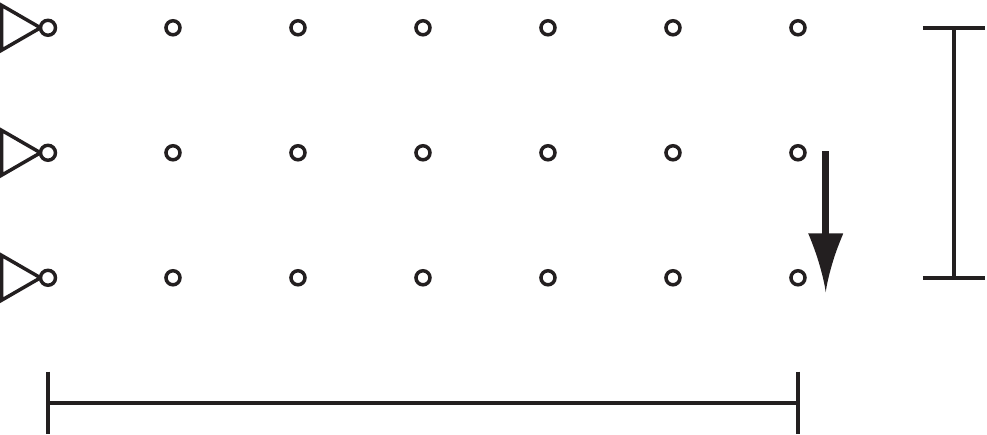}
  \begin{picture}(0,0)
    \put(-150,-50){
    \put(53,44){{\footnotesize $N_{X}${\,}@$1\,\mathrm{m}$}}
    \put(144,88){{\footnotesize $N_{Y}${\,}@$1\,\mathrm{m}$}}
    \put(127,78){{\footnotesize $\tilde{\bi{p}}$}}
    }
  \end{picture}
  \bigskip
  \caption{Example (III). 
  The problem setting for $(N_{X},N_{Y})=(6,2)$. }
  \label{fig:gs6x2}
\end{figure}
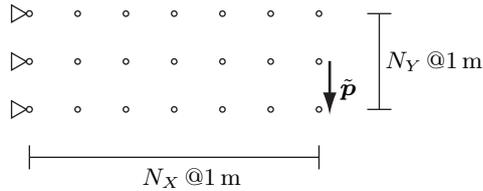

\begin{figure}[tp]
  %%%% C:\doc\robust\topology\load_ccp\eva3\opt_design.m
  \centering
  \begin{subfigure}[b]{0.275\textwidth}
    \centering
    \includegraphics[scale=0.30]{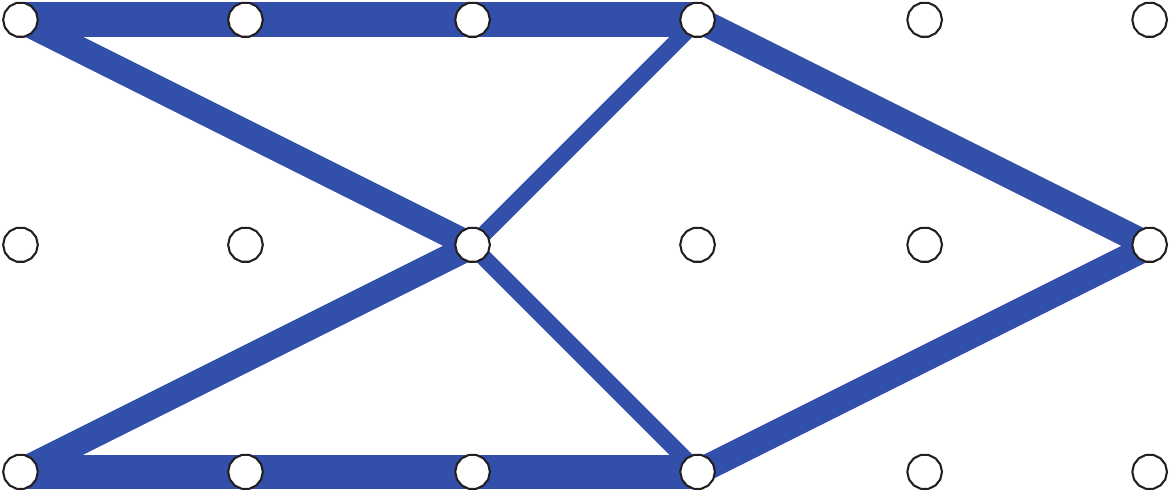}
    \caption{}
    \label{fig:x5_y2_nominal}
  \end{subfigure}
  \hfill
  \begin{subfigure}[b]{0.32\textwidth}
    \centering
    \includegraphics[scale=0.36]{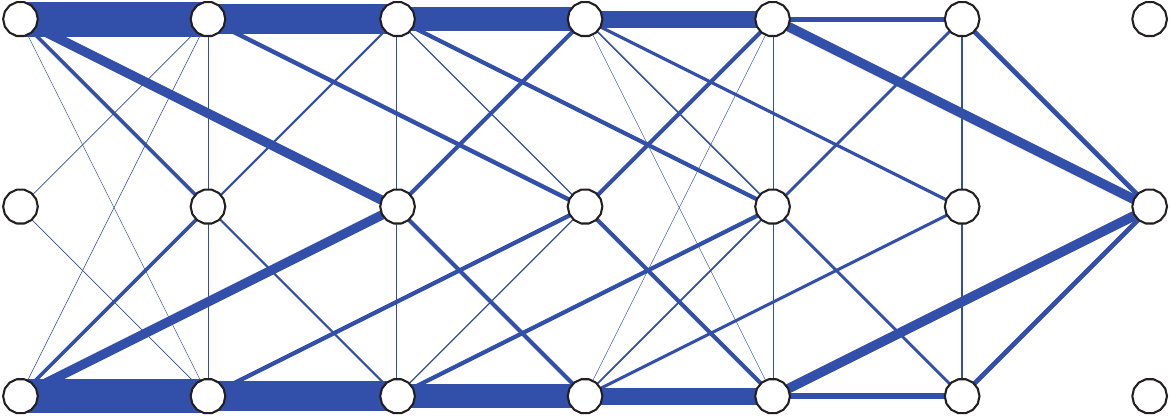}
    \caption{}
    \label{fig:x6_y2_nominal}
  \end{subfigure}
  \hfill
  \begin{subfigure}[b]{0.37\textwidth}
    \centering
    \includegraphics[scale=0.42]{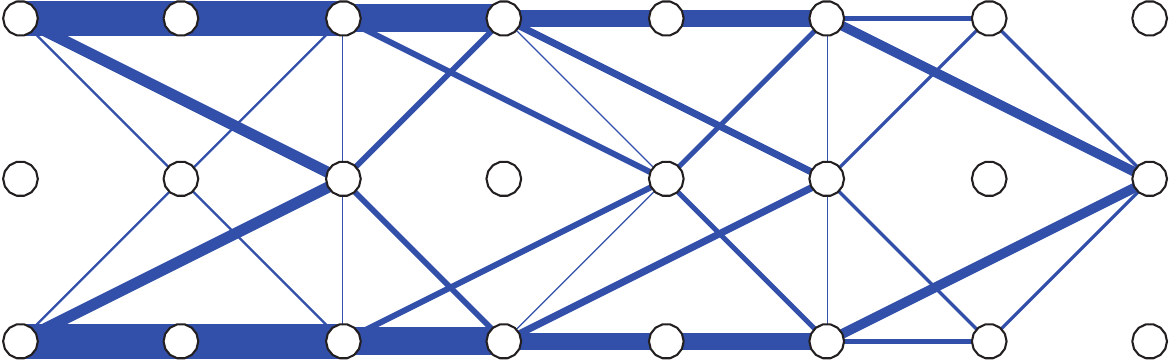}
    \caption{}
    \label{fig:x7_y2_nominal}
  \end{subfigure}
  \par\medskip
  \begin{subfigure}[b]{0.42\textwidth}
    \centering
    \includegraphics[scale=0.48]{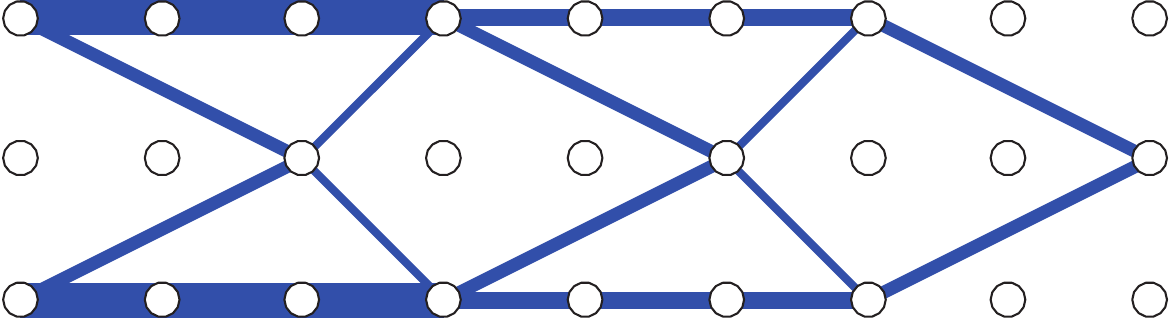}
    \caption{}
    \label{fig:x8_y2_nominal}
  \end{subfigure}
  \hfill
  \begin{subfigure}[b]{0.47\textwidth}
    \centering
    \includegraphics[scale=0.54]{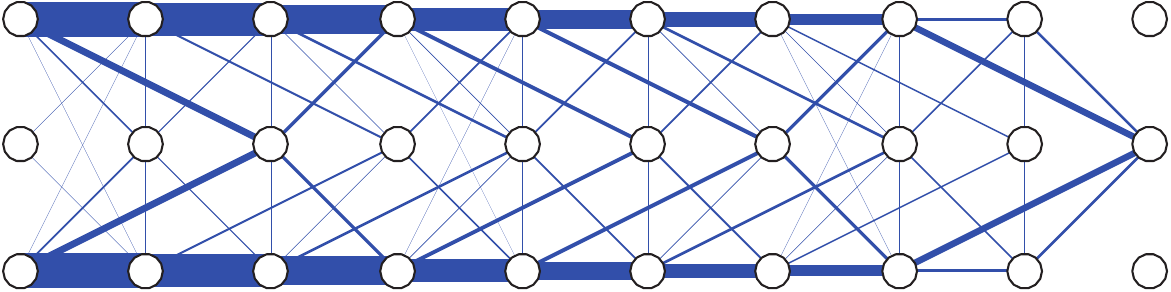}
    \caption{}
    \label{fig:x9_y2_nominal}
  \end{subfigure}
  \caption{Example (III)-1. 
  The optimal solutions of the compliance minimization for the nominal 
  external load. 
  \subref{fig:x5_y2_nominal} $(N_{X},N_{Y})=(5,2)$; 
  \subref{fig:x6_y2_nominal} $(6,2)$; 
  \subref{fig:x7_y2_nominal} $(7,2)$; 
  \subref{fig:x8_y2_nominal} $(8,2)$; and 
  \subref{fig:x9_y2_nominal} $(9,2)$. 
  }
  \label{fig:center_load_y2_nominal}
%\end{figure}
%
\bigskip
%\begin{figure}[tp]
  %%%% C:\doc\robust\topology\load_ccp\eva3\opt_design.m
  \centering
  \begin{subfigure}[b]{0.275\textwidth}
    \centering
    \includegraphics[scale=0.30]{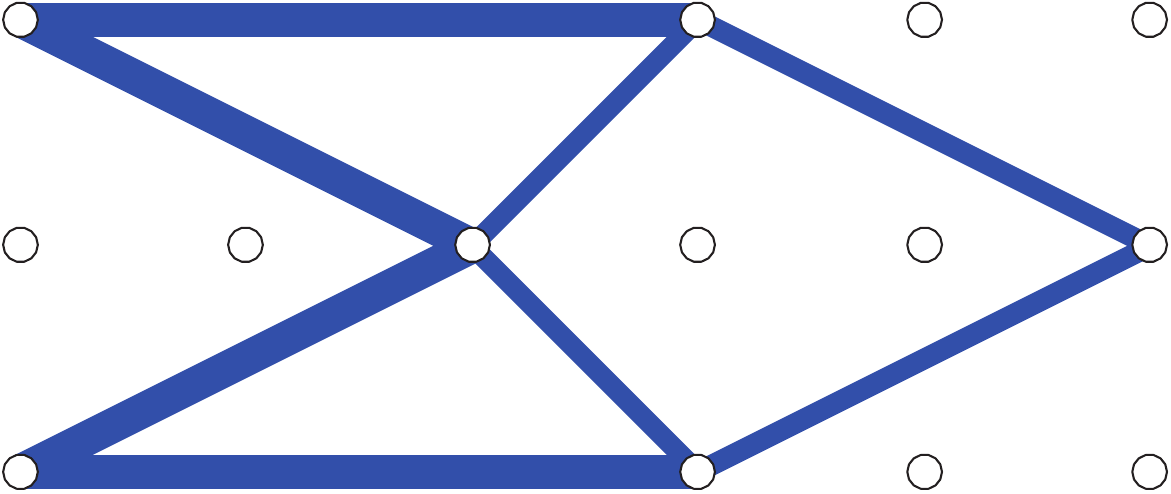}
    \caption{}
    \label{fig:x5_y2_robust}
  \end{subfigure}
  \hfill
  \begin{subfigure}[b]{0.32\textwidth}
    \centering
    \includegraphics[scale=0.36]{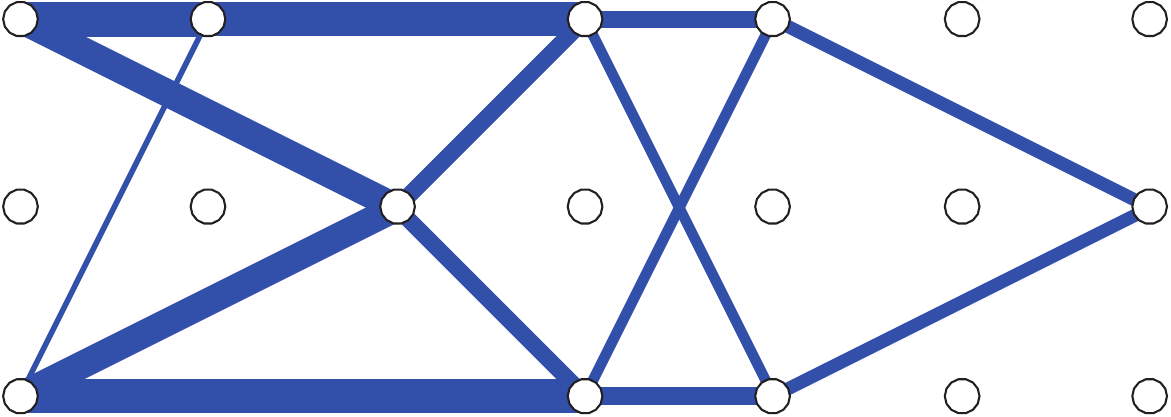}
    \caption{}
    \label{fig:x6_y2_robust}
  \end{subfigure}
  \hfill
  \begin{subfigure}[b]{0.37\textwidth}
    \centering
    \includegraphics[scale=0.42]{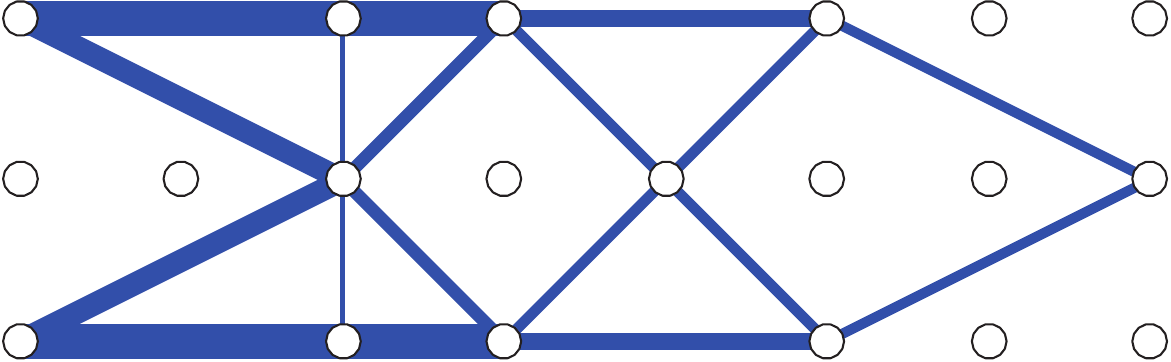}
    \caption{}
    \label{fig:x7_y2_robust}
  \end{subfigure}
  \par\medskip
  \begin{subfigure}[b]{0.42\textwidth}
    \centering
    \includegraphics[scale=0.48]{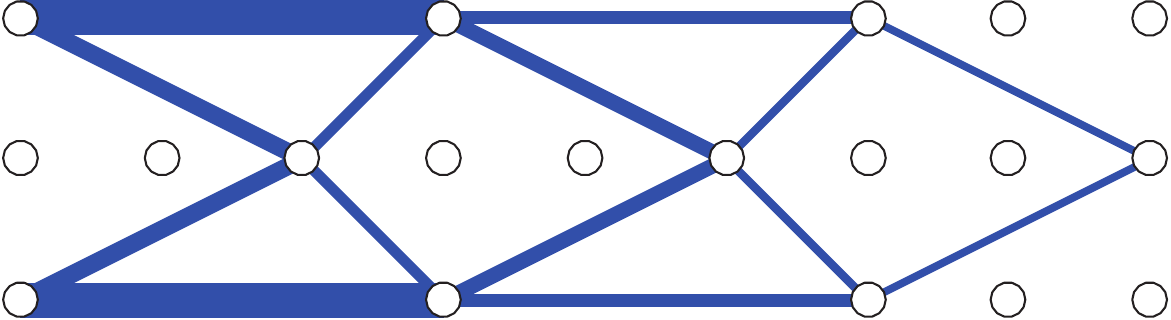}
    \caption{}
    \label{fig:x8_y2_robust}
  \end{subfigure}
  \hfill
  \begin{subfigure}[b]{0.47\textwidth}
    \centering
    \includegraphics[scale=0.54]{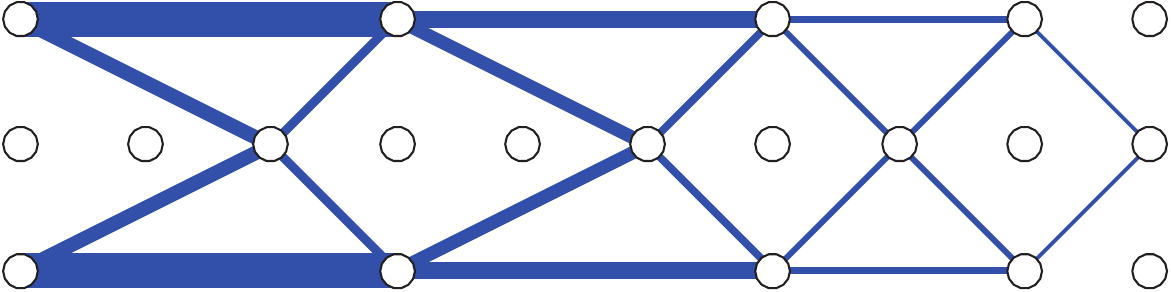}
    \caption{}
    \label{fig:x9_y2_robust}
  \end{subfigure}
  \caption{Example (III)-1. 
  The solutions obtained by the proposed method for the robust 
  optimization under the load uncertainty. 
  \subref{fig:x5_y2_robust} $(N_{X},N_{Y})=(5,2)$; 
  \subref{fig:x6_y2_robust} $(6,2)$; 
  \subref{fig:x7_y2_robust} $(7,2)$; 
  \subref{fig:x8_y2_robust} $(8,2)$; and 
  \subref{fig:x9_y2_robust} $(9,2)$. 
  }
  \label{fig:center_load_y2_robust}
\end{figure}

\begin{figure}[tp]
  %%%% C:\doc\robust\topology\load_ccp\eva3\opt_design.m
  \centering
  \begin{subfigure}[b]{0.275\textwidth}
    \centering
    \includegraphics[scale=0.30]{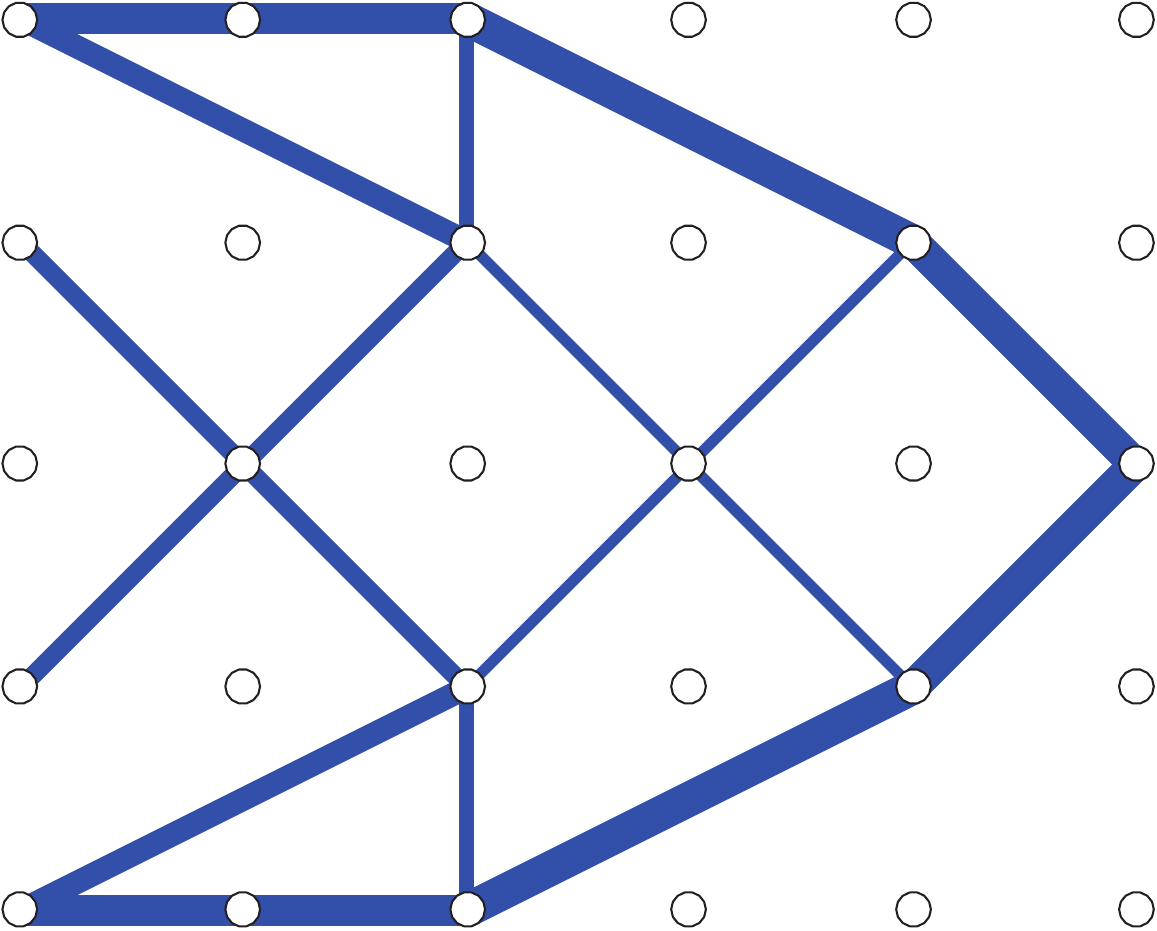}
    \caption{}
    \label{fig:x5_y4_nominal}
  \end{subfigure}
  \hfill
  \begin{subfigure}[b]{0.32\textwidth}
    \centering
    \includegraphics[scale=0.36]{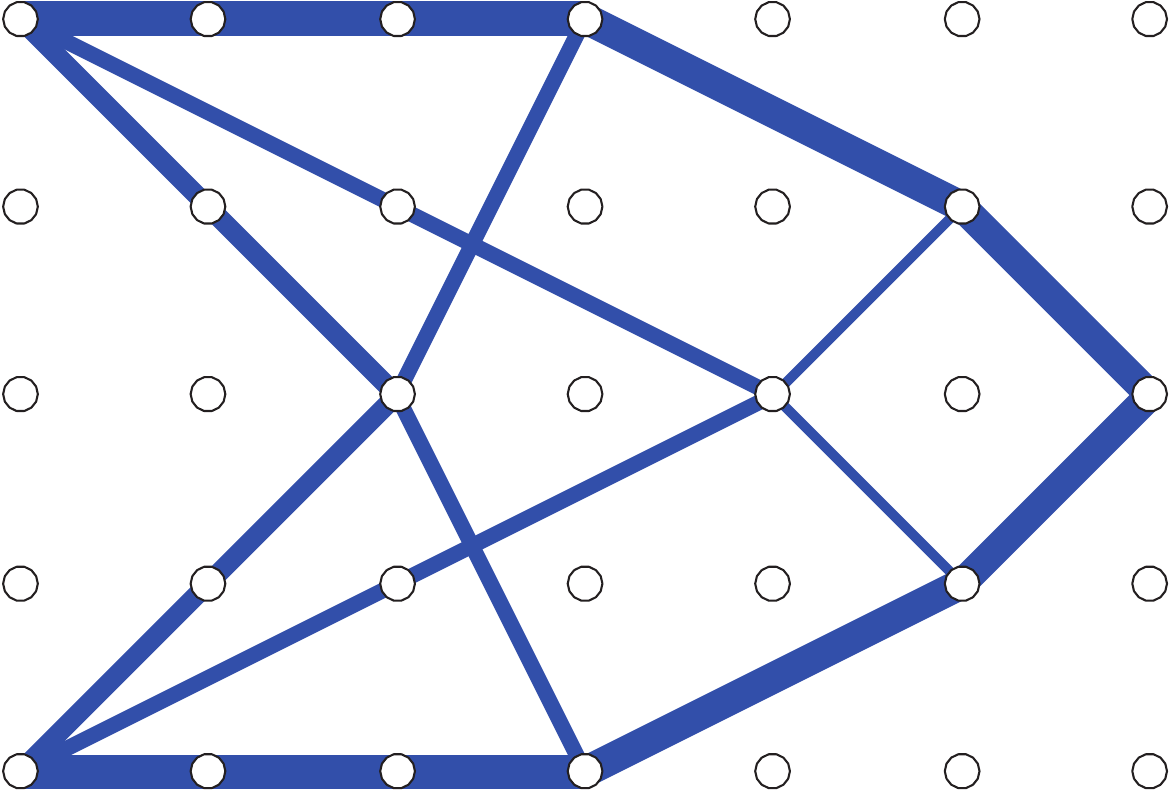}
    \caption{}
    \label{fig:x6_y4_nominal}
  \end{subfigure}
  \hfill
  \begin{subfigure}[b]{0.37\textwidth}
    \centering
    \includegraphics[scale=0.42]{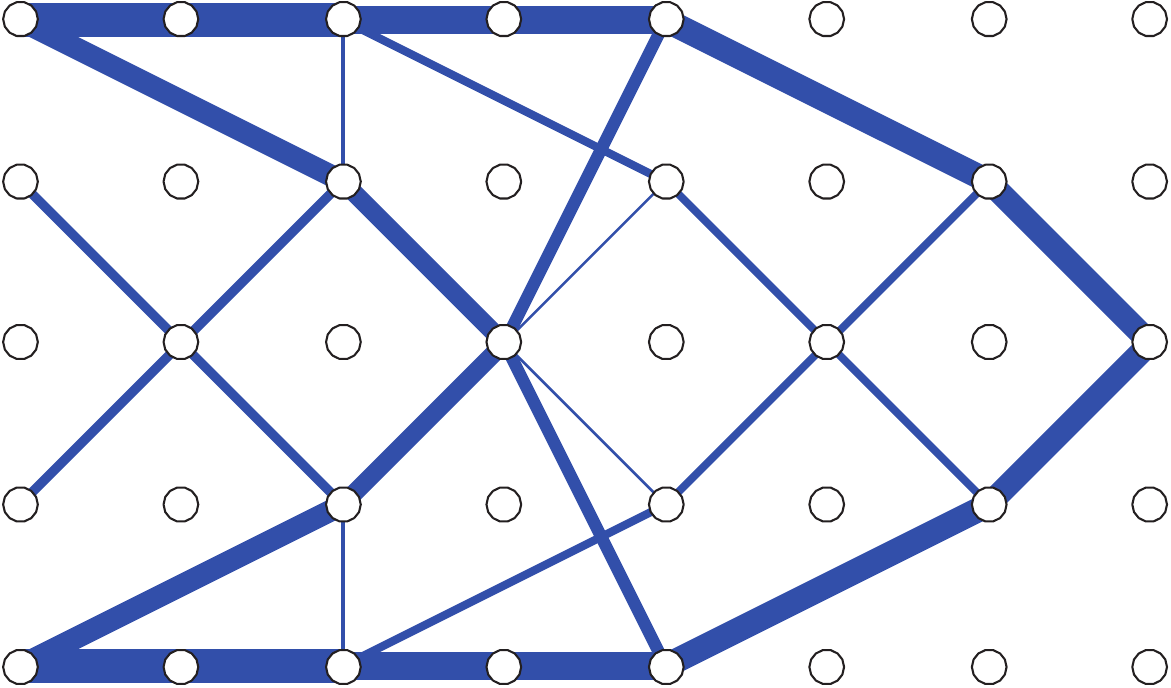}
    \caption{}
    \label{fig:x7_y4_nominal}
  \end{subfigure}
  \par\medskip
  \begin{subfigure}[b]{0.42\textwidth}
    \centering
    \includegraphics[scale=0.48]{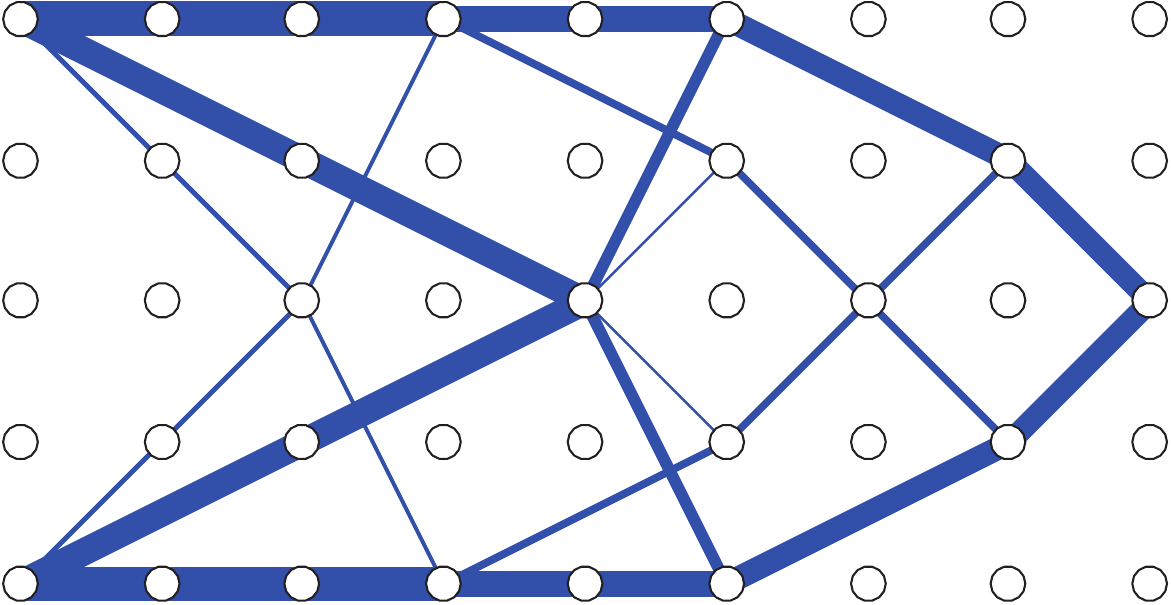}
    \caption{}
    \label{fig:x8_y4_nominal}
  \end{subfigure}
  \hfill
  \begin{subfigure}[b]{0.47\textwidth}
    \centering
    \includegraphics[scale=0.54]{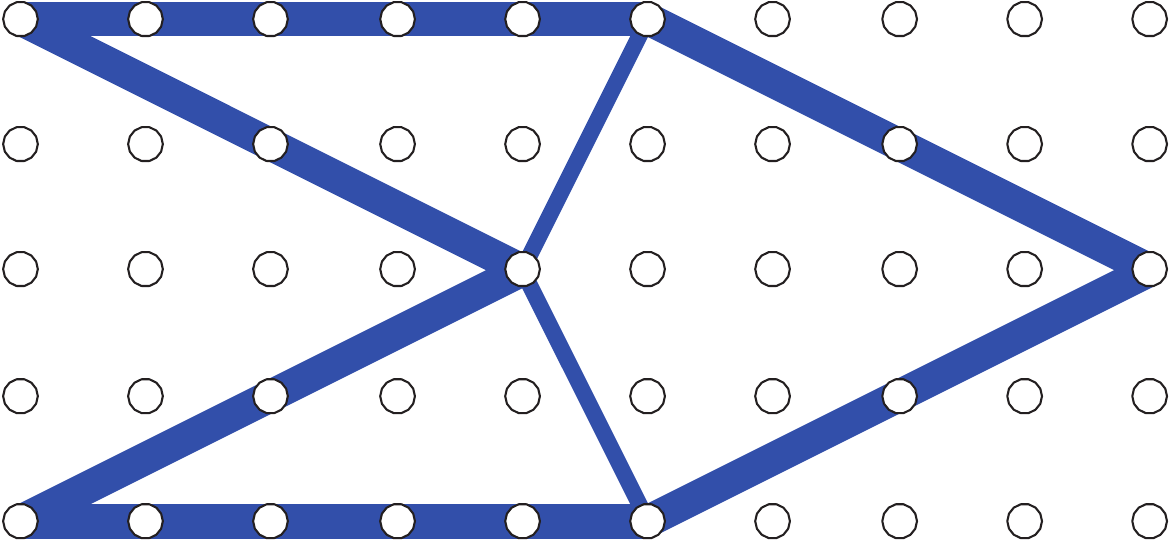}
    \caption{}
    \label{fig:x9_y4_nominal}
  \end{subfigure}
  \caption{Example (III)-2. 
  The optimal solutions of the compliance minimization for the nominal 
  external load. 
  \subref{fig:x5_y4_nominal} $(N_{X},N_{Y})=(5,4)$; 
  \subref{fig:x6_y4_nominal} $(6,4)$; 
  \subref{fig:x7_y4_nominal} $(7,4)$; 
  \subref{fig:x8_y4_nominal} $(8,4)$; and 
  \subref{fig:x9_y4_nominal} $(9,4)$. 
  }
  \label{fig:center_load_y4_nominal}
%\end{figure}
%
\bigskip
%\begin{figure}[tp]
  %%%% C:\doc\robust\topology\load_ccp\eva3\opt_design.m
  \centering
  \begin{subfigure}[b]{0.275\textwidth}
    \centering
    \includegraphics[scale=0.30]{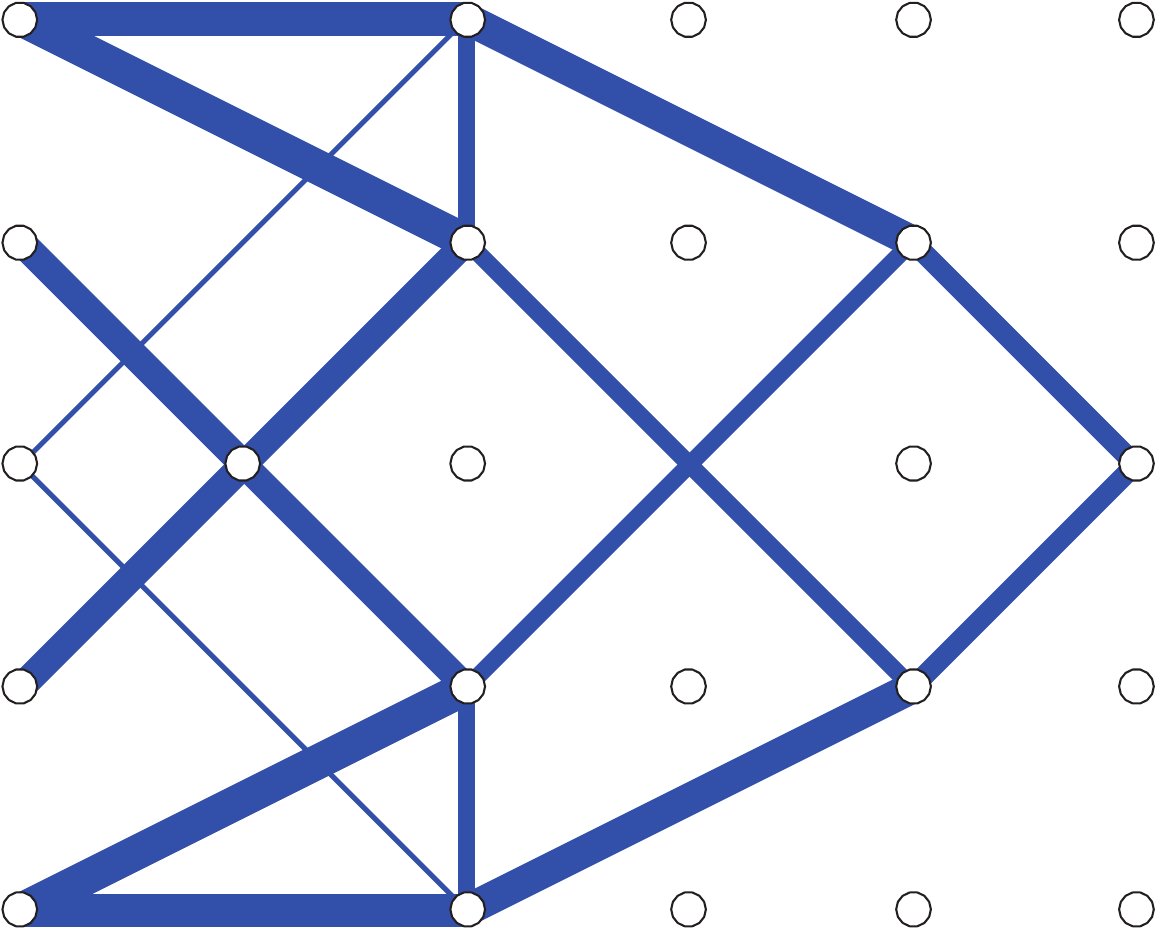}
    \caption{}
    \label{fig:x5_y4_robust}
  \end{subfigure}
  \hfill
  \begin{subfigure}[b]{0.32\textwidth}
    \centering
    \includegraphics[scale=0.36]{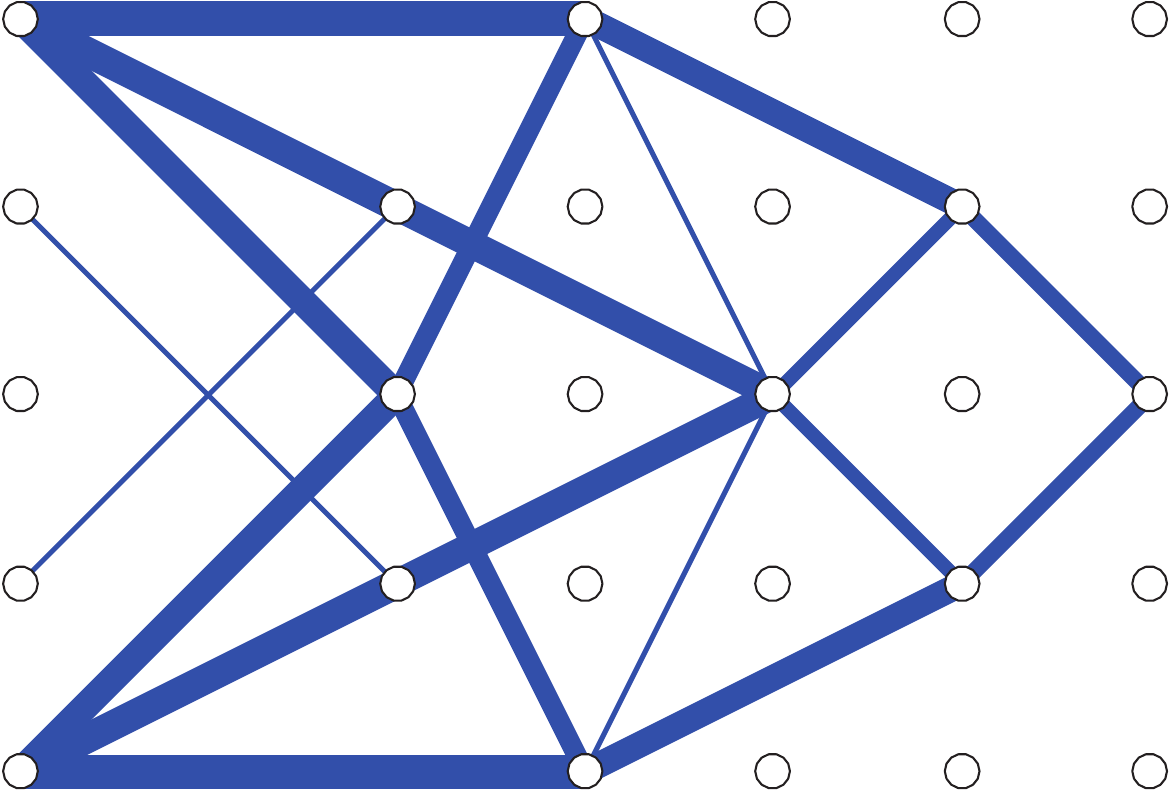}
    \caption{}
    \label{fig:x6_y4_robust}
  \end{subfigure}
  \hfill
  \begin{subfigure}[b]{0.37\textwidth}
    \centering
    \includegraphics[scale=0.42]{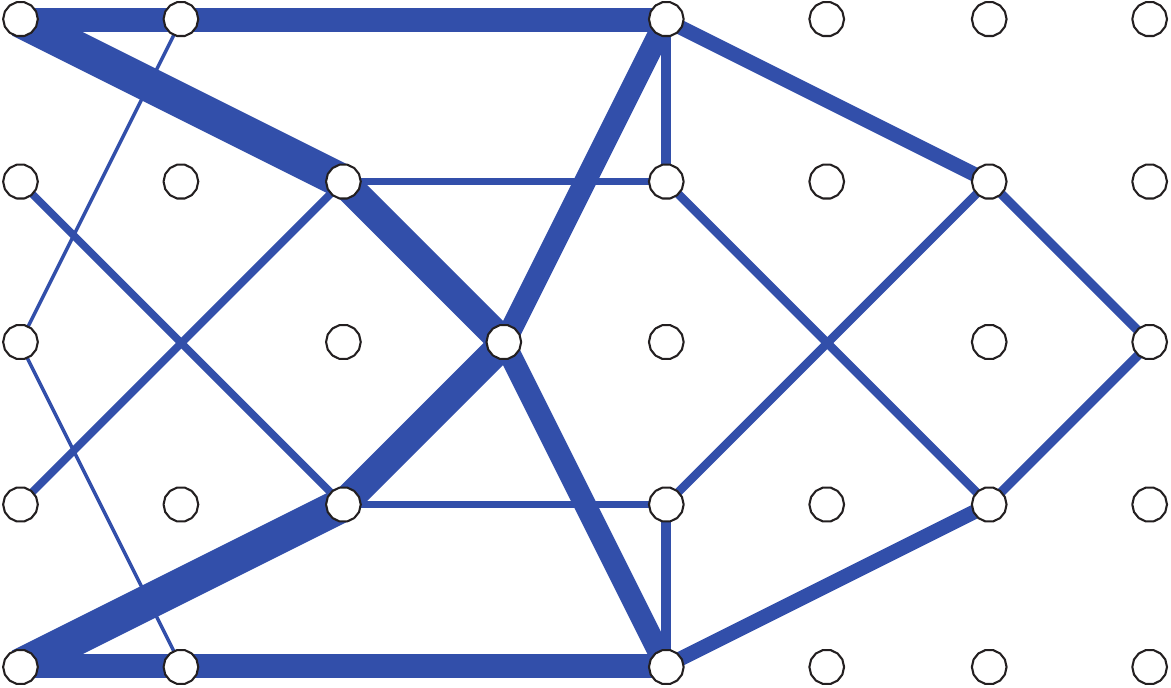}
    \caption{}
    \label{fig:x7_y4_robust}
  \end{subfigure}
  \par\medskip
  \begin{subfigure}[b]{0.42\textwidth}
    \centering
    \includegraphics[scale=0.48]{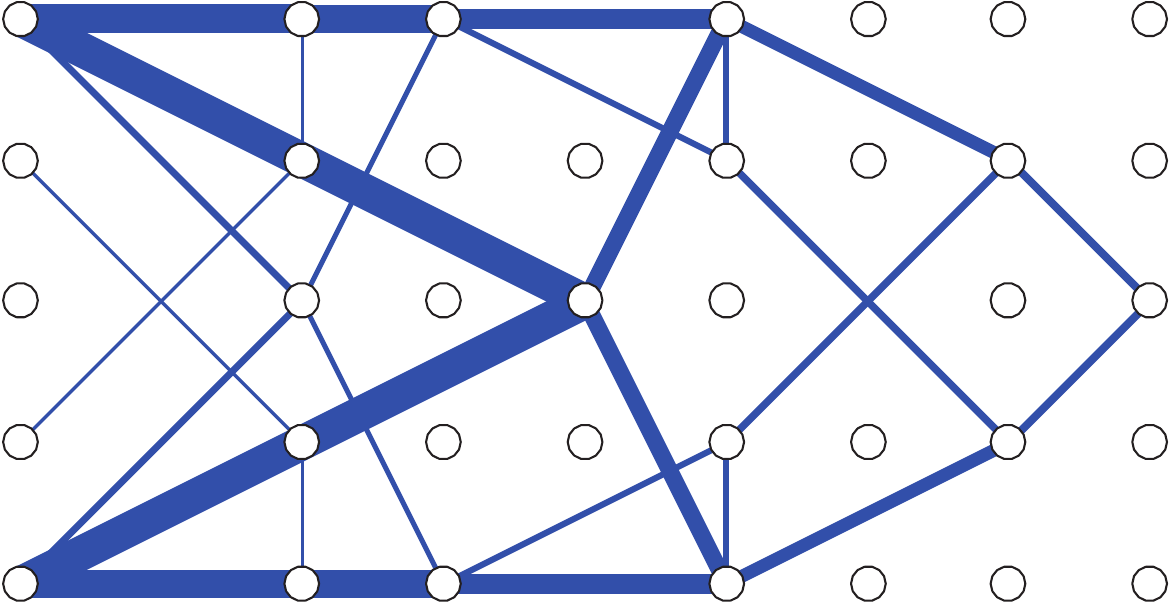}
    \caption{}
    \label{fig:x8_y4_robust}
  \end{subfigure}
  \hfill
  \begin{subfigure}[b]{0.47\textwidth}
    \centering
    \includegraphics[scale=0.54]{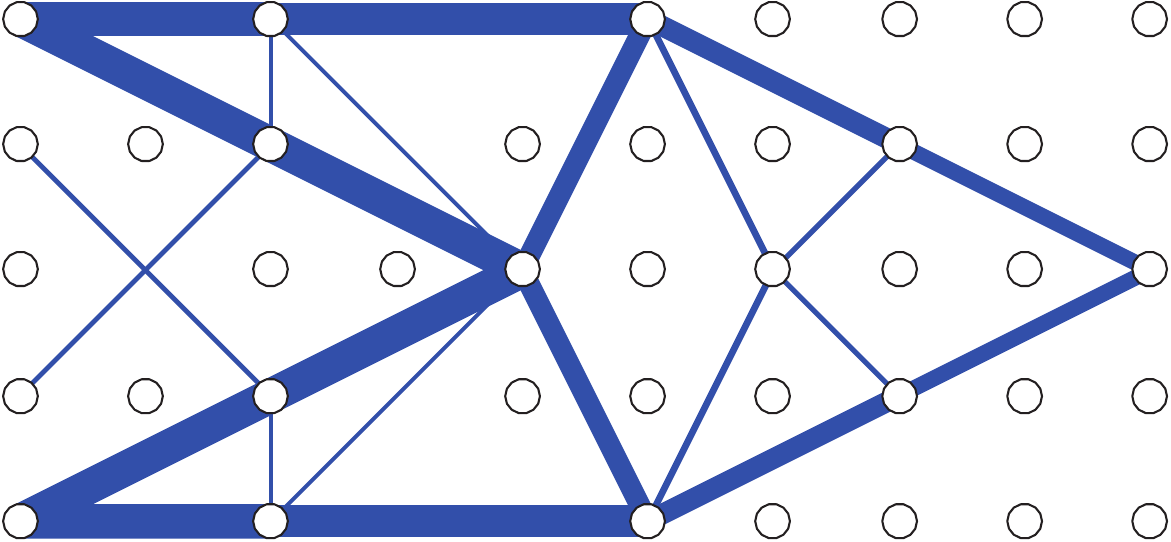}
    \caption{}
    \label{fig:x9_y4_robust}
  \end{subfigure}
  \caption{Example (III)-2. 
  The solutions obtained by the proposed method for the robust 
  optimization under the load uncertainty. 
  \subref{fig:x5_y4_robust} $(5,4)$; 
  \subref{fig:x6_y4_robust} $(6,4)$; 
  \subref{fig:x7_y4_robust} $(7,4)$; 
  \subref{fig:x8_y4_robust} $(8,4)$; and 
  \subref{fig:x9_y4_robust} $(9,4)$. 
  }
  \label{fig:center_load_y4_robust}
\end{figure}

\begin{figure}[tp]
  %%%% C:\doc\robust\topology\load_ccp\eva3\opt_design.m
  \centering
  \begin{subfigure}[b]{0.275\textwidth}
    \centering
    \includegraphics[scale=0.30]{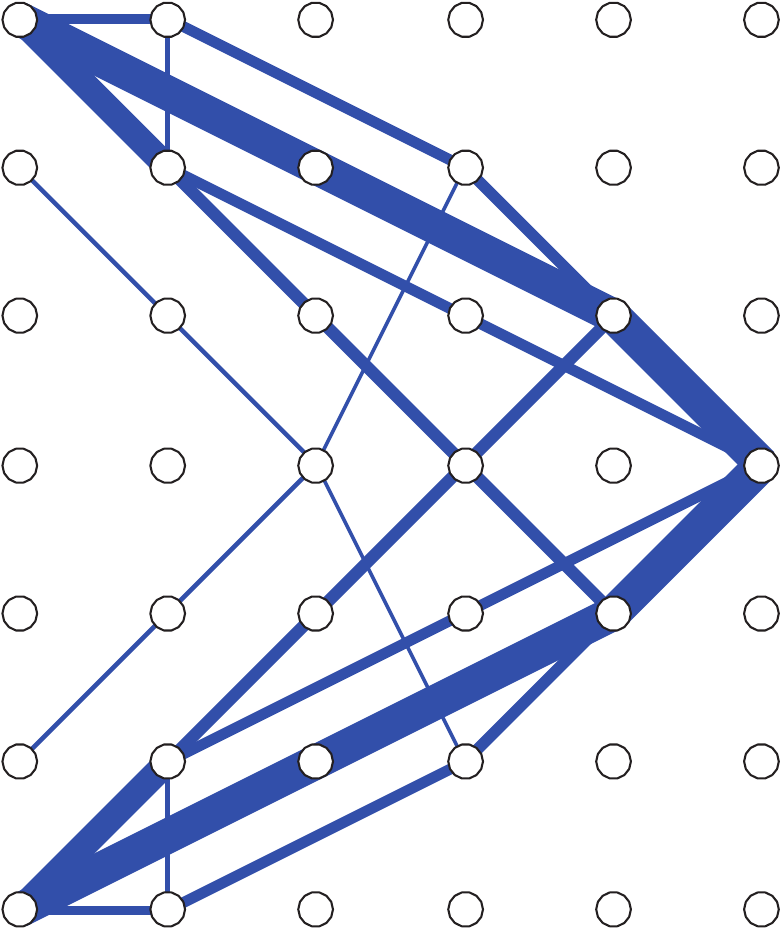}
    \caption{}
    \label{fig:x5_y6_nominal}
  \end{subfigure}
  \hfill
  \begin{subfigure}[b]{0.32\textwidth}
    \centering
    \includegraphics[scale=0.30]{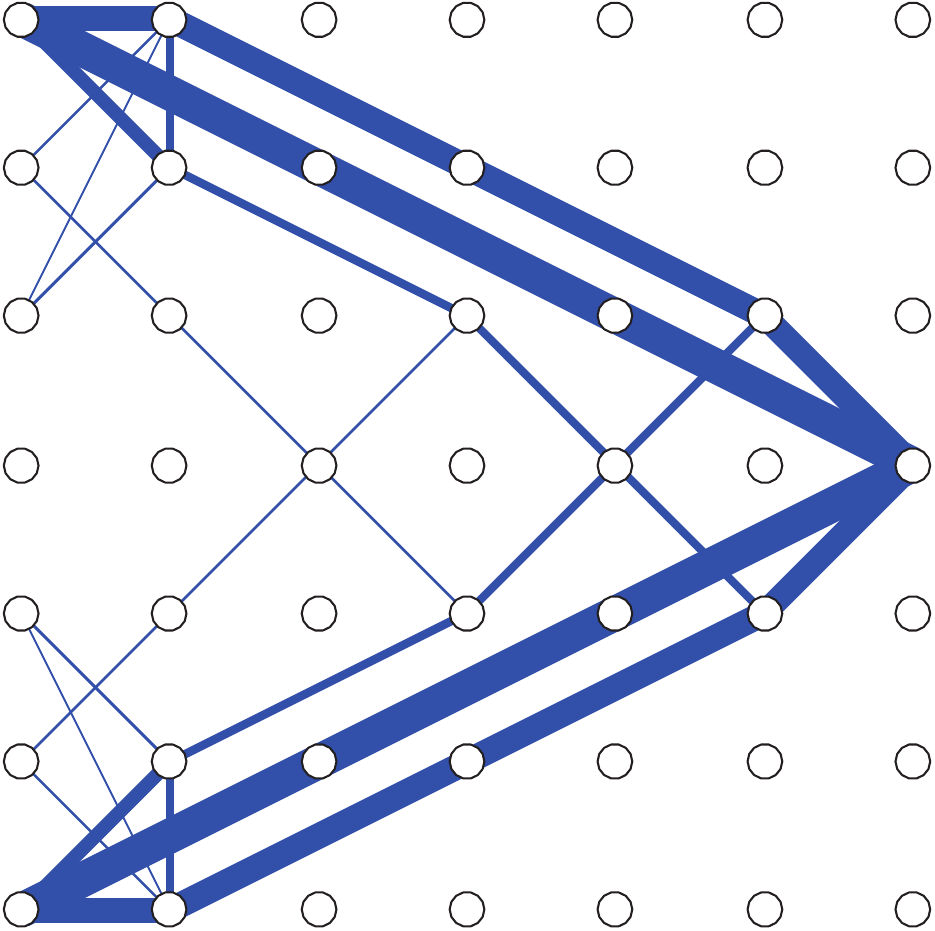}
    \caption{}
    \label{fig:x6_y6_nominal}
  \end{subfigure}
  \hfill
  \begin{subfigure}[b]{0.37\textwidth}
    \centering
    \includegraphics[scale=0.30]{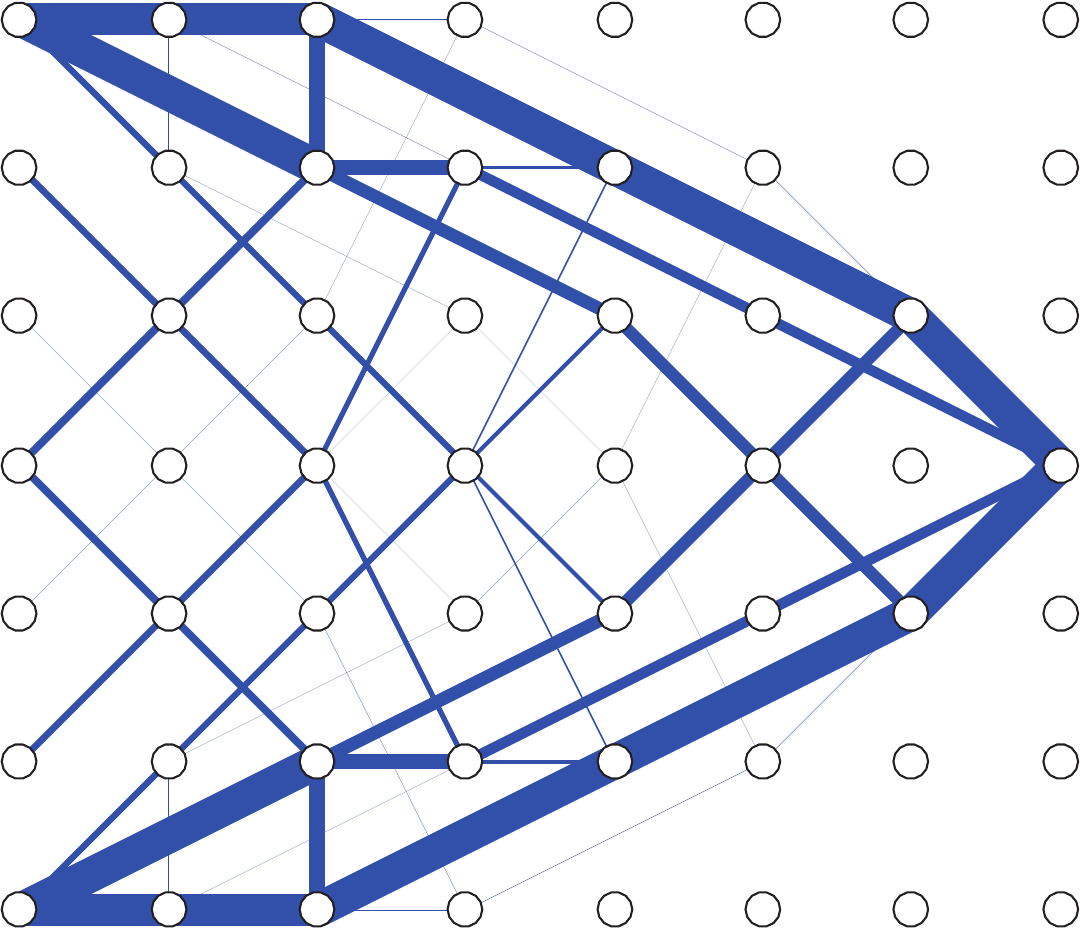}
    \caption{}
    \label{fig:x7_y6_nominal}
  \end{subfigure}
  \par\medskip
  \begin{subfigure}[b]{0.42\textwidth}
    \centering
    \includegraphics[scale=0.30]{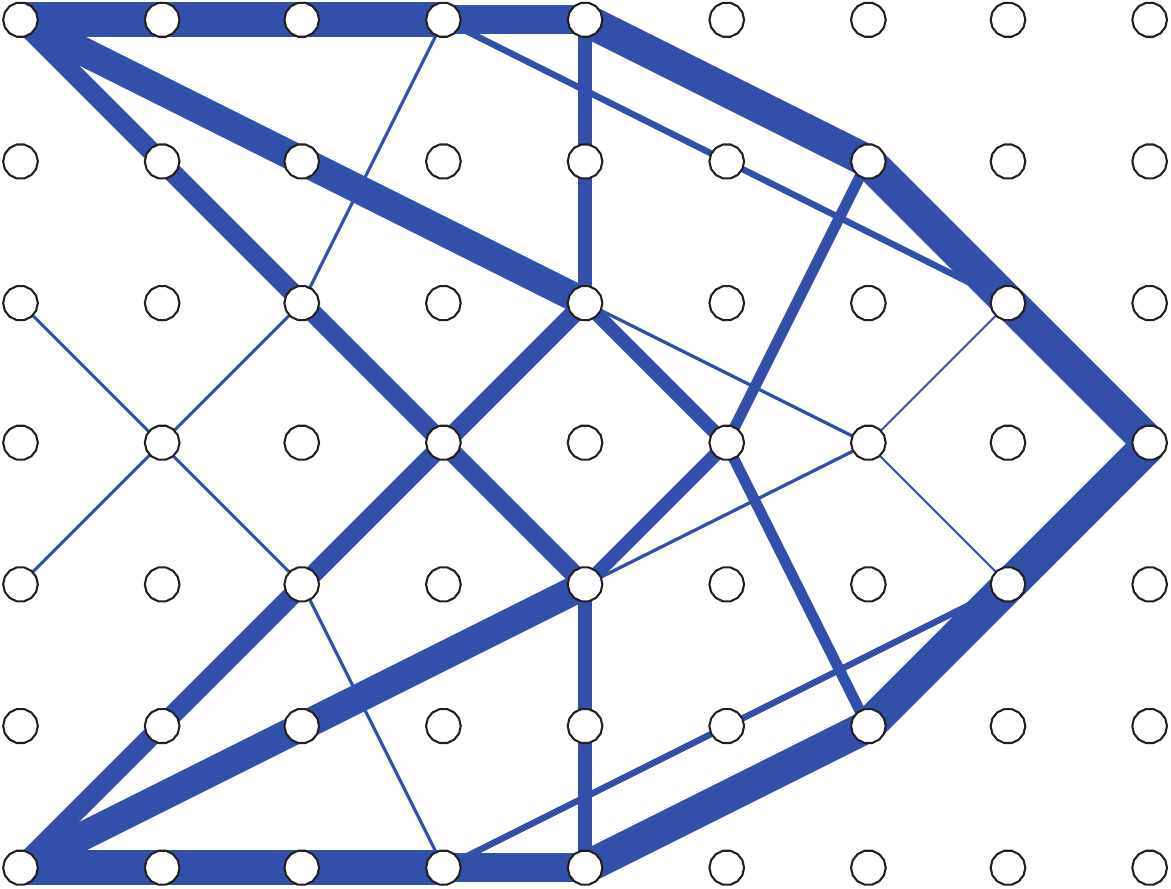}
    \caption{}
    \label{fig:x8_y6_nominal}
  \end{subfigure}
  \hfill
  \begin{subfigure}[b]{0.47\textwidth}
    \centering
    \includegraphics[scale=0.34]{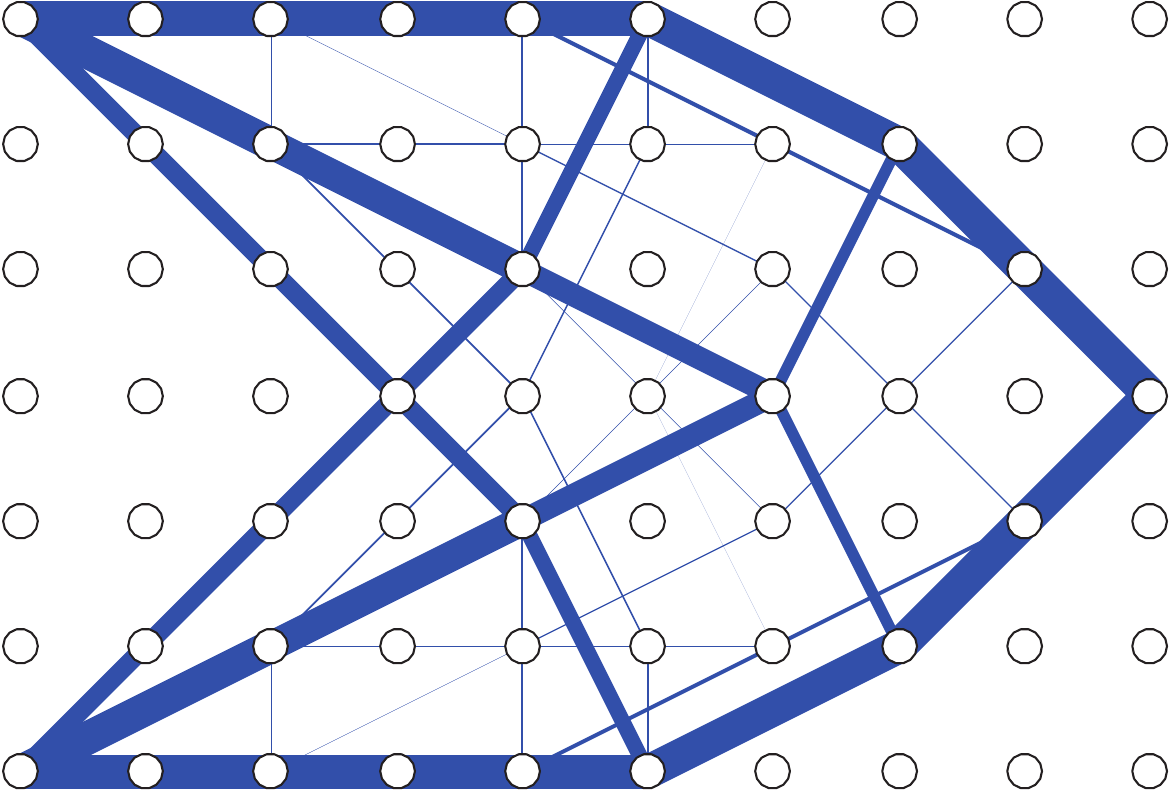}
    \caption{}
    \label{fig:x9_y6_nominal}
  \end{subfigure}
  \caption{Example (III)-3. 
  The optimal solutions of the compliance minimization for the nominal 
  external load. 
  \subref{fig:x5_y6_nominal} $(5,6)$; 
  \subref{fig:x6_y6_nominal} $(6,6)$; 
  \subref{fig:x7_y6_nominal} $(7,6)$; 
  \subref{fig:x8_y6_nominal} $(8,6)$; and
  \subref{fig:x9_y6_nominal} $(9,6)$. 
  }
  \label{fig:center_load_y6_nominal}
%\end{figure}
%
\bigskip
%\begin{figure}[tp]
  %%%% C:\doc\robust\topology\load_ccp\eva3\opt_design.m
  \centering
  \begin{subfigure}[b]{0.275\textwidth}
    \centering
    \includegraphics[scale=0.30]{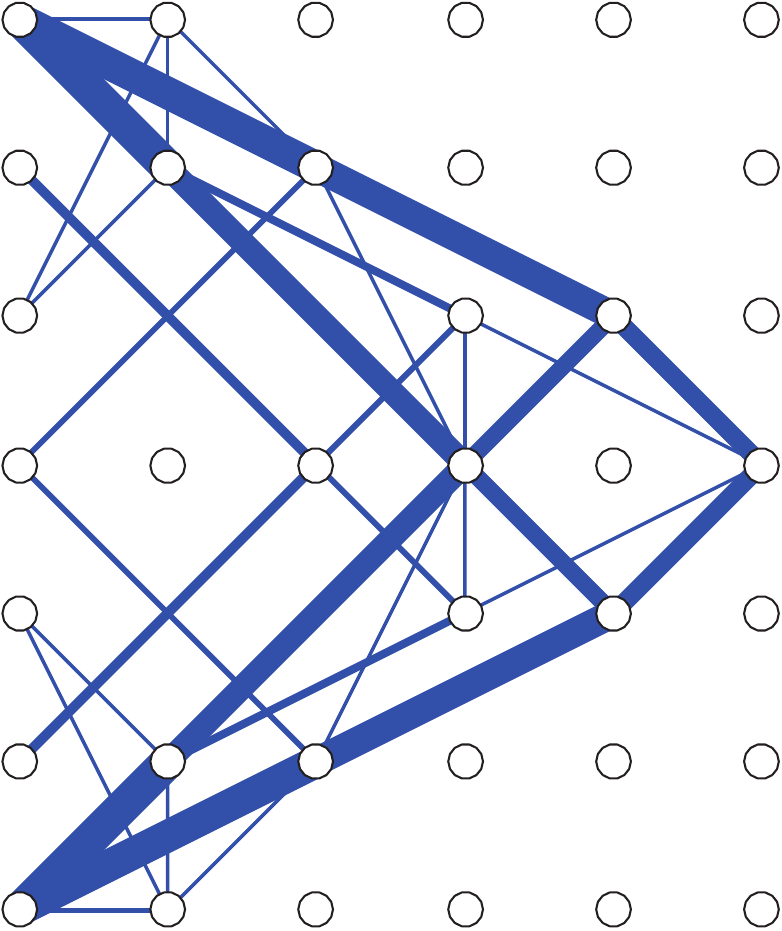}
    \caption{}
    \label{fig:x5_y6_robust}
  \end{subfigure}
  \hfill
  \begin{subfigure}[b]{0.32\textwidth}
    \centering
    \includegraphics[scale=0.30]{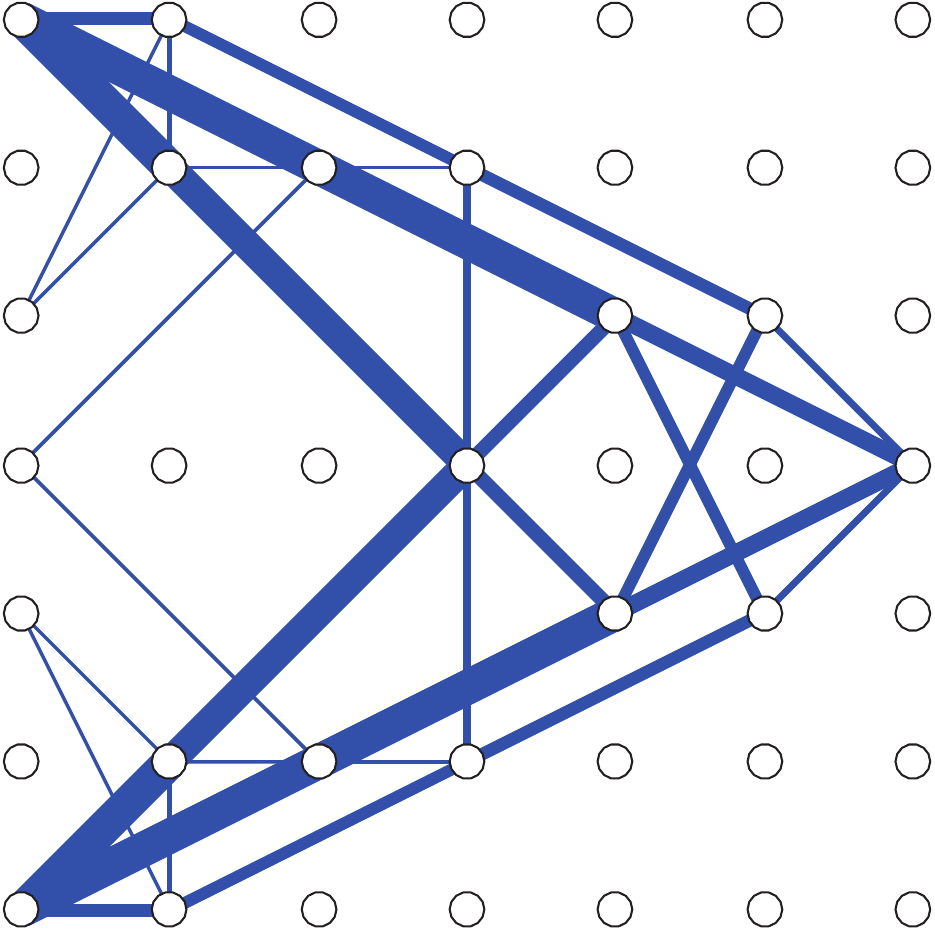}
    \caption{}
    \label{fig:x6_y6_robust}
  \end{subfigure}
  \hfill
  \begin{subfigure}[b]{0.37\textwidth}
    \centering
    \includegraphics[scale=0.30]{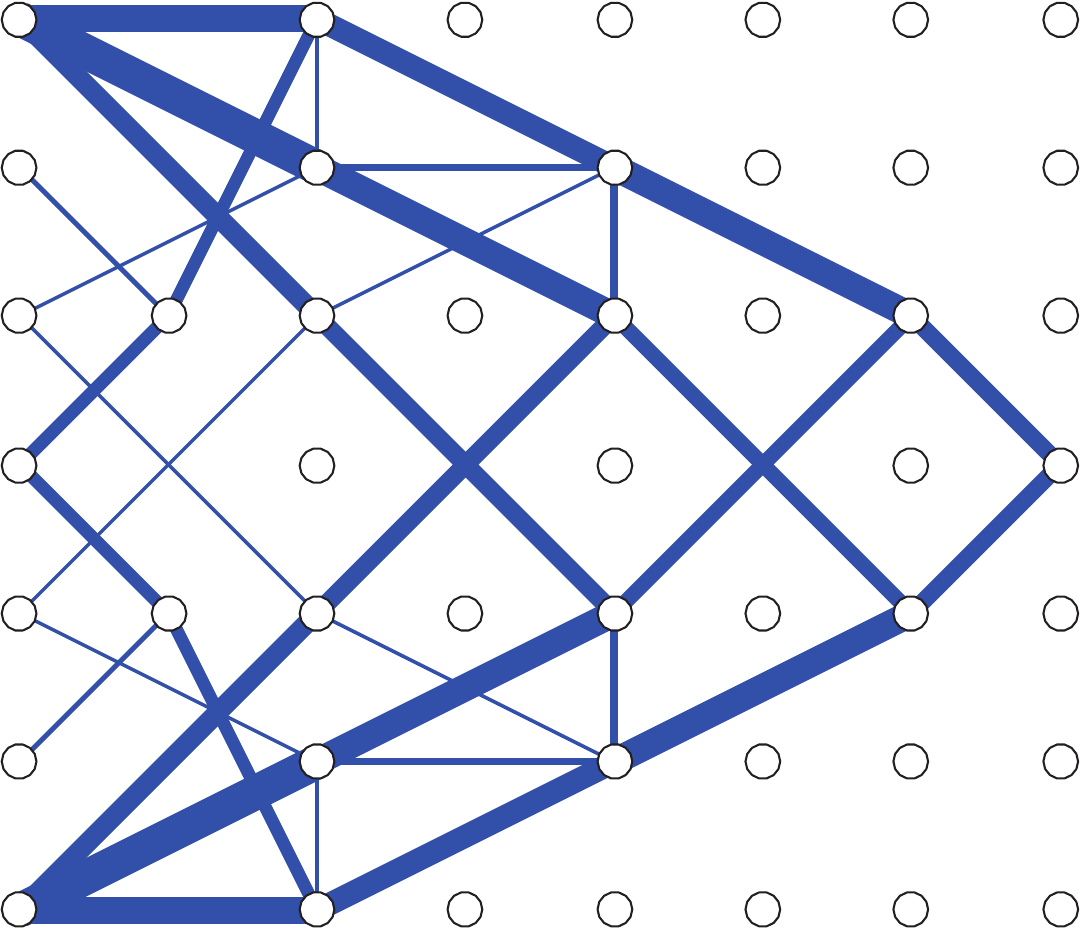}
    \caption{}
    \label{fig:x7_y6_robust}
  \end{subfigure}
  \par\medskip
  \begin{subfigure}[b]{0.42\textwidth}
    \centering
    \includegraphics[scale=0.30]{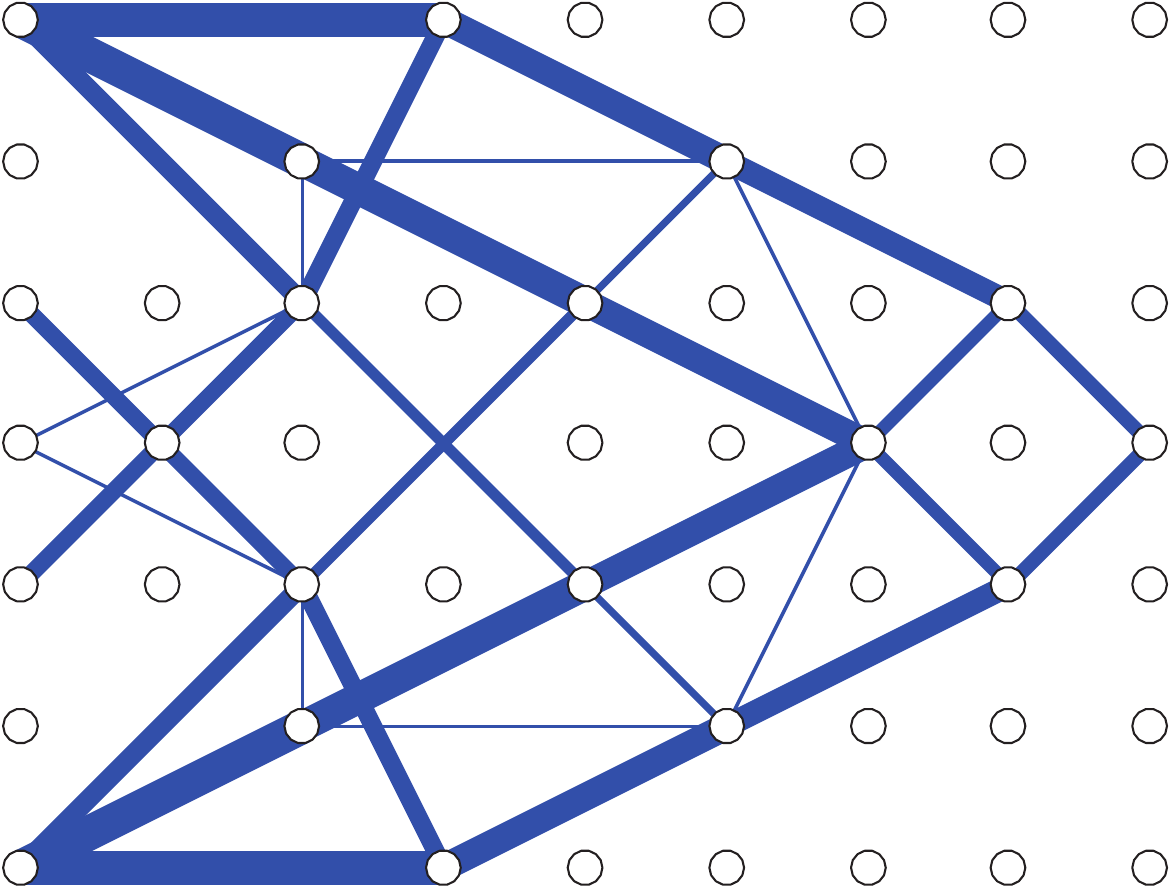}
    \caption{}
    \label{fig:x8_y6_robust}
  \end{subfigure}
  \hfill
  \begin{subfigure}[b]{0.47\textwidth}
    \centering
    \includegraphics[scale=0.34]{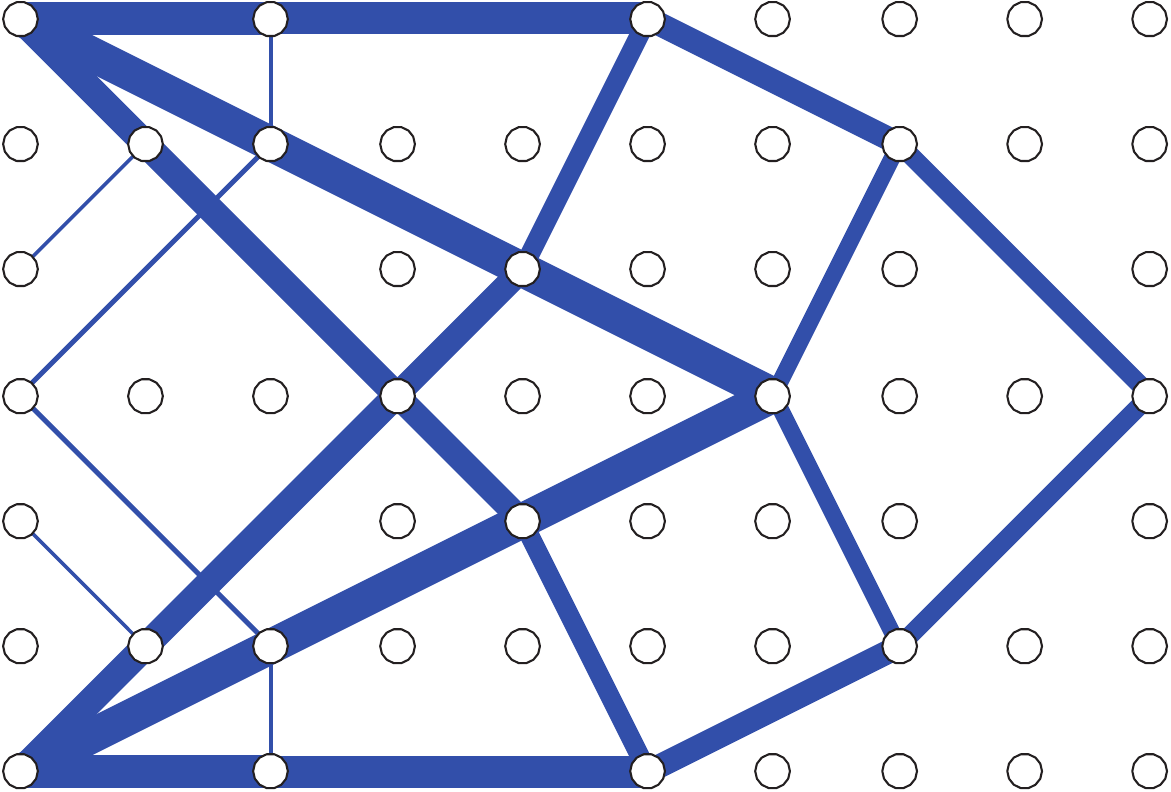}
    \caption{}
    \label{fig:x9_y6_robust}
  \end{subfigure}
  \caption{Example (III)-3. 
  The solutions obtained by the proposed method for the robust 
  optimization under the load uncertainty. 
  \subref{fig:x5_y6_robust} $(5,6)$; 
  \subref{fig:x6_y6_robust} $(6,6)$; 
  \subref{fig:x7_y6_robust} $(7,6)$; 
  \subref{fig:x8_y6_robust} $(8,6)$; and
  \subref{fig:x9_y6_robust} $(9,6)$. 
  }
  \label{fig:center_load_y6_robust}
\end{figure}

\newcommand{\FourTwoM}{\ensuremath{84}}
\newcommand{\FourTwoD}{\ensuremath{24}}
\newcommand{\FourTwoIter}{\ensuremath{16}}
\newcommand{\FourTwoTime}{\ensuremath{43.5}}
\newcommand{\FourTwoObj}{\ensuremath{4438.477}}
\newcommand{\FourTwoTilde}{\ensuremath{3673.056}}
\newcommand{\FourTwoNominal}{\ensuremath{3515.625}}
\newcommand{\FourTwoRatio}{\ensuremath{1.045}}
\newcommand{\FourTwoSpecified}{\ensuremath{6565.430}}
\newcommand{\FourTwoVolume}{\ensuremath{3.2 \times 10^{6}}}

\newcommand{\FourFourM}{\ensuremath{188}}
\newcommand{\FourFourD}{\ensuremath{40}}
\newcommand{\FourFourIter}{\ensuremath{20}}
\newcommand{\FourFourTime}{\ensuremath{52.4}}
\newcommand{\FourFourObj}{\ensuremath{1196.609}}
\newcommand{\FourFourTilde}{\ensuremath{1118.034}}
\newcommand{\FourFourNominal}{\ensuremath{810.826}}
\newcommand{\FourFourRatio}{\ensuremath{1.379}}
\newcommand{\FourFourSpecified}{\ensuremath{1628.906}}
\newcommand{\FourFourVolume}{\ensuremath{6.4 \times 10^{6}}}

\newcommand{\FourSixM}{\ensuremath{292}}
\newcommand{\FourSixD}{\ensuremath{56}}
\newcommand{\FourSixIter}{\ensuremath{15}}
\newcommand{\FourSixTime}{\ensuremath{98.6}}
\newcommand{\FourSixObj}{\ensuremath{666.380}}
\newcommand{\FourSixTilde}{\ensuremath{562.722}}
\newcommand{\FourSixNominal}{\ensuremath{417.598}}
\newcommand{\FourSixRatio}{\ensuremath{1.348}}
\newcommand{\FourSixSpecified}{\ensuremath{891.971}}
\newcommand{\FourSixVolume}{\ensuremath{9.6 \times 10^{6}}}

\newcommand{\FiveTwoM}{\ensuremath{108}}
\newcommand{\FiveTwoD}{\ensuremath{30}}
\newcommand{\FiveTwoIter}{\ensuremath{18}}
\newcommand{\FiveTwoTime}{\ensuremath{25.0}}
\newcommand{\FiveTwoObj}{\ensuremath{7221.094}}
\newcommand{\FiveTwoTilde}{\ensuremath{5708.559}}
\newcommand{\FiveTwoNominal}{\ensuremath{5512.500}}
\newcommand{\FiveTwoRatio}{\ensuremath{1.036}}
\newcommand{\FiveTwoSpecified}{\ensuremath{7924.219}}
\newcommand{\FiveTwoVolume}{\ensuremath{4.0 \times 10^{6}}}

\newcommand{\FiveFourM}{\ensuremath{240}}
\newcommand{\FiveFourD}{\ensuremath{50}}
\newcommand{\FiveFourIter}{\ensuremath{16}}
\newcommand{\FiveFourTime}{\ensuremath{62.1}}
\newcommand{\FiveFourObj}{\ensuremath{2514.685}}
\newcommand{\FiveFourTilde}{\ensuremath{1469.158}}
\newcommand{\FiveFourNominal}{\ensuremath{1304.012}}
\newcommand{\FiveFourRatio}{\ensuremath{1.127}}
\newcommand{\FiveFourSpecified}{\ensuremath{2714.082}}
\newcommand{\FiveFourVolume}{\ensuremath{8.0 \times 10^{6}}}

\newcommand{\FiveSixM}{\ensuremath{372}}
\newcommand{\FiveSixD}{\ensuremath{70}}
\newcommand{\FiveSixIter}{\ensuremath{26}}
\newcommand{\FiveSixTime}{\ensuremath{312.8}}
\newcommand{\FiveSixObj}{\ensuremath{1420.620}}
\newcommand{\FiveSixTilde}{\ensuremath{828.309}}
\newcommand{\FiveSixNominal}{\ensuremath{575.268}}
\newcommand{\FiveSixRatio}{\ensuremath{1.440}}
\newcommand{\FiveSixSpecified}{\ensuremath{1698.426}}
\newcommand{\FiveSixVolume}{\ensuremath{12.0 \times 10^{6}}}

\newcommand{\SixTwoM}{\ensuremath{132}}
\newcommand{\SixTwoD}{\ensuremath{36}}
\newcommand{\SixTwoIter}{\ensuremath{43}}
\newcommand{\SixTwoTime}{\ensuremath{70.4}}
\newcommand{\SixTwoObj}{\ensuremath{13698.325}}
\newcommand{\SixTwoTilde}{\ensuremath{9514.907}}
\newcommand{\SixTwoNominal}{\ensuremath{8760.417}}
\newcommand{\SixTwoRatio}{\ensuremath{1.086}}
\newcommand{\SixTwoSpecified}{\ensuremath{24434.245}}
\newcommand{\SixTwoVolume}{\ensuremath{4.8 \times 10^{6}}}

\newcommand{\SixFourM}{\ensuremath{292}}
\newcommand{\SixFourD}{\ensuremath{60}}
\newcommand{\SixFourIter}{\ensuremath{19}}
\newcommand{\SixFourTime}{\ensuremath{113.9}}
\newcommand{\SixFourObj}{\ensuremath{4063.725}}
\newcommand{\SixFourTilde}{\ensuremath{2160.081}}
\newcommand{\SixFourNominal}{\ensuremath{1814.815}}
\newcommand{\SixFourRatio}{\ensuremath{1.190}}
\newcommand{\SixFourSpecified}{\ensuremath{4200.883}}
\newcommand{\SixFourVolume}{\ensuremath{9.6 \times 10^{6}}}

\newcommand{\SixSixM}{\ensuremath{452}}
\newcommand{\SixSixD}{\ensuremath{84}}
\newcommand{\SixSixIter}{\ensuremath{19}}
\newcommand{\SixSixTime}{\ensuremath{208.0}}
\newcommand{\SixSixObj}{\ensuremath{2168.363}}
\newcommand{\SixSixTilde}{\ensuremath{1197.005}}
\newcommand{\SixSixNominal}{\ensuremath{811.665}}
\newcommand{\SixSixRatio}{\ensuremath{1.475}}
\newcommand{\SixSixSpecified}{\ensuremath{2351.388}}
\newcommand{\SixSixVolume}{\ensuremath{14.4 \times 10^{6}}}

\newcommand{\SevenTwoM}{\ensuremath{156}}
\newcommand{\SevenTwoD}{\ensuremath{42}}
\newcommand{\SevenTwoIter}{\ensuremath{42}}
\newcommand{\SevenTwoTime}{\ensuremath{87.7}}
\newcommand{\SevenTwoObj}{\ensuremath{19198.058}}
\newcommand{\SevenTwoTilde}{\ensuremath{13108.577}}
\newcommand{\SevenTwoNominal}{\ensuremath{12223.214}}
\newcommand{\SevenTwoRatio}{\ensuremath{1.072}}
\newcommand{\SevenTwoSpecified}{\ensuremath{33964.844}}
\newcommand{\SevenTwoVolume}{\ensuremath{5.6 \times 10^{6}}}

\newcommand{\SevenFourM}{\ensuremath{344}}
\newcommand{\SevenFourD}{\ensuremath{70}}
\newcommand{\SevenFourIter}{\ensuremath{40}}
\newcommand{\SevenFourTime}{\ensuremath{320.6}}
\newcommand{\SevenFourObj}{\ensuremath{6563.146}}
\newcommand{\SevenFourTilde}{\ensuremath{3074.707}}
\newcommand{\SevenFourNominal}{\ensuremath{2484.871}}
\newcommand{\SevenFourRatio}{\ensuremath{1.237}}
\newcommand{\SevenFourSpecified}{\ensuremath{6939.107}}
\newcommand{\SevenFourVolume}{\ensuremath{11.2 \times 10^{6}}}

\newcommand{\SevenSixM}{\ensuremath{532}}
\newcommand{\SevenSixD}{\ensuremath{98}}
\newcommand{\SevenSixIter}{\ensuremath{15}}
\newcommand{\SevenSixTime}{\ensuremath{226.7}}
\newcommand{\SevenSixObj}{\ensuremath{2811.740}}
\newcommand{\SevenSixTilde}{\ensuremath{1420.604}}
\newcommand{\SevenSixNominal}{\ensuremath{1123.393}}
\newcommand{\SevenSixRatio}{\ensuremath{1.265}}
\newcommand{\SevenSixSpecified}{\ensuremath{4848.057}}
\newcommand{\SevenSixVolume}{\ensuremath{16.8 \times 10^{6}}}

\newcommand{\EightTwoM}{\ensuremath{180}}
\newcommand{\EightTwoD}{\ensuremath{48}}
\newcommand{\EightTwoIter}{\ensuremath{23}}
\newcommand{\EightTwoTime}{\ensuremath{64.5}}
\newcommand{\EightTwoObj}{\ensuremath{27245.117}}
\newcommand{\EightTwoTilde}{\ensuremath{17498.409}}
\newcommand{\EightTwoNominal}{\ensuremath{16531.250}}
\newcommand{\EightTwoRatio}{\ensuremath{1.059}}
\newcommand{\EightTwoSpecified}{\ensuremath{34548.828}}
\newcommand{\EightTwoVolume}{\ensuremath{6.4 \times 10^{6}}}

\newcommand{\EightFourM}{\ensuremath{396}}
\newcommand{\EightFourD}{\ensuremath{80}}
\newcommand{\EightFourIter}{\ensuremath{37}}
\newcommand{\EightFourTime}{\ensuremath{440.0}}
\newcommand{\EightFourObj}{\ensuremath{9111.073}}
\newcommand{\EightFourTilde}{\ensuremath{3946.944}}
\newcommand{\EightFourNominal}{\ensuremath{3260.031}}
\newcommand{\EightFourRatio}{\ensuremath{1.211}}
\newcommand{\EightFourSpecified}{\ensuremath{10059.252}}
\newcommand{\EightFourVolume}{\ensuremath{12.8 \times 10^{6}}}

\newcommand{\EightSixM}{\ensuremath{612}}
\newcommand{\EightSixD}{\ensuremath{112}}
\newcommand{\EightSixIter}{\ensuremath{32}}
\newcommand{\EightSixTime}{\ensuremath{632.0}}
\newcommand{\EightSixObj}{\ensuremath{4059.438}}
\newcommand{\EightSixTilde}{\ensuremath{1890.087}}
\newcommand{\EightSixNominal}{\ensuremath{1468.478}}
\newcommand{\EightSixRatio}{\ensuremath{1.287}}
\newcommand{\EightSixSpecified}{\ensuremath{6081.969}}
\newcommand{\EightSixVolume}{\ensuremath{19.2 \times 10^{6}}}

\newcommand{\NineTwoM}{\ensuremath{204}}
\newcommand{\NineTwoD}{\ensuremath{54}}
\newcommand{\NineTwoIter}{\ensuremath{28}}
\newcommand{\NineTwoTime}{\ensuremath{100.7}}
\newcommand{\NineTwoObj}{\ensuremath{49880.911}}
\newcommand{\NineTwoTilde}{\ensuremath{23510.654}}
\newcommand{\NineTwoNominal}{\ensuremath{22562.500}}
\newcommand{\NineTwoRatio}{\ensuremath{1.042}}
\newcommand{\NineTwoSpecified}{\ensuremath{86183.594}}
\newcommand{\NineTwoVolume}{\ensuremath{7.2 \times 10^{6}}}

\newcommand{\NineFourM}{\ensuremath{448}}
\newcommand{\NineFourD}{\ensuremath{90}}
\newcommand{\NineFourIter}{\ensuremath{24}}
\newcommand{\NineFourTime}{\ensuremath{286.4}}
\newcommand{\NineFourObj}{\ensuremath{8988.477}}
\newcommand{\NineFourTilde}{\ensuremath{5057.029}}
\newcommand{\NineFourNominal}{\ensuremath{4255.319}}
\newcommand{\NineFourRatio}{\ensuremath{1.188}}
\newcommand{\NineFourSpecified}{\ensuremath{8759.291}}
\newcommand{\NineFourVolume}{\ensuremath{14.4 \times 10^{6}}}

\newcommand{\NineSixM}{\ensuremath{692}}
\newcommand{\NineSixD}{\ensuremath{126}}
\newcommand{\NineSixIter}{\ensuremath{35}}
\newcommand{\NineSixTime}{\ensuremath{863.0}}
\newcommand{\NineSixObj}{\ensuremath{4242.989}}
\newcommand{\NineSixTilde}{\ensuremath{2122.988}}
\newcommand{\NineSixNominal}{\ensuremath{1829.790}}
\newcommand{\NineSixRatio}{\ensuremath{1.160}}
\newcommand{\NineSixSpecified}{\ensuremath{6148.389}}
\newcommand{\NineSixVolume}{\ensuremath{21.6 \times 10^{6}}}

\begin{table}[bp]
  %%%% C:\doc\robust\topology\load_ccp\eva3\opt_design.m
  \centering
  \caption{Characteristics of the problem instances in example (III).}
  \label{tab:ex.II.data}
  \begin{tabular}{lrrr}
    \toprule
    $(N_{X},N_{Y})$ & $m$ & $d$ & $\overline{c}$ $(\mathrm{mm}^{3})$ \\
    \midrule
%    (4,2) 
%    & \FourTwoM & \FourTwoD & \FourTwoNominal & \FourTwoSpecified \\
    (5,2) 
    & \FiveTwoM & \FiveTwoD & \FiveTwoVolume \\
    (6,2) 
    & \SixTwoM & \SixTwoD & \SixTwoVolume \\
    (7,2) 
    & \SevenTwoM & \SevenTwoD & \SevenTwoVolume \\
    (8,2) 
    & \EightTwoM & \EightTwoD & \EightTwoVolume \\
    (9,2) 
    & \NineTwoM & \NineTwoD & \NineTwoVolume \\
    \midrule
%    (4,4) 
%    & \FourFourM & \FourFourD & \FourFourNominal & \FourFourSpecified \\
    (5,4) 
    & \FiveFourM & \FiveFourD & \FiveFourVolume \\
    (6,4) 
    & \SixFourM & \SixFourD & \SixFourVolume \\
    (7,4) 
    & \SevenFourM & \SevenFourD & \SevenFourVolume \\
    (8,4) 
    & \EightFourM & \EightFourD & \EightFourVolume \\
    (9,4) 
    & \NineFourM & \NineFourD & \NineFourVolume \\
    \midrule
%    (4,6) 
%    & \FourSixM & \FourSixD & \FourSixNominal & \FourSixSpecified \\
    (5,6) 
    & \FiveSixM & \FiveSixD & \FiveSixVolume \\
    (6,6) 
    & \SixSixM & \SixSixD & \SixSixVolume \\
    (7,6) 
    & \SevenSixM & \SevenSixD & \SevenSixVolume \\
    (8,6) 
    & \EightSixM & \EightSixD & \EightSixVolume \\
    (9,6) 
    & \NineSixM & \NineSixD & \NineSixVolume \\
    \bottomrule
  \end{tabular}
%\end{table}
%
%\begin{table}[bp]
  %%%% C:\doc\robust\topology\load_ccp\eva3\opt_design.m
  \centering
  \caption{Computational results of example (III).}
  \label{tab:ex.II.result}
  \begin{tabular}{lrrrrr}
    \toprule
    $(N_{X},N_{Y})$ & Obj.\ (J) & {\#}iter. 
    & Time (s) & $\tilde{w}$ (J) & Nom.\ opt.\ (J) \\
    \midrule
%    (4,2)
%    & \FourTwoObj & \FourTwoIter & \FourTwoTime & \FourTwoTilde & \FourTwoRatio \\
    (5,2)
    & \FiveTwoObj & \FiveTwoIter & \FiveTwoTime & \FiveTwoTilde & \FiveTwoNominal \\
    (6,2)
    & \SixTwoObj & \SixTwoIter & \SixTwoTime & \SixTwoTilde & \SixTwoNominal \\
    (7,2)
    & \SevenTwoObj & \SevenTwoIter & \SevenTwoTime & \SevenTwoTilde & \SevenTwoNominal \\
    (8,2)
    & \EightTwoObj & \EightTwoIter & \EightTwoTime & \EightTwoTilde & \EightTwoNominal \\
    (9,2)
    & \NineTwoObj & \NineTwoIter & \NineTwoTime & \NineTwoTilde & \NineTwoNominal \\
    \midrule
%    (4,4)
%    & \FourFourObj & \FourFourIter & \FourFourTime & \FourFourTilde & \FourFourNominal \\
    (5,4)
    & \FiveFourObj & \FiveFourIter & \FiveFourTime & \FiveFourTilde & \FiveFourNominal \\
    (6,4)
    & \SixFourObj & \SixFourIter & \SixFourTime & \SixFourTilde & \SixFourNominal \\
    (7,4)
    & \SevenFourObj & \SevenFourIter & \SevenFourTime & \SevenFourTilde & \SevenFourNominal \\
    (8,4)
    & \EightFourObj & \EightFourIter & \EightFourTime & \EightFourTilde & \EightFourNominal \\
    (9,4)
    & \NineFourObj & \NineFourIter & \NineFourTime & \NineFourTilde & \NineFourNominal \\
    \midrule
%    (4,6)
%    & \FourSixObj & \FourSixIter & \FourSixTime & \FourSixTilde & \FourSixNominal \\
    (5,6)
    & \FiveSixObj & \FiveSixIter & \FiveSixTime & \FiveSixTilde & \FiveSixNominal \\
    (6,6)
    & \SixSixObj & \SixSixIter & \SixSixTime & \SixSixTilde & \SixSixNominal \\
    (7,6)
    & \SevenSixObj & \SevenSixIter & \SevenSixTime & \SevenSixTilde & \SevenSixNominal \\
    (8,6)
    & \EightSixObj & \EightSixIter & \EightSixTime & \EightSixTilde & \EightSixNominal \\
    (9,6)
    & \NineSixObj & \NineSixIter & \NineSixTime & \NineSixTilde & \NineSixNominal \\
    \bottomrule
  \end{tabular}
\end{table}

Consider the problem setting shown in \reffig{fig:gs6x2}. 
Ground structures are generated in the manner explained in 
section~\ref{sec:ex.eva4}. 
The maximum length of the members in a ground structure is 
$3\,\mathrm{m}$. 
The uncertainty model of the external load is defined by using 
\eqref{eq.ex.uncertainty.1} with 
$\tilde{p}_{1}=100\,\mathrm{kN}$ and $\alpha=0.75\tilde{p}_{1}$. 
The lower and upper bounds for the member cross-sectional areas are 
$\underline{x}=50\,\mathrm{mm^{2}}$ and 
$\overline{x}=500\,\mathrm{mm^{2}}$, respectively. 
The upper bound for the structural volume is 
$\overline{c}=4N_{X}N_{Y} \times 10^{5}\,\mathrm{mm}^{3}$. 

For problem instances with $(N_{X},N_{Y})=(5,2)$, $(6,2),\dots,(9,2)$, 
\reffig{fig:center_load_y2_nominal} shows the optimal solutions of the 
compliance minimization for the nominal external load. 
For the robust optimization, the solutions obtained by the proposed 
method are collected in \reffig{fig:center_load_y2_robust}. 
The instance sizes are listed in \reftab{tab:ex.II.data}. 
The computational results are listed in \reftab{tab:ex.II.result}. 
It is observed in \reffig{fig:x5_y2_robust} and \reffig{fig:x8_y2_robust} 
that the intermediate nodes on chains in \reffig{fig:x5_y2_nominal} and 
\reffig{fig:x8_y2_nominal} are removed as a result of robust 
optimization. 
Although the nominal optimal solutions in 
\reffig{fig:x6_y2_nominal}, 
\reffig{fig:x7_y2_nominal}, and \reffig{fig:x9_y2_nominal} have very 
complicated forms, the robust solutions in 
\reffig{fig:x6_y2_robust}, 
\reffig{fig:x7_y2_robust}, and \reffig{fig:x9_y2_robust} are simple and 
practically preferable. 
Thus, it is often that robustness against uncertain loads and the 
minimal cross-sectional area constraints for the existing 
members yield simple truss 
topology. 

The solutions obtained for the instances with 
$(N_{X},N_{Y})=(5,4)$, $(6,4),\dots,(9,4)$ are collected in 
\reffig{fig:center_load_y4_nominal} and 
\reffig{fig:center_load_y4_robust}. 
The robust optimal solution obtained by the proposed method has a form 
similar to the corresponding nominal optimal solution, but many chains 
in the nominal optimal solution are replaced with single members. 

The solutions obtained for the instances with 
$(N_{X},N_{Y})=(5,6)$, $(6,6),\dots,(9,6)$ are collected in 
\reffig{fig:center_load_y6_nominal} and 
\reffig{fig:center_load_y6_robust}. 
It is observed that the nominal optimal solutions in 
\reffig{fig:x7_y6_nominal} and 
\reffig{fig:x9_y6_nominal} have so many thin members. 
In contrast, the robust solutions in \reffig{fig:x7_y6_robust} and 
\reffig{fig:x9_y6_robust} have fewer members. 
The layout of thick members in \reffig{fig:x6_y6_robust} is different 
from that in \reffig{fig:x6_y6_nominal}. 
Also, the layout of thick members in \reffig{fig:x8_y6_robust} is 
different from that in \reffig{fig:x8_y6_nominal}. 

It is observed in \reftab{tab:ex.II.result} that the proposed method 
converged mostly within 40 iterations. 
The computational time for the instance with about $600$ members is 
about $10$ minutes. 
Thus, the proposed method finds a reasonable feasible 
solution with relatively 
small computational cost.

\section{Conclusions}
\label{sec:conclude}

In this paper, we have presented a new formulation and algorithm for 
robust truss topology optimization considering uncertainty in the 
external load. 
Specifically, combinatorial aspects of the problem has been dealt with 
in the framework of the complementarity constraints. 
As for the uncertainty model of the external load, we have supposed that 
uncertain external forces can 
be applied to all the nodes of a truss. 
This model depends on the set of existing nodes, and hence on the set of 
existing members, while the member cross-sectional areas are the design 
variables to be optimized. 
Thus, the robust optimization problem involves design-dependent 
constraints. 
Also, it has been explained that overlapping members should be 
incorporated to a ground structure. 
In the final truss design, however, presence of overlapping members is 
not allowed from a practical point of view. 
In this paper, the set of existing nodes, the selection among 
overlapping members, and the lower bound constraints for the 
cross-sectional areas of the existing members are treated by using the 
complementarity constraints. 

In the conventional truss topology optimization, it is often that an 
optimal solution has a sequence of parallel consecutive members, called 
a chain. 
To stabilize a truss, a chain is replaced with a longer single member. 
Special consideration, such as the local buckling constraints, is needed 
to avoid presence of a too long member converted from a chain. 
In contrast, the solution obtained with the proposed method does not 
have a member which is longer than the maximum member length of the 
ground structure, because a truss design including a chain is infeasible 
for the presented robust optimization problem. 

This paper has presented an SDPCC (semidefinite programming with 
complementarity constraints) formulation of a structural optimization 
problem. 
Then, its DC programming reformulation has been solved with a 
convex-concave procedure. 
This algorithm is a version of MM algorithms and EM algorithms, which 
are widely used in machine learning, image processing, etc. 
It has been shown through the numerical experiments that the proposed 
heuristic can converge to a high-quality solution within relatively 
small computational cost. 
The method can certainly 
handle complementarity constraints other than the ones 
presented in this paper. 
An example is a set of constraints that prohibits the presence of 
mutually crossing members in a truss design.

\paragraph{Acknowledgments}

This work is partially supported by 
JSPS KAKENHI 26420545 and 17K06633.

\end{document}